\overfullrule=0pt
\input amssym
\input xypic
\def\Qxy{{{\Bbb Q}\langle\langle x,y\rangle\rangle}}
\def\Qab{{{\Bbb Q}\langle\langle a,b\rangle\rangle}}
\def\Rab{{R\langle\langle a,b\rangle\rangle}}

\def\ZZab{{\overline{\cal Z}\langle\langle a,b\rangle\rangle}}
\def\RRab{{\bar{R}\langle\langle a,b\rangle\rangle}}
\def\QC{{{\Bbb Q}\langle\langle c_1,c_2,\ldots\rangle\rangle}}

\def\grt{{\frak{grt}}}
\def\krv{{\frak{krv}}}

\def\fu{{\frak{u}}}

\def\ba{{\frak{a}}}
\def\be{{\bf{e}}}
\def\ba{{\bf{a}}}
\def\uk{{\underline{k}}}
\def\Ufu{{{\cal U}{\fu}}}
\def\F{{{\cal F}}}
\def\G{{{\cal G}}}
\def\QQ{{{\Bbb Q}}}
\def\Q{{{\bf Q}}}
\def\M{{\bf M}}
\def\P{{\bf P}}
\def\Zi{{{\Bbb Z}}}
\def\Z{{{\cal Z}}}

\centerline{\bf The Fay relations satisfied by the elliptic associator}
\vskip .4cm
\centerline{Leila Schneps}
\vskip 1cm
\centerline{\it To Moshe Jarden on the occasion of his 80th birthday}
\vskip 1cm
\noindent {\bf Abstract.} We recall the construction by B.~Enriquez 
of the elliptic associator $A_\tau$, a power series in two non-commutative 
variables $a,b$ defined as an iterated integral of the 
Kronecker function, and turn our attention to a family of 
{\it Fay relations} satisfied by $A_\tau$, derived from the original well-known
Fay relation satisfied by the Kronecker function. The Fay relations of
$A_\tau$ were studied by Broedel, Matthes and Schlotterer, and determined up to 
non-explicit correction terms that arise from the necessity of regularizing the 
non-convergent integral.  
In this article, we study a reduced version $\bar{A}_\tau$ of the elliptic 
associator mod $2\pi i$. We recall a different construction of $\bar{A}_\tau$ in three 
steps, due to Matthes, Lochak and the author: first one defines the 
reduced {\it elliptic generating series} $\bar{E}_\tau$ which comes 
from the reduced Drinfel'd associator $\overline{\Phi}_{KZ}$ and whose 
coefficients generate the same ring $\bar{R}$ as those of $\bar{A}_\tau$;
then one defines $\Psi$ to be the automorphism of the free associative ring 
$\RRab$ defined by $\Psi(a)=\bar{E}_\tau$ and $\Psi([a,b])=[a,b]$; 
finally one shows that the reduced elliptic associator $\bar{A}_\tau$ is equal 
to $\Psi\bigl({{ad(b)}\over{e^{ad(b)}-1}}(a)\bigr)$.
Using this construction and mould theory and working with Lie-like versions
of the elliptic generating series and associator, we prove the following 
results: first, a mould satisfies the Fay relations if and only if a
closely related  mould satisfies the well-known ``swap circ-neutrality'' 
relations defining the elliptic Kashiwara-Vergne Lie algebra $\krv_{ell}$,
second, the reduced elliptic generating series satisfies a 
family of Fay relations with extremely simple correction terms coming directly
from those of the Drinfel'd associator, and third, the correction terms
for the Fay relations satisfied by the reduced elliptic associator can be 
deduced explicitly from these. 
\vskip 1cm
\noindent {\bf \S 0. Introduction}
\vskip .3cm
The present article investigates the family of {\it Fay relations} satisfied
by the elliptic associator $A_\tau$ defined by Enriquez in [E], arising from
the definition of $A_\tau$ as an iterated integral of the classical 
Kronecker function (cf.~\S 2.1) which satisfies the well-known three-term Fay
relation.  The series $A_\tau$ is a group-like power series in
two non-commutative variables $a$ and $b$, with coefficients in a subring 
$R$ of ${\cal O}({\cal H})$ that was determined in [LMS] and whose definition 
is recalled below (cf.~\S 2.2).
The Fay relations for the elliptic associator were studied in [BMS]. They 
are expressed in the notation of \'Ecalle's mould theory, which we 
now briefly introduce (see [Ec], especially the beginning, for an 
overall introduction, or alternatively the introductory article
[S2] which contains complete proofs).

Fix a base ring $R$. In this article, the word {\it mould} will always refer 
to a rational-function mould, which is a family $A=A_r)_{r\ge 0}$ such that 
each $A_r$ lies in the ring of rational functions on an infinite number of
commutative variables $u_1,u_2,\ldots$, satisfying the unique property that
each $A_r$ is a function of $r$ variable; in particular, $A(\emptyset)\in R$.
We generally drop the subscript $r$ and write 
$A(u_1,\ldots,u_r)$ for $A_r$, which is called the {\it depth $r$ part} of the
mould $A$.  

A {\it polynomial mould} is one
for which each $A(u_1,\ldots,u_r)$ lies in the subring of polynomials
$R[u_1,u_2,\ldots,u_r]$.  
The set of 
moulds forms an $R$-module by adding moulds and multiplying them by scalars 
componentwise, i.e.~in each depth.

Consider the subspace of power series in $\Rab$ that lie in the kernel of 
the derivation mapping $a\mapsto 1$ and $b\mapsto 0$ (which is the case for 
all Lie-like and group-like power series and in particular all those 
considered in this article). These are precisely the power series that can be
written as power series in the variables $c_i=ad(a)^{i-1}(b)$ for $i\ge 1$.
There is a canonical map $ma$ from these power series to polynomial moulds
given by linearly extending the map on monomials
$$ma:c_{a_1}\cdots c_{a_r}\mapsto u_1^{a_1-1}\cdots u_r^{a_r-1}.$$

For any mould $A$, let $A'$ denote the mould defined by
$$A'(u_1,\ldots,u_r)={{1}\over{u_1\cdots u_r}}A(u_1,\ldots,u_r).\eqno(0.1)$$
Define the {\it Fay operator} $\F$ on a mould $B$ by the formula
$$\F(B)(u_1,\ldots,u_r)=
B(u_1,\ldots,u_r)+ B(u_2,\ldots,u_r,-\overline{u}_r)
+\sum_{i=1}^{r-1} B(u_2,\ldots,u_i,-\overline{u}_i,\overline{u}_{i+1},u_{i+2},\ldots,u_r)$$
where $\overline{u}_i=u_1+\cdots+u_i.$

Let $A(\tau)=ma(A_\tau)$ be the polynomial mould associated to the elliptic
associator $A_\tau\in \Rab$. 
The family of Fay relations satisfied by $A_\tau$ and studied in [BMS]
can be written in the following form: for each $r\ge 2$, we have
the equality
$$\F\bigl(A'(\tau)\bigr)(u_1,\ldots,u_r)
=\hbox{an undetermined correction term in lower depth.}$$

Generally speaking, we will say that a mould $A$ satisfies the {\it strict Fay relations} if
$\F(A')=0$, or {\it corrected Fay relations} if there is a mould $C_{A'}$
such that $\F(A')=C_{A'}$ and in each depth $r$, $C_{A'}$ is a sum of
linear terms in the $u_i$ and polynomial expressions in the parts of $A'$ 
of depth lower than $r$.  The main result of [BMS] is that the elliptic
associator $A_\tau$ satisfies a corrected family of Fay relations, although the
correction terms are not explicitly known.

In the present article we use mould theory to give a different interpretation
of how the Fay relations arise, which allows the explicit determination of the
correction terms.  There is however one caveat, which is that we can only
apply the main results of mould theory to this situation on the condition of
working modulo $2\pi i$, in the following sense.  The ring $R\subset {\cal O}({\cal H})$ generated by the coefficients of the power series $A_\tau$
contains the element $2\pi i$ (viewed as a constant function on ${\cal H}$).
The quotient ring $R'=R/\langle 2\pi i\rangle$ is non-trivial; its 
structure was determined completely in [LMS] and is recalled in \S 2.2 below.
Unfortunately, it so happens that reducing the coefficients of 
$A_\tau$ from $R$ to $\bar{R}$ gives zero, however this can be rectified 
by showing that the power series $A_\tau^{1/2\pi i}$
also has coefficients in $R$, but its reduction $\bar{A}_\tau$ 
is highly non-trivial.  

The main focus of this article is to determine the explicit family of Fay relations 
satisfied by the {\it reduced elliptic associator} $\bar{A}_\tau$.  On the way there,
we deduce the explicit (and much simpler) Fay relations satisfied by the
{\it reduced elliptic generating series} $\bar{E}_\tau$, and we prove the close
relationship between the Fay relations in the linearized context and the 
defining relations of the elliptic Kashiwara-Vergne Lie algebra.
\vskip .3cm
\noindent {\bf Outline of the article.}
\vskip .1cm
In \S 1, we recall the elements of mould theory necessary to define the 
strict and corrected Fay relations, and establish the following
mould-theoretic statement (Corollary of Theorem 1.3): 
\vskip .2cm
\noindent {\bf Theorem A.} {\it The space of polynomial alternal push-invariant
moulds $A$ satisfying a family of corrected Fay relations of the form
$$\F(A')(u_1,\ldots,u_r)=k_r(u_2+\cdots+u_r)$$
for a set of constants $(k_r)_{r\ge 2}$ is
isomorphic to the elliptic Kashiwara-Vergne Lie 
algebra $\krv_{ell}$}. 
\vskip .1cm
\noindent This theorem will be crucial in determining the Fay relations of
the Lie-like elliptic generating series.
\vskip .2cm
In \S 2, we recall the definition and construction of the reduced elliptic
generating series and elliptic associator. Let $\Z$ denote the $\QQ$-algebra of multiple zeta values, and $\overline{\Z}$
its quotient modulo the ideal generated by $\zeta(2)$. 
Let $\overline{\Phi}_{KZ}\in \overline{\Z}\langle\langle a,b\rangle\rangle$ denote the 
Drinfel'd associator reduced mod $\zeta(2)$.
This power series lies in the Grothendieck-Teichm\"uller group 
$GRT(\overline{\Z})$.  Enriquez in [En1] defined 
an elliptic version $GRT_{ell}$ with a canonical surjection to $GRT$,
and he also defined a section map
$$\Gamma:GRT\rightarrow GRT_{ell}.$$
In [En2], he defined a certain key automorphism $g_\tau$ of the power series
ring ${\cal O}({\cal H})\langle\langle a,b\rangle\rangle$. In [LMS], we 
defined a ring $R\subset {\cal O}({\cal H})$ with the following properties
(cf.~\S 2.1): $R$ contains $2\pi i$ and the ring $\bar{R}=R/\langle 2\pi i
\rangle$ is non-trivial, $g_\tau$ restricts to an automorphism of $\Rab$ that 
passes to an automorphism of $\RRab$. We then 
defined the {\it reduced elliptic generating series} by
$$\bar{E}_\tau:=g_\tau\Bigl(\Gamma(\overline{\Phi}_{KZ})\Bigr)\in GRT_{ell}
(\bar{R}).\eqno(0.2)$$
We further showed in [LMS] that there exists a unique automorphism $\Psi$ of
$\RRab$ such that $\Psi(a)=\bar{E}_\tau$ and $\Psi([a,b])=[a,b]$.
We define the {\it Lie-like elliptic generating series} by setting
$\bar{\be}_\tau=\psi(a)$ where $\psi$ is the derivation ${\rm log}(\Psi)$
of $\RRab$.

Let $t_{01}$ be the power series defined by
$$t_{01}={{ad(b)}\over{e^{ad(b)}-1}}(a)\in \RRab.\eqno(0.3)$$ 
One of the main results of [LMS] is the following alternative construction of
the reduced elliptic associator $\bar{A}_\tau$:
$$\bar{A}_\tau=\Psi(e^{t_{01}}).$$  
We conclude \S 2 by defining the {\it Lie-like reduced elliptic 
associator} $\bar{\ba}_\tau:=\psi(t_{01})$, which is 
obtained from $\bar{A}_\tau$ by a {\it double linearization}, by first
passing from $e^{t_{01}}$ to $t_{01}$, i.e.~taking 
${\rm log}\bar{A}_\tau=\Psi(t_{01})$, and then setting
$\bar{\ba}_\tau=\psi(t_{01})$ where $\psi={\rm log}(\Psi)$ as above.
\vskip .5cm
Finally, \S 3 is devoted to the Fay relations. The main results of \S 3.1
and \S 3.2 are stated in the following theorem, which regroups the statements
of Theorems 3.1 and 3.4.
Let $C=(c_r)_{r\ge 0}$ be the well-known constant correction mould
associated to the Drinfel'd associator $\overline{\Phi}_{KZ}$ given by
$$c_r=\cases{{{\zeta(r)}\over{r}}&for $r\ge 3$ odd\cr
0&otherwise.}\eqno(0.4)$$
and let $T=ma(t_{01}+a+{{1}\over{2}}[a,b])$ where $t_{01}$
is defined in (0.2). 
\vskip .3cm
\noindent {\bf Theorem B.} {\it The Lie-like reduced
elliptic generating series
satisfies the corrected Fay relations
$$\F(\bar{\be}'_\tau\bigr)(u_1,\ldots,u_r)=-rc_r(u_2+\cdots+u_r).$$
The Lie-like reduced elliptic associator satisfies the corrected
Fay relations
$$\F\bigl(\bar{\ba}'(\tau)\bigr)(u_1,\ldots,u_r)=
-rc_r(u_2+\cdots+u_r)
+\sum_{{{i=3}\atop{i\ {\rm odd}}}}^{r-1} \zeta(i)\Bigl(T'(u_2,u_3,\ldots,u_{r-i+1})-T'(u_{i+1},u_{i+2},\ldots,u_r)\Bigr).$$
In particular, we have
$$\F\bigl(\bar{\be}(\tau)\bigr)(u_1,\ldots,u_r)=\F\bigl(\bar{\ba}'(\tau)\bigr)(u_1,\ldots,u_r)=0\ \ \ \hbox{for all even values of}\ r.$$}
The first statement is proved in \S 3.1, using the identity $\bar{\be}_\tau
=\psi(a)$, mould theory and the equivalence between the Fay relation and the
elliptic Kashiwara-Vergne of Theorem A (Corollary of Theorem 1.3). The second statement is proved 
in \S 3.2, using the identity $\bar{\ba}_\tau=\psi(t_{01})$ and mould theory. 
The third statement follows from the definition of $C$ and the fact that
$T$ is zero in all odd depths.

Finally, \S 3.3 is devoted to the group-like situation: we show how to use the
closed forms of the correction terms for the Lie-like series to obtain 
a step-by-step computation of the explicit Fay relations satisfied by the
the group-like elliptic generating series $\bar{E}_\tau$ and elliptic 
associator $\bar{A}_\tau$, and illustrate the algorithm up to depth 4.
\vskip .3cm
\noindent {\bf Acknowledgments.} These results are based on the essential work of
Benjamin Enriquez bringing the well-known ``genus zero'' Drinfel'd associator 
into the elliptic framework.  I owe thanks to Nils Matthes for introducing me to
the notion of the Fay relations satisfied by the elliptic associator, and to both
Nils and Pierre Lochak for a fruitful collaboration summarized in the early parts
of this article, without which I would not have been able to make the connections
between the Fay relations, mould theory and the Kashiwara-Vergne Lie algebra.  

Above all, though, I am truly happy to take advantage of these acknowledgments to 
thank one more person, not directly related to this research but whose influence 
counted immensely in my early development as a Galois-and-number theorist.  This is,
of course, Moshe Jarden, to whom I am ever grateful for the role he played in my life 
as a mathematician, the encouragement he unfailingly gave me even after my research 
interests strayed from straightforward inverse Galois theory off into the distant 
realms of Grothendieck-Teichm\"uller theory, and the way in which he and his family 
were so welcomingly ready to throw open the doors of their home and their 
country to me at any time. Thanks to Moshe and his colleagues, my visits to Tel
Aviv University have counted among the most productive and inspiring periods of my life.
\vskip .5cm
\noindent {\bf \S 1. The Fay relations on moulds}
\vskip .3cm
\noindent {\bf \S 1.1. The Fay operator on moulds}
\vskip .2cm
For the purposes of this article, a {\it mould} is nothing other than a
family $(A_r)_{r\ge 0}$ of rational functions in commutative variables
over a given ring, such that for each $r$, $A_r$ is a rational function
of $r$ variables $u_1,\ldots,u_r$; in particular $A_0(\emptyset)$ is a
constant.  When there is no ambiguity, we drop the subscript and just
write $A(u_1,\ldots,u_r)$ instead of $A_r(u_1,\ldots,u_r)$, called {\it depth
$r$ part of $A$}. All the results in this section are general and hold over 
$\QQ$, which can be assumed to the base ring, but it can be replaced by
any larger ring, and indeed it will be in the subsequent sections.
Moulds can be multiplied 
by scalars and added componentwise (depth by depth).  We write $ARI$
for the vector space of moulds $A$ having $A(\emptyset)=0$.
We write $ARI^{pol}$ for the vector subspace of polynomial-valued moulds.

Let $\Qab$ denote the free completed associative algebra on two non-commutative
variables $a$, $b$, and consider the subspace ${\rm Ker}\partial$, where
$\partial$ is the derivation mapping $a\mapsto 1$, $b\mapsto 0$.  This subspace,
which we denote by $\QQ[C]$, is precisely the space of power series in the 
free variables $c_i=ad(a)^{i-1}(b)$.
There is a standard linear map from $\QQ[C]$ to $ARI^{pol}$
given by linearly extending the map on monomials
$$ma:c_{a_1}\cdots c_{a_r}\mapsto u_1^{a_1-1}\cdots u_r^{a_r-1}.$$

A mould $A\in ARI$ is said to be {\it alternal} if 
$$A\Bigl(sh\bigl(v_1,\ldots,v_i),(v_{i+1},\ldots,v_r)\bigr)\Bigr)=0,$$
for $1\le i\le r-1$,
where $sh$ denotes the shuffle operator, considered additively in the sense
that with our notation, we have for example
$$A\bigl(sh((u_1),(u_2))\bigr)=A(u_1,u_2)+A(u_2,u_1).$$
 We write $ARI_{al}$ for the
subspace of alternal moulds in $ARI$. It is well-known that the restriction
of $ma$ to polynomial alternal moulds induces an isomorphism 
$$ma:Lie[C]\buildrel\sim\over\rightarrow ARI^{pol}_{al},$$ 
where ${\rm Lie}[C]$ denotes the free Lie algebra on the
letters $c_i$, $i\ge 1$ viewed as a subspace of $\QQ[C]$.

The usefulness of moulds is that they easily allow us to work with denominators
and not only polynomials, thus generalizing the situation of power series in the
$c_i$, which is necessary, as we will see, to express the Fay relations 
and many other important properties that cannot be expressed within the
power series ring. Our first use of denominators will appear via the mould 
operator $dar$ and its inverse $dar^{-1}$ defined as follows:
fixed $A(\emptyset$), and for $r>0$ we set
$$darA(u_1,\ldots,u_r)=u_1\cdots u_r\,A(u_1,\ldots,u_r)
\ \ {\rm and}\ \ dar^{-1}A(u_1,\ldots,u_r)={{A(u_1,\ldots,u_r)}\over{u_1\cdots u_r}}.$$
Throughout this article we will use the notation $A':=dar^{-1}A$.  As
studied in [BMS], the Fay operator associated to the elliptic associator
$A_\tau$ actually acts on $A'_\tau$, which is the motivation for the
following definition.
\vskip .3cm
\noindent {\bf Definition 1.1.} The {\it Fay operator} on moulds, denoted 
$\F$, is defined by
$$\F(A)(u_1,\ldots,u_r):=A(u_1,\ldots,u_r)+
A(u_2,\ldots,u_r,-\overline{u}_r)
+\sum_{i=1}^{r-1} A(u_2,\ldots,u_i,-\overline{u}_i,\overline{u}_{i+1},u_{i+2},\ldots,u_r)$$
where $\overline{u}_i=u_1+\cdots+u_i.$
A mould $A\in ARI$ is said to {\it satisfy the strict Fay relations} if
$$\F(A')(u_1,\ldots,u_r)=0\eqno(1.1.1)$$
for $r\ge 2$. The mould $A$ is said to satisfy {\it a family of corrected
Fay relations} if there exists a ``correction'' mould $C_{A'}$ whose
depth $r$ part is a polynomial expression in the parts of $A'$ of depth $<r$ (possibly plus a linear expression in the $u_i$), such that
$$\F(A')(u_1,\ldots,u_r)=C_{A'}(u_1,\ldots,u_r).\eqno(1.1.2)$$
\vskip .4cm
\noindent {\bf \S 1.2. The Fay operator and swap circ-neutrality}
\vskip .3cm
We begin this subsection by introducing several more simple but important
mould operators, each of which preserves the value $A(\emptyset)$. Firstly, for
$r>0$, we have
$$durA(u_1,\ldots,u_r)=(u_1+\cdots +u_r)A(u_1,\ldots,u_r)$$
$$\Delta A(u_1,\ldots,u_r)=u_1\cdots u_r(u_1+\cdots+u_r)A(u_1,\ldots,u_r)$$
$$push_uA(u_1,\ldots,u_r)=A(-u_1-\cdots-u_r,u_1,\ldots,u_{r-1}).$$
We also have the swap operator
$$swapA(v_1,\ldots,v_r)=A(v_r,v_{r-1}-v_r,\ldots,v_1-v_2),$$
which sends moulds $A\in ARI$ to moulds $swapA\in \overline{ARI}$,
where $\overline{ARI}$ is identical to $ARI$ except that the moulds are
rational functions of commutative variables $v_1,v_2,\ldots$  (The use of
differently named variables for $\overline{ARI}$ is purely a convenience to
distinguish between moulds and their swaps.)

We now introduce three operators defined on $\overline{ARI}$. The first
is the inverse of the $swap$ (also denoted $swap$), defined by
$$swapB(u_1,\ldots,u_r)=B(u_1+\cdots+u_r,u_1+\cdots+u_{r-1},\ldots,u_1+u_2,u_1)$$
for $B\in \overline{ARI}$.
The second is the push-operator $push_v$ defined by
$$push_vB(v_1,\ldots,v_r)=B(-v_r,v_1-v_r,\ldots,v_{r-1}-v_r).$$
Finally, we introduce the $circ$-operator on $\overline{ARI}$, defined by
$$circB(v_1,v_2,\ldots,v_r)=B(v_2,\ldots,v_r,v_1).$$
Observe that the operators $push_u$ and $push_v$ are both of order $r+1$ in depth
$r$, and $circ$ is of order $r$.  Let $ARI^\Delta$ denote the subspace
of moulds $A\in ARI$ such that $\Delta A\in ARI^{pol}$; in other words,
$ARI^\Delta$ is the subspace of rational moulds whose denominator in depth 
$r$ is ``at worst'' $u_1\cdots u_r(u_1+\cdots+u_r)$.

\vskip .3cm
\noindent {\bf Lemma 1.1.} {\it We have the following relations between
mould operators:
\vskip .1cm
(i) $\Delta=dar\cdot dur=dur\cdot dar$
\vskip .1cm
(ii) $swap\cdot push_u^{-1}=push_v\cdot swap$
\vskip .1cm
(iii) A mould $A\in ARI$ is $push_u$-invariant if and only if $swapA$ is
$push_v$-invariant.}
\vskip .3cm
\noindent Proof. The first statement is obvious.  For the second, we
use the definitions to compute
$$\eqalign{(push_v\cdot swapA)(v_1,\ldots,v_r)&=
swapA(-v_r,v_1-v_r,\ldots,v_{r-1}-v_r)\cr
&=A(v_{r-1}-v_r,v_{r-2}-v_{r-1},\ldots,v_1-v_2,-v_1)\cr
&=push_u^{-1}A(v_r,v_{r-1}-v_r,\ldots,v_1-v_2)\cr
&=(swap\cdot push_u^{-1}A)(v_1,\ldots,v_r).}$$
For the third, we note that by (ii), $push_v\cdot swapA$ is equal to 
$swap\cdot push_u^{-1}A$ which is itself equal to $swapA$ thanks to the
$push_u$-invariance of $A$.\hfill{$\diamondsuit$}
\vskip .3cm
\noindent {\bf Theorem 1.2.} {\it Let $M\in ARI^\Delta$ be alternal and push-invariant. 
Let $A=\Delta M$ and $A'=dar^{-1}A=durM$.  Then 
the following are equivalent:
\vskip .2cm
(i) $A$ satisfies the strict Fay relations 
$$\F(A')(u_1,\ldots,u_r)=0;\eqno(1.2.1)$$
\vskip .1cm
(ii) $swapM$ satisfies the first alternality relation 
$$swapM\bigl(sh\bigl((v_1),(v_2,\ldots,v_r)\bigr)\bigr)=0;\eqno(1.2.2)$$
\vskip .1cm
(iii) $swapM$ satisfies the {\bf circ-neutrality relation}
$$swapM(v_1,\ldots,v_r)+swapM(v_2,\ldots,v_r,v_1)+\cdots+
swapM(v_r,v_1,\ldots,v_{r-1})=0.\eqno(1.2.3)$$}
\noindent Proof.  Let us show that (i) and (ii) are equivalent.
We do this by swapping the Fay relation (1.2.1):
$$\eqalign{swap&\Bigl(\F(A')(u_1,\ldots,u_r)\Bigr)=swap\Bigl(A'(u_1,\ldots,u_r)
+\sum_{i=1}^r A'(u_2,\ldots,u_i,-\overline{u}_i,\overline{u}_{i+1},u_{i+2},\ldots,u_r)\Bigr)\cr
&=A'(v_r,v_{r-1}-v_r,\ldots,v_2-v_3,v_1-v_2)\cr
&\qquad\sum_{i=1}^r A'(v_{r-1}-v_r,\ldots,v_{r-i+1}-v_{r-i+1},-v_{r-i+1},v_{r-i},v_{r-i-1}-v_{r-i},
\ldots,v_1-v_2)\cr
&=swap(A')(v_1,\ldots,v_r)
+\sum_{i=1}^r swap(A')(v_1-v_r,\ldots,v_{r-i}-v_r,-v_r,v_{r-i+1}-v_r,\ldots,v_{r-1}-v_r)\cr
&=swap(A')(v_1,\ldots,v_r)+swap(A')\Bigl(sh\bigl((-v_r),(v_1-v_r,\ldots,v_{r-1}-v_r)\bigr)\Bigr)\cr
&=swap(A')(v_1,\ldots,v_r)+\bigl(push_v\cdot swap(A')\bigr)\Bigl(sh\bigl((v_1),(v_2,\ldots,v_r)\bigr)\Bigr).}\eqno(1.2.4)$$
Thus the condition (1.2.1) of (i) is equivalent to
$$swap(A')(v_1,\ldots,v_r)+\bigl(push_v\cdot swap(A')\bigr)\Bigl(sh\bigl((v_1),(v_2,\ldots,v_r)\bigr)\Bigr)=0.\eqno(1.2.5)$$
We now rewrite (1.2.5) as the equivalent condition (1.2.6) obtained by 
applying $push_v^{-1}$ to (1.2.5):
$$push_v^{-1}\cdot swap(A')(v_1,\ldots,v_r)+swap(A')\Bigl(sh\bigl((v_1),(v_2,\ldots,v_r)\bigr)\Bigr)=0.\eqno(1.2.6)$$
Since we have $A'=durM$ and
$$swap\,durM(v_1,\ldots,v_r)=v_1\,swapM(v_1,\ldots,v_r),$$
we can use this to rewrite (1.2.6) in the equivalent form (1.2.7) as 
follows:
$$\eqalign{\bigl(push_v^{-1}&\cdot swap\,durM\bigr)(v_1,\ldots,v_r)+swap\,durM\Bigl(sh\bigl((v_1),(v_2,\ldots,v_r)\bigr)\Bigr)\cr
&=push_v^{-1}\bigl(v_1\,swapM(v_1,\ldots,v_r)\bigr)+v_1\,swapM(v_1,\ldots,v_r)+
v_2\,swapM\Bigl(v_2,sh\bigl((v_1),(v_3,\ldots,v_r)\bigr)\Bigr)\cr
&=(v_2-v_1)push_v^{-1}swapM(v_1,\ldots,v_r)+v_1\,swapM(v_1,\ldots,v_r)+
v_2\,swapM\Bigl(v_2,sh\bigl((v_1),(v_3,\ldots,v_r)\bigr)\Bigr)\cr
&=(v_2-v_1)swapM(v_1,\ldots,v_r)+v_1\,swapM(v_1,\ldots,v_r)
+ v_2\,swapM\Bigl(v_2,sh\bigl((v_1),(v_3,\ldots,v_r)\bigr)\Bigr)}$$
$$=v_2\, swapM\Bigl(sh\bigl((v_1),(v_2,\ldots,v_r)\bigr)\Bigr)=0,\eqno(1.2.7)$$
where the third inequality comes from the $push_v$-invariance of $swapM$
which follows from the $push_u$-invariance of $M$ and Lemma 1 (iii).
But clearly (1.2.7) is equivalent to the first alternality relation
(1.2.2). This proves the equivalence of (i) and (ii).

We now prove the equivalence of (ii) and (iii). Since $swapM$ is
$push_v$-invariant, the circ-neutrality expression
$$swapM(v_1,\ldots,v_r)+circ\,swapM(v_1,\ldots,v_r)+
\cdots+circ^{r-1}swapM(v_1,\ldots,v_r)\eqno(1.2.8)$$
can be written 
$$push_v^r\cdot swapM(v_1,\ldots,v_r)+circ\,push_vswapM(v_1,\ldots,v_r)
\qquad \qquad 
\qquad \qquad \qquad \qquad \qquad $$ 
$$+circ^2push_v^2swapM(v_1,\ldots,v_r)+\cdots+circ^{r-1}push_v^{r-1}\,swapM(v_1,\ldots,v_r).\eqno(1.2.9)$$
Let us compute each term of (1.2.9). For the first term, we have
$$\bigl(push_v^r\cdot swapM\bigr)(v_1,\ldots,v_r)=
swapM(v_2-v_1,\ldots,v_r-v_1,-v_1).\eqno(1.2.10)$$
For the second term, we have
$$(circ\,push_vswapM)(v_1,\ldots,v_r)=
push_vswapM(v_2,\ldots,v_r,v_1)=
swapM(-v_1,v_2-v_1,\ldots,v_{r-1}-v_1).\eqno(1.2.11)$$
Using the general formula for $push_v^i$ for $i>2$ 
$$push_v^iA(v_1,\ldots,v_r)=A(v_{r-i+2}-v_{r-i+1},v_{r-i+3}-v_{r-i+1},
\ldots,v_r-v_{r-i+1},-v_{r-i+1},v_1-v_{r-i+1},\ldots,v_{r-i}-v_{r-i+1}),$$
adding $i$ to every index (modulo $r$) we have
$$push_v^iA(v_{i+1},\ldots,v_r,v_1,\ldots,v_i)=
A(v_2-v_1,v_3-v_1,
\ldots,v_i-v_1,-v_1,v_{i+1}-v_1,\ldots,v_r-v_1),\eqno(1.2.12)$$
which is thus equal to
$$\bigl(circ^ipush_v^iswapM\bigr)(v_1,\ldots,v_r)$$
for $i=2,\ldots,r-1$.
But the sum of the terms (1.2.10), (1.2.11) and (1.2.12) for $i=2,\ldots,r-1$ is
exactly equal to
$$swapM\bigl(sh\bigl((-v_1),(v_2-v_1,v_3-v_1,\ldots,v_r-v_1)\bigr)\bigr),$$
yet at the same time to the circ-neutrality expression (1.2.8).
Therefore $swapM$ satisfies the first alternality relation if and only
if it is satisfies the circ-neutrality relation (1.2.3).
This concludes the proof of Theorem 1.2.
\hfill{$\diamondsuit$}
\vfill\eject
\noindent {\bf \S 1.3. Correction terms}
\vskip .3cm
In this subsection we generalize the statement of Theorem 1.2 to moulds
whose swaps only satisfy the first alternality relation and/or the 
circ-neutrality relation up to addition of a constant mould (i.e. a mould
whose value is a constant in each depth).
\vskip .2cm
\noindent {\bf Theorem 1.3.} {\it Let $M\in ARI^\Delta$ be alternal and push-invariant.  Let $A=\Delta M$ and $A'=durM=dar^{-1}(A)$. Then the following are equivalent for the same constant mould $C=(c_r)_{r\ge 0}$:
\vskip .1cm
(i) $swapM+C$ satisfies the first alternality relation;
\vskip .1cm
(ii) $swapM+C$ satisfies the circ-neutrality relation;
\vskip .1cm
(iii) $\F(A')(u_1,\ldots,u_r)=-rc_r(u_2+\cdots+u_r).$\hfill$(1.3.1)$}
\vskip .2cm
\noindent Proof. We saw in the proof of Theorem 1.2 that 
$$swapM(v_1,\ldots,v_r)+circ(swapM)(v_1,\ldots,v_r)+\cdots+
circ^{r-1}(swapM)(v_1,\ldots,v_r)$$
$$=swapM\bigl(sh\bigl((-v_1),
(v_2-v_1,v_3-v_1,\ldots,v_r-v_1)\bigr)\bigr).$$
If $swapM+C$ is circ-neutral, this means that
$$swapM(v_1,\ldots,v_r)+circ\,swapM(v_1,\ldots,v_r)+\cdots+
circ^{r-1}swapM(v_1,\ldots,v_r)=-rc_r$$
for $r\ge 2$; similarly, if $swapM+C$ satisfies the first alternality
relation, then we must have
$$=swapM\bigl(sh\bigl((v_1),
(v_2,v_3,\ldots,v_r)\bigr)\bigr)=-rc_r,$$
but this must then also hold for any variable change in the $v_i$, so we
must also have
$$=swapM\bigl(sh\bigl((-v_1),
(v_2-v_1,v_3-v_1,\ldots,v_r-v_1)\bigr)\bigr)=-rc_r.$$
Therefore circ-neutrality and the first alternality relation are still 
equivalent in the presence of a constant correction.

Now assume that $swapM+C$ satisfies the first alternality relation.
We saw in (1.2.7) that
$$v_2\, swapM\bigl(sh\bigl((v_1),(v_2,\ldots,v_r)\bigr)\bigr)
=\bigl(push_v^{-1}swap\,durM\bigr)(v_1,\ldots,v_r)+swap\,durM\bigl(sh\bigl((v_1),(v_2,\ldots,v_r)\bigr)\bigr).$$
Thus if 
$$swapM\bigl(sh\bigl((v_1),(v_2,\ldots,v_r)\bigr)\bigr)=-rc_r,$$
we have
$$\bigl(push_v^{-1}\cdot swap\,durM\bigr)(v_1,\ldots,v_r)+swap\,durM\bigl(sh\bigl((v_1),(v_2,\ldots,v_r)\bigr)\bigr)=-rc_rv_2.$$
Applying $push_v$ to both sides gives
$$swap\,durM(v_1,\ldots,v_r)+push_v\cdot swap\,durM\bigl(sh\bigl((v_1),(v_2,\ldots,v_r)\bigr)\bigr)=-r(v_1-v_r)c_r$$
or writing $A'=durM$,
$$swap(A')(v_1,\ldots,v_r)+push_v\cdot swap(A')\bigl(sh\bigl((v_1),(v_2,\ldots,v_r)\bigr)\bigr)=-r(v_1-v_r)c_r.$$
But by (1.2.4), this is then equivalent to
$$\bigl(swap\,\F(A')\bigr)(v_1,\ldots,v_r)=-r(v_1-v_r)c_r.$$
Therefore, swapping both sides, we end up with
$$\F(A')(u_1,\ldots,u_r)=-rc_r(u_2+u_3+\cdots+u_r).$$
\vskip .4cm
\noindent {\bf \S 1.4. The elliptic Kashiwara-Vergne Lie algebra}
\vskip .3cm
In [En1], Enriquez defined an elliptic version $\grt_{ell}$ of the 
Grothendieck-Teichm\"uller Lie algebra $\grt$, contained in the vector
space ${\rm Der}^0{\rm Lie}[a,b]$ of derivations of the free Lie algebra
${\rm Lie}[a,b]$ that annihilate the bracket $[a,b]$; Enriquez explicitly 
constructed
a section map $\gamma:\grt\rightarrow \grt_{ell}$.  In [AT], a simple
proof of the Lie algebra inclusion $\iota:\grt\hookrightarrow\krv$ was found.
The article [RS] is devoted to the construction of an elliptic version 
$\krv_{ell}$ of the Kashiwara-Vergne Lie algebra, also as a subspace
of ${\rm Der}^0{\rm Lie}[a,b]$, that has various 
good properties including a section map $\gamma$ from the Kashiwara-Vergne Lie 
algebra $\krv$ to $\krv_{ell}$ extending Enriquez's section, making the diagram
$$\xymatrix{\grt\ar[d]_\gamma\ar[rr]^\iota&&\krv\ar[d]^\gamma\\
\grt_{ell}\ar[rd]&&\krv_{ell}\ar[ld]\\
&{\rm Der}^0{\rm Lie}[a,b]&}$$
commute.  The following definition of the elliptic Kashiwara-Vergne Lie
algebra given in [RS] shows how closely the Fay relations are 
related to $\krv_{ell}$.
\vskip .3cm
\noindent {\bf Definition 1.2.} The underlying vector space of the 
elliptic Kashiwara-Vergne Lie algebra $\krv_{ell}$ is the space of power series 
$f\in \QQ\langle\langle a,b\rangle\rangle$ whose associated mould 
$F\in ARI^{pol}$ satisfies the following properties:
$M:=\Delta^{-1}(F)$ is alternal, push-invariant and there exists
a constant mould $C$ such that $swapM+C$ is circ-neutral.
\vskip .3cm
In view of Theorem 1.3, however, we now have another equivalent definition 
using the Fay relations, as follows.
\vskip .3cm
\noindent {\bf Corollary of Theorem 1.3.} {\it The underlying vector space of 
the elliptic Kashiwara-Vergne Lie algebra $\krv_{ell}$ is isomorphic to
the vector space of power series $f\in \QQ\langle\langle a,b\rangle\rangle$ 
whose associated mould $F\in ARI^{pol}$ satisfies the following properties:
$F$ is alternal and push-invariant and there exists a constant mould 
$C=(c_r)_{r\ge 2}$
such that for $r\ge 2$, we have
$$\F(A')(u_1,\ldots,u_r)=-rc_r(u_2+\cdots+u_r).$$}

\noindent {\bf Remark.} It is proved in [RS] that the image of the inclusion 
map of $\krv_{ell}$ into ${\rm Der}^0{\rm Lie}[a,b]$ is closed under the 
natural Lie bracket on ${\rm Der}^0{\rm Lie}[a,b]$, making $\krv_{ell}$
into a Lie algebra. 
\vskip .8cm
\noindent {\bf \S 2. The elliptic associator $A_\tau$ and its reduced, Lie and 
mould versions}
\vskip .4cm
\noindent {\bf \S 2.1. Enriquez' elliptic associator $A_\tau$}
\vskip .3cm
Let $\Z$ denote the $\QQ$-algebra of multizeta values, and let
$\overline{\Z}$ be the quotient of this algebra by the ideal generated by
$\zeta(2)$.  Let $\Qxy$ denote the free completed polynomial algebra 
on two non-commutative variables $x$ and $y$\footnote{$^*$}{In this introductory
section containing definitions, we use the variables $x$ and $y$ so as to 
conform with the usual notation in articles on the subject.  In the next
section we will identify these variables with $a$ and $b$ via $x=b$, $y=a$.}.  Let $\Phi_{KZ}(x,y)\in
\Z\langle\langle x,y\rangle\rangle$ denote the Drinfel'd associator, and let us define 
three Bernoulli power
series inside $\QQ\langle\langle x,y\rangle\rangle$ by
$$t_{01} := Ber_x(-y),\ \ \ \ t_{02}:=Ber_{-x}(y),\ \ \ \ t_{12}:=[y,x],\eqno(2.1.1)$$
where 
$$Ber_x(y):={{ad(x)}\over{e^{ad(x)}-1}}(y).$$
Let $\epsilon_{2i}$, $i\ge 0$, denote the well-known derivations defined on
$\Qxy$ by
$$\epsilon_{2i}(x)=ad(x)^{2i}y,\ \ \ \epsilon_{2i}([x,y])=0.$$
Note that the second condition determines the value of $\epsilon_{2i}(y)$
uniquely.  In particular, $\epsilon_0(x)=y$, $\epsilon_0(y)=0$.
We write $\fu$ for the Lie subalgebra of ${\rm Der}\,{\rm Lie}[x,y]$
generated by the $\epsilon_{2i}$, and
$\Ufu$ for its enveloping algebra.

\vskip .2cm
\noindent {\bf Definition 2.1.} Let 
$g_\tau$ be the automorphism of $\Qxy$ defined in [En2, Prop. 5.1] (see also
[LMS, \S 2.3]).
It is a solution of the differential equation
$${{1}\over{2\pi i}}{{\partial}\over{\partial \tau}}g_\tau=
-\Bigl(\epsilon_0+\sum_{k\ge 1} {{2}\over{(2k-2)!}}G_{2k}(\tau)\epsilon_{2k}\Bigr)g_\tau, $$
where
$$G_{2k}(\tau)=-{{B_{2k}}\over{4k}}+\sum_{n\ge 1} \sigma_{2k-1}(n)q^n$$
is the Hecke-normalized Eisenstein series, with
$q=e^{2\pi i\tau}$. Enriquez singles out the solution $g_\tau$
by specifying its asymptotic behavior as $\tau\rightarrow i\infty$.

An explicit expression for $g_\tau$ can be given in the form
$$g_\tau={\rm id}+\sum_{\uk=(2k_1,\ldots,2k_r)} {\cal G}_\uk \epsilon_{\uk},\eqno(2.1.2)$$
where the sum runs over all tuples $(2k_1,\ldots,2k_r)$ of even integers
$\ge 0$, with $r\ge 1$, and ${\cal G}_\uk$ is defined as a certain
regularized iterated integral of $G_{2k_1},\ldots,G_{2k_n}$ [LMS, \S 2.2].
The operators $\epsilon_{\uk}$ are not linearly independent. Choose
a basis for the $\QQ$-vector space they generate, write $g_\tau$ in 
that basis, and let ${\cal E}^{\rm geom}\subset {\cal O}({\cal H})$ denote the $\QQ$-vector space spanned by
the coefficients in $g_\tau$ of the basis elements. These coefficients are all linear 
combinations of the $\G_{\uk}$, so ${\cal E}^{\rm geom}$ is a subspace of the $\QQ$-vector space spanned by the $\G_{\uk}$, and it is independent of the
choice of basis.  It was shown in [LMS] that ${\cal E}^{\rm geom}$ is not
merely a vector space but actually forms a subring of the ring of functions
${\cal O}(\cal H)$ on the Poincar\'e upper half-plane (this is simply 
a direct consequence of the definition of $g(\tau)$ which shows it to be a
group-like power series).  We can also view the multizeta value 
algebra $\Z[2\pi i]$ as a subring of constant functions of 
${\cal O}(\cal H)$.  It is shown in [LMS] (Cor.~2.10) that the subring of ${\cal O}(\cal H)$ generated by ${\cal E}^{\rm geom}$ and by $\Z[2\pi i]$ is isomorphic to the tensor product
${\cal E}^{\rm geom}\otimes_\QQ \Z[2\pi i]$.  
\vskip .2cm
\noindent {\bf Definition 2.2.} Let $R$ denote the subring 
${\cal E}^{\rm geom}\otimes_\QQ \Z[2\pi i]\subset {\cal O}({\cal H})$. Set
$$A=\Phi_{KZ}(t_{01},t_{12})e^{2\pi it_{01}}\Phi_{KZ}(t_{01},t_{12})^{-1}
\in \Z[2\pi i]\langle\langle x,y\rangle\rangle.\eqno(2.1.3)$$
The {\it elliptic associator} is the group-like power series 
$$A_\tau=g_\tau(A)\in R\langle\langle x,y\rangle\rangle.\eqno(2.1.4)$$
\vskip .2cm
This elliptic associator was introduced by Enriquez in [En1].
The formula (2.1.4) is not the original definition; Enriquez
views it rather as a property of the elliptic associator, which can be
defined in various more intrinsic ways, in particular via iterated integrals
of the Kronecker function as in (2.1.5) and (2.1.6) below  ([En2, \S 6.2]). 
 However, (2.1.4) emerges as the most apt definition for our purposes.  The 
definition by iterated integrals, on the other hand, which we now describe, 
is that studied in [BMS] which gave rise to the discovery
of the Fay relations satisfied by $A_\tau$.
\vskip .2cm
\noindent {\bf Definition 2.3.} Let $F\pmatrix{u\cr v}$ denote the
Kronecker function 
$$F\pmatrix{u\cr v}={{\theta(u+v)}\over{\theta(u)\theta(v)}}$$
where $\theta(u)$ is the only Jacobi theta-function that is odd in $z$,
$$\theta_1(z,q)=2\sum_{n=0}^\infty (-1)^nq^{(n+1/2)^2}e^{(2n+1)iz}$$
with $q=e^{i\pi\tau}$.  This function satisfies
$$\theta(z+\tau)=qe^{-2i\pi z}\theta(z),\ \ \ 
\theta(z+1)=-\theta(z),\ \ \ {{d\theta}\over{dz}}(0)=1,$$
and the only zeros of $\theta$ are the points of the lattice $\Zi+\Zi\tau$.
\vskip.2cm
The Kronecker function satisfies the famous elementary {\it Fay relation}
$$F\pmatrix{u_1\cr v_2}F\pmatrix{u_2\cr v_2}=
F\pmatrix{u_1+u_2\cr v_1}F\pmatrix{u_2\cr v_2-v_1}
+F\pmatrix{u_1+u_2\cr v_2}F\pmatrix{u_1\cr v_1-v_2}.\eqno(2.1.5)$$
Let 
$$\widetilde{A}_\tau(u_1,\ldots,u_r)=\int_{0<v_r<\cdots<v_1<1}F\pmatrix{u_1\cr v_1}\cdots
F\pmatrix{u_r\cr v_r} dv_r\cdots dv_1$$
be the iterated integral, regularized as usual.  
The elliptic associator is given by the formula
$$A_\tau=e^{{{i\pi}\over{2}}[x,y]}\ (2\pi)^{-[x,y]}\
\widetilde{A}_\tau(x/2\pi i,2\pi iy)\ (2\pi)^{[x,y]}\
e^{{{i\pi}\over{2}}[x,y]}.\eqno(2.1.6)$$
It is shown in [En1] that $A_\tau$ satisfies (2.1.4).

The elementary Fay relation satisfied by the Kronecker function 
``propagates'', in a manner described in [BMS]. 
The result shown there is that, writing $A(\tau):=ma(A_\tau)$ for the
polynomial-valued mould associated to the power series $A_\tau$, and 
$A'(\tau)=dar^{-1}A(\tau)$, the mould
$A(\tau)$ satisfies a family of Fay relations of the form (1.1.2),
in which the correction mould $C_{A'_\tau}$ arises from the need to regularize 
the iterated integrals defining $A_\tau$. The correction term was 
not explicitly determined in [BMS].
\vskip .4cm
\noindent {\bf \S 2.2. Reducing the elliptic associator mod $2\pi i$}
\vskip .3cm
The goal of this article is to show that, at least modulo $2\pi i$, the mould
construction of the elliptic associator based on the Drinfel'd associator
reveals a different explanation for the Fay relations, and furthermore allows 
us to explicitly determine the correction term, which arises 
directly from that of the Drinfel'd associator. 

From now until the end of this article, we make the variable change 
$a:=y$, $b:=x$, and consider $t_{01}$, $t_{02}$, $t_{12}$, $A$, and $A_\tau$ 
from (2.1.1), (2.1.3) and (2.1.4) as power series in the variables $a,b$.  
The reason for this is that it is a more convenient notation for the
application of mould theory starting in \S 2.3 below.

Let $R={\cal E}^{\rm geom}\otimes_\QQ \Z[2\pi i]$ as in Definition 2.2, and let
$\bar{R}$ denote the quotient of this ring by the ideal generated by
$1\otimes 2\pi i$; thus
$$\bar{R}\simeq {\cal E}^{\rm geom}\otimes_\QQ \overline{\Z}$$
where $\overline{\Z}=\Z/\bigl(\zeta(2)\bigr)$.

Since $A$ is a conjugate of $e^{2\pi it_{01}}$, the same holds for $A_\tau=
g_\tau(A)$, and therefore every coefficient in the power series $A_\tau$ 
vanishes in $\bar{R}$. 
\vskip .3cm
\noindent {\bf Lemma 2.1.} {\it The power series $A_\tau^{1/2\pi i}$ lies in 
$\Rab$,}
\vskip .2cm
\noindent Proof. We first show that the power series $A^{1/2\pi i}$ lies
in $\Z\langle\langle a,b\rangle\rangle$.  Indeed, this follows directly from
(2.1.3), since we have
$$A^{1/2\pi i}=\Phi_{KZ}(t_{01},t_{12})e^{t_{01}}\Phi_{KZ}(t_{01},t_{12})^{-1},
\eqno(2.2.1)$$
and all coefficients of the Drinfel'd associator lie in $\Z$. In order to
show that $A_\tau^{1/2\pi i}\in \Rab$, we note that 
$$A_\tau^{1/2\pi i}=g_\tau(A^{1/2\pi i}).$$ 
Since the coefficients of
$g_\tau$ generate ${\cal E}^{\rm geom}$, all the coefficients of 
$g_\tau(A^{1/2\pi i})$ lie in the ring generated by ${\cal E}^{\rm geom}$ and
$\Z$, which is a subring of $R$.\hfill{$\diamondsuit$}
\vskip .3cm
\noindent {\bf Definition 2.4.} The {\it reduced elliptic associator} 
$\bar{A}_\tau\in \RRab$ is the power series obtained by reducing the 
coefficients of $A_\tau^{1/2\pi i}$ from $R$ to $\bar{R}$.
\vskip .3cm
It is easy to see that $\bar{A}_\tau$ is highly non-trivial; indeed,
(2.2.1) shows that $A^{1/2\pi i}$ is already non-trivial.  Indeed, it is shown
in [LMS] (Thm.~3.6) that the coefficients of $\bar{A}_\tau$ together
with the single element $2\pi i\tau$ generate the entire ring $\bar{R}$.
\vskip .3cm
\noindent {\bf \S 2.3. Automorphism construction of the reduced group-like and Lie-like elliptic associators}
\vskip .3cm
The main result that we will use for our analysis of the power series
$\bar{A}$ and $\bar{A}_\tau$ is the construction
of these series as the images of $e^{t_{01}}$ by
specific automorphisms of the ring $\RRab$. Let
$$\Gamma:GRT\rightarrow GRT_{ell}$$
denote the section map defined by Enriquez in [En1] from the 
prounipotent Grothendieck-Teichm\"uller group to its elliptic version. 
Enriquez also defined a graded Lie algebra version $\grt_{ell}$ of $GRT_{ell}$
equipped with a natural surjection $\grt_{ell}\rightarrow \grt$ and
a Lie algebra section 
$$\gamma:\grt\rightarrow \grt_{ell}$$ 
(cf.~\S 1.4).  The reduced Drinfel'd associator $\overline{\Phi}_{KZ}$ lies in
$GRT(\overline{\Z})$, and its image $\Gamma(\overline{\Phi}_{KZ})$ lies in
the group $GRT_{ell}(\overline{\Z})$, which Enriquez considers 
as a subgroup of automorphisms of
the ring $\ZZab$ fixing the bracket $[a,b]$.
Let $\Psi_{KZ}$ denote the automorphism of $\ZZab$ fixing $[a,b]$ given
by $\Gamma(\overline{\Phi}_{KZ})\in GRT_{ell}(\overline{\Z})$.
Like $\Psi_{KZ}$, the automorphism $g_\tau$ of defined in \S 2.1 
can also be viewed as an automorphism of the ring
$\ZZab$ fixing $[a,b]$, since each derivation $\epsilon_{2i}$ annihilates
$[a,b]$ and $g_\tau$ is a group-like power series in the $\epsilon_{2i}$.

Let $r_\tau$ and $\psi_{KZ}$ denote the derivations of $\ZZab$
$$r_\tau=log(g_\tau),\ \ \ \ \psi_{KZ}=log(\Psi_{KZ}),$$
and set
$$\Psi:=g_\tau\circ \Psi_{KZ}=
e^{r_\tau}\circ e^{\psi_{KZ}}=e^{\psi}
\in {\rm Aut}(\RRab),\eqno(2.3.1)$$
where letting $CH$ denote the Campbell-Hausdorff law in the derivation
Lie algebra ${\rm Der}\,{\rm Lie}[a,b]$, we have 
$$\psi:=CH(r_\tau,\psi_{KZ}).\eqno(2.3.2)$$
\vskip .2cm
\noindent {\bf Definition 2.5.} 
The {\it reduced elliptic generating series} is the power series 
$\bar{E}_\tau$ given by
$$\bar{E}_\tau:=\Psi(a)\in \RRab.\eqno(2.3.3)$$
The {\it reduced Lie-like elliptic generating series} is defined by
$$\bar{\be}_\tau:=\psi(a).\eqno(2.3.4)$$
By the previous section, the {\it reduced elliptic associator} is 
given by
$$\bar{A}_\tau=\Psi(e^{t_{01}}).\eqno(2.3.5)$$
Finally, the {\it reduced Lie-like elliptic associator} is defined by
$$\bar{\ba}_\tau:=\psi(t_{01})\eqno(2.3.6)$$
is called the {\it reduced Lie-like elliptic associator}.
\vskip .3cm
\noindent The power series $\bar{E}_\tau$ was introduced in [LMS].
Indeed, it was shown there that up to adjoining the single element $2\pi i\tau$,
the coefficients of $\bar{E}_\tau$ generate the full ring $\bar{R}$ of reduced 
elliptic multiple zeta values, and the same statement holds for the coefficients
of $\bar{A}_\tau$.  
\vskip .2cm
\noindent {\bf Lemma 2.2.} {\it The power series $\psi_{KZ}(a)$ lies in
the elliptic Kashiwara-Vergne Lie algebra $\krv_{ell}$.  More precisely,
setting $\bar{\bf e}:=ma\bigl(\psi_{KZ}(a)\bigr)$, the mould $\Delta^{-1}\bar{\bf e}$ is alternal, push-invariant and $swap(\bar{\bf e})+C$ is circ-neutral,
where $C$ is the mould in (0.4).}
\vskip .2cm
\noindent Proof. We have $\psi_{KZ}(a)=\gamma(\varphi_{KZ})$, where
$\varphi_{KZ}\in \grt\subset \krv$ denotes the Lie-like Drinfel'd associator and
$\gamma:\krv\rightarrow\krv_{ell}$ is the Enriquez section extended to
$\krv$, as recalled in \S 1.4. Let $\nabla:=\Delta^{-1}\circ \gamma$.
Let $\bar{\bf e}=ma(\varphi_{KZ})$.  Then $\bar{\bf e}$ is alternal since $\grt\subset
{\rm Lie}[a,b]$.  Taking $F=\bar{\bf e}$ in [S1, Theorem 1.3.2] (in which the morphism 
$\nabla$ is called $Ad_{ari}(invpal)$ following \'Ecalle's notation), we see 
that $swap(\nabla(\bar{\bf e})\bigr)+C$ is alternal for the $C$ in (0.4), which is
the usual constant correction mould associated to $\varphi_{KZ}$.  
Thus $\bar{\bf e}$ is bialternal, i.e.~alternal and with alternal swap up to adding
the constant mould $C$. But all bialternal moulds that are even in depth 1
(which is the case for $\bar{\be}$) are push-invariant by [S2], Lemma 2.5.5.
Therefore $\bar{\be}$ is alternal and push-invariant and $swap(\bar{\bf e})$
is alternal, so in particular it satisfies the first alternality relation,
and therefore by Theorem 1.3, $swap(\bar{\bf e})+C$ satisfies the
circ-neutrality relation. This concludes the proof.\hfill{$\diamondsuit$}
\vskip .5cm
\noindent {\bf \S 2.4. Mould versions of the Lie-like reduced elliptic 
generating series and associator}
\vskip .3cm
To express the action of derivations in mould terms, we need to introduce
a few more basic mould theory notions.  To start with,
the standard formula for mould multiplication is given by
$$mu(A,B)(u_1,\ldots,u_r)=\sum_{i=0}^r A(u_1,\ldots,u_i)B(u_{i+1},\ldots,u_r).$$
It is easy to see that the $mu$ multiplication extends the usual multiplication of 
power series in that
$$ma(fg)=mu\bigl(ma(f),ma(g)\bigr)$$
for $f,g\in \QQ[C]$. Recall that $ARI$ denotes the vector space of 
(rational) moulds satisfying $A(\emptyset)=0$.
One can make $ARI$ into a Lie algebra by equipping it with the
standard Lie bracket 
$$lu(A,B)=mu(A,B)-mu(B,A).$$
The vector space $ARI$ equipped with the $lu$ bracket is a Lie algebra
denoted $ARI_{lu}$.  We can add a certain very important yet ``trivial'' mould, called $a$, to 
$ARI_{lu}$; this mould takes the value $a$ on the empty set and $0$ in all 
depths $\ge 1$.  We write $ARI^a$ for the vector space generated by $ARI$ 
and the mould $a$, and $ARI^a_{lu}$ for the Lie algebra obtained by extending 
the bracket to $ARI^a$ via the formula
$$lu(Q,a)=durQ.$$
We note that this formula merely extends to all moulds in ARI the known 
identity for a power series $f$:
$$ma\bigl([f,a]\bigr)=lu(ma(f),a)=dur\bigl(ma(f)\bigr)$$
([R], App. A, Prop. 5.1, see also [S2], Lemma 3.3.1).  The advantage of
adding the mould $a$ to $ARI$ is that extending $ma$ to the letter $a$
by mapping it to the mould also denoted by $a$ makes $ma$ extend to isomorphisms
$$ma:\Qab\rightarrow ARI^{a,pol}$$
and restrict to an isomorphism
$$ma:{\rm Lie}[a,b]\rightarrow (ARI^a)^{pol}_{al}.$$
Thus, using a bracket notation for moulds associated to power series, we
may set 
$$\cases{\bar{E}(\tau):=ma(\bar{E}_\tau)\cr
\bar{A}(\tau):=ma(\bar{A}_\tau)\cr
g(\tau):=ma\bigl(g_\tau(a)\bigr)\cr
\overline{\be}(\tau):=ma(\overline{\be}_\tau)\cr
\bar{\ba}(\tau):=ma(\bar{\ba}_\tau).}$$
\vskip .3cm

Set $U_1=ma([a,b])$, so that $U_1(u_1)=-u_1$.
For every mould $P\in ARI$, there exists an associated derivation
$Darit(P)$ of $ARI^a_{lu}$ (defined explicitly in Appendix A)
having the two following properties:
\vskip .1cm
i) $Darit(P)\cdot a=P$;
\vskip .1cm
ii) $Darit(P)\cdot U_1=0$.\hfill(2.4.1)
\vskip .3cm
Thanks to these two properties, we see that knowing the value of the derivation
$Darit(P)$ on $a$ yields the mould $P$ uniquely.  In particular, therefore,
the derivation $Darit\bigl(\bar{\be}(\tau)\bigr)$ restricted to 
$(ARI^a)^{pol}_{al}$ is equal to the derivation $\psi$ defined in (2.3.2), 
in the sense that for every $f\in {\rm Lie}[a,b]$, we have
$$ma\bigl(\psi(f)\bigr)=Darit\bigl(\bar{\be}(\tau)\bigr)\cdot ma(f).\eqno(2.4.2)$$
Since by Definition 2.5 we have $\bar{\be}_\tau=\psi(a)$ and 
$\bar{\ba}_\tau=\psi(t_{01})$,
we can translate these into the mould equalities
$$\bar{\be}(\tau):=ma\bigl(\bar{\be}_\tau\bigr)=ma\bigl(\psi(a)\bigr)=Darit\bigl(\bar{\be}(\tau)\bigr)\cdot a.\eqno(2.4.3)$$
$$\bar{\ba}(\tau):=ma\bigl(\bar{\ba}_\tau\bigr)=ma\bigl(\psi(t_{01})\bigr)=Darit\bigl(\bar{\be}(\tau)\bigr)\cdot ma(t_{01}).\eqno(2.4.4)$$
Finally, since $\psi=log(\Psi)$ with $\Psi(a)=\bar{E}_\tau$ and
$\bar{A}_\tau=\Psi(e^{t_{01}})$,
we also have the mould versions
$$\bar{E}(\tau)={\rm exp}\Bigl(Darit\bigl(\bar{\be}(\tau)\bigr)\Bigr)\cdot a=\sum_{n\ge 0} {{1}\over{n!}}Darit\bigl(\bar{\be}(\tau)\bigr)^n\cdot a\eqno(2.4.5)$$
and
$$\bar{A}(\tau)={\rm exp}\Bigl(Darit\bigl(\bar{\be}(\tau)\bigr)\Bigr)\cdot ma(e^{t_{01}})=\sum_{n\ge 0} {{1}\over{n!}}Darit\bigl(\bar{\be}(\tau)\bigr)^n\cdot 
ma(e^{t_{01}}).\eqno(2.4.6)$$

In the next section \S 3, we will use the results of \S 1 and the expressions
(2.4.3) and (2.4.4) to compute the Fay correction of the Lie-like elliptic 
generating series and the Lie-like elliptic associators, and then
based on the expressions (2.4.5) and (2.4.6), we will show how to 
deduce the Fay corrections for the group-like versions.
\vskip .6cm
\noindent {\bf \S 3. The Fay relations of the elliptic generating
series and elliptic associator}
\vskip .3cm
The results of the previous sections make it quite natural and easy to
determine the Fay relations satisfied by the elliptic generating series
and the elliptic associator, particularly in their Lie-like forms
(and mod $2\pi i$).
\vskip .3cm
\noindent {\bf \S 3.1. The Fay relations for the Lie-like elliptic
generating series}
\vskip .3cm
\noindent {\bf Theorem 3.1.} {\it The Lie-like elliptic generating series
$\overline{\be}(\tau)$ satisfies the Fay relations 
$$\F\bigl(\overline{\be}'(\tau)\bigr)=\cases{-\zeta(r)(u_2+\cdots+u_r)&$r\ge 2$ odd\cr
0&$r\ge 2$ even.}\eqno(3.1.1)$$}

\noindent {\bf Proof.}  Recall that $\bar{\be}_\tau=\psi(a)$ where
$\psi=CH(r_\tau,\psi_{KZ})$ from (2.3.2).  Let $\varphi_{KZ}\in\grt$ be the
Lie-like version of the Drinfel'd associator, which lies in the 
Grothendieck-Teichm\"uller Lie algebra.  Then we have
$\psi_{KZ}=\gamma(\varphi_{KZ}),$
where $\gamma:\grt\rightarrow\grt_{ell}$ is the Enriquez section on the
level of the Lie algebras.
As recalled in \S 1.4, it was shown in [AT]
that there is an inclusion $\iota:\grt\rightarrow \krv$, and in [RS] that the 
section map $\gamma:\grt\rightarrow\grt_{ell}$ extends to a section map $\gamma:\krv\rightarrow
\krv_{ell}$.  Therefore $\gamma(\varphi_{KZ})\in \krv_{ell}$, i.e.
$\psi_{KZ}\in \krv_{ell}$.  Furthermore, it is shown in [RS] that all the
derivations $\epsilon_{2i}$ defined in \S 2.1 lie in $\krv_{ell}$.
Let $\frak{u}$ denote the Lie subalgebra of $\krv_{ell}$ generated by
the $\epsilon_{2i}$.  It is shown in [LMS] (just following Lemma 3.1, based
on a result in [HM]) that any Lie bracket of
an element of $\frak{u}$ with $\gamma(\varphi_{KZ})$ lies in $\frak{u}$.
In particular, any Lie bracket of the derivations 
$r_\tau$ and $\gamma(\varphi_{KZ})=\psi_{KZ}$ lies in $\frak{u}$.
Therefore the Campbell-Hausdorff product $\psi=CH(r_\tau,\psi_{KZ})$ can be
written $\psi=\psi_{KZ}+D$ with $D\in \frak{u}$.
Since $\psi_{KZ}\in
\krv_{ell}$ and the $\epsilon_{2i}$, $i\ge 0$ are in $\krv_{ell}$, this shows that
$\psi\in \krv_{ell}$, so in terms of moulds we have
$$ma\bigl(\psi(a)\bigr)=\bar{\be}(\tau)\in ma(\krv_{ell}).$$
Thus by Definition 1.2, $\Delta^{-1}\overline{\be}(\tau)$ is alternal, 
push-invariant and its swap is circ-neutral up to addition of a constant
mould that we now determine.

Let $\bar{\bf e}:=\psi_{KZ}(a)$ as in Lemma 2.2, and let ${\cal D}:=ma\bigl(D(a)\bigr)$.  Since
$\psi=\psi_{KZ}+D$, and $\bar{\be}(\tau)=ma\bigl(\psi(a)\bigr)$ and
$\bar{\bf e}=ma\bigl(\psi_{KZ}(a)\bigr)$, we have
$$swap\bigl(\Delta^{-1}\bar{\be}(\tau)\bigr)=swap\bigl(\Delta^{-1}\bar{\bf e}\bigr)+swap\bigl(\Delta^{-1}{\cal D}\bigr).$$
To determine the constant correction mould for 
$swap\bigl(\Delta^{-1}\bar{\be}(\tau)\bigr)$, we determine the constant
corrections of $swap\bigl(\Delta^{-1}{\cal D}\bigr)$ and $swap\bigl(\Delta^{-1}\bar{\bf e}\bigr)$ separately.

We first show that $swap\bigl(\Delta^{-1}{\cal D}\bigr)$ is 
strictly circ-neutral.  To see this, let $U_{2i}=ma\bigl(\epsilon_{2i}(a)\bigr)$ be the
mould associated to each derivation $\epsilon_{2i}$; then $U_{2i}(u_1)=
u_1^{2i}$ and $U_{2i}$ is zero in every depth different from 1. We
have $\Delta^{-1}U_{2i}=U_{2i-2}$. In particular, we see that
the moulds $swap(U_{2i-2})$ are all automatically circ-neutral simply
because circ-neutrality is a property that only concerns the parts of a 
mould in depths $\ge 2$.  The moulds $U_{2i-2}$ are also trivially alternal
for the same reason, and they are also push-invariant since they are
even in depth 1.  It is shown in [RS] that the vector space of alternal, 
push-invariant moulds with strictly circ-neutral swap forms a Lie algebra;
in particular, if $\delta$ denotes the image under $ma$ of any bracket of 
the derivations $\epsilon_{2i}$, then $\Delta^{-1}\delta$ has strictly 
circ-neutral swap.  This holds in particular for $\delta=D$, which is a sum 
of Lie brackets of the $\epsilon_{2i}$. Thus 
$\Delta^{-1}{\cal D})$ has strictly circ-neutral swap. 

We now need to determine the constant correction mould of 
$swap(\Delta^{-1}\bar{\bf e})$; however, this was done in Lemma 2.2,
where we showed that $swap(\Delta^{-1}\bar{\bf e})+C$ is circ-neutral
for the constant mould $C$ defined in (0.4).  Therefore 
$swap\bigl(\bar{\be}(\tau)\bigr)+C$ is circ-neutral for the same mould $C$.
Thus, by Theorem 1.3, for each $r\ge 2$ we have the Fay relation
$$\F\bigl(\bar{\be}'(\tau)\bigr)(u_1,\ldots,u_r)=-rc_r(u_2+\cdots+u_r),$$
which is exactly (3.1.1) for this particular constant mould $C$ This concludes
the proof.\hfill{$\diamondsuit$}
\vskip .4cm
\noindent {\bf \S 3.2. The Fay relations for the Lie-like elliptic
	associator}
	\vskip .3cm
	In this section, we use the results of \S 2 to explicitly deduce the Fay 
	relations satisfied by the {\it reduced Lie-like elliptic associator }
	$\bar{\ba}(\tau)$ from the simple family 
	satisfied by the Lie-like elliptic generating series $\bar{\be}(\tau)$
	given in Theorem 3.1 and the characteristics of the mould $T_{01}:=ma(t_{01})$.

	We begin by separating $t_{01}$ into the sum 
	$t_0+\hat t_{01}$ where $t_0=-a-{{1}\over{2}}[a,b]$ and $\hat{t}_{01}$ is 
	of minimal degree $2$ in $b$. Write the associated moulds 
	$T_{01}=T_0+\hat{T}_{01}$ where
	$T_0=ma(t_0)=-a-{{1}\over{2}}U_1$ and $\hat{T}_{01}=ma(\hat t_{01})$,
	which is given explicitly by 
	$$\hat{T}_{01}(u_1,\ldots,u_r)={{B_r}\over{r!}}\sum_{i=1}^r (-1)^{i-1}\Bigl({{r-1}\atop{i-1}}\Bigr)u_i\eqno(3.2.1)$$
	where $B_r$ is the $r$-th Bernoulli number. In particular,
	$\hat{T}_{01}(u_1,\ldots,u_r)=0$ for all even $r$, and also for $r=1$.
	Let $\hat{T}'_{01}=dar^{-1}\hat{T}_{01}$.
	\vskip .1cm
	We saw in (2.4.4) that
	$$\bar{\ba}(\tau)=Darit\bigl(\overline{\be}(\tau)\bigr)\cdot T_{01}.$$ 
	Since $T_{01}=T_0+\hat{T}_{01}$, we have
	$$\bar{\ba}(\tau)=Darit\bigl(\overline{\be}(\tau)\bigr)\cdot T_0+ 
	Darit\bigl(\overline{\be}(\tau)\bigr)\cdot \hat{T}_{01}.\eqno(3.2.2)$$ 
	But $T_0=-a-(1/2)U_1$, and by the properties of the $Darit$-derivation
	given in (2.4.1), we thus have
	$$Darit\bigl(\overline{\be}(\tau)\bigr)\cdot T_0=-\overline{\be}(\tau),$$
	so by (3.2.2), we have
	$$\bar{\ba}(\tau)=-\bar{\be}(\tau)+Darit\bigl(\bar{\be}(\tau)\bigr)\cdot \hat{T}_{01}.\eqno(3.2.3)$$
	Thus the Fay relations satisfied by $\bar{\ba}(\tau)$ come from 
	the Fay relations satisfied by $\overline{\be}(\tau)$ on the one hand, 
	and those coming from $Darit\bigl(\bar{\be}(\tau)\bigr)\cdot \hat{T}_{01}$ 
	on the other.  
	The Fay relations satisfied by $\overline{\be}(\tau)$ are given in 
	Theorem 3.1, and Theorem 3.2 below gives a formula for the Fay relations of a 
	mould of the type $P=Darit(N)\cdot R$ under certain hypothesis on $N$ and
	$R$ that do hold when $N=\bar{\be}(\tau)$ and $R=\hat{T}_{01}$.
	In Lemma 3.3, we show that the Fay correction mould for $\hat{T}_{01}$
is zero. Based on these results, we can compute the explicit Fay correction 
for $\bar{\ba}(\tau)$; 
it is given in Theorem 3.4.
\vskip .3cm
\noindent {\bf Theorem 3.2.} {\it Let $M$ be an alternal and push-invariant
mould such that $swapM+C$ is circ-neutral for some constant-valued
mould $C=(c_r)_{r\ge 2}$. Set $N=\Delta M$.
Let $R\in ARI$ be even in depth 1 and satisfy a family of
Fay relations with correction mould $C_{R'}$, i.e.
$$\F(R')(u_1,\ldots,u_r)=C_{R'}(u_1,\ldots,u_r).$$
Then the mould $P=Darit(N)\cdot R$ satisfies the Fay relations 
$$\F(P')(u_1,\ldots,u_r)=C_{P'}(u_1,\ldots,u_r)$$
with correction mould $C_{P'}$ given for $r\ge 2$ by
$$\eqalign{C_{P'}(u_1,\ldots,u_r)
&=\sum_{1<i<j\le r} C_{R'}(u_1,\ldots,u_{i-1},u_i+\cdots+u_j,u_{j+1},\ldots,u_r)M(u_{i+1},\ldots,u_j)\cr
&-\sum_{1<i<j\le r} C_{R'}(u_1,\ldots,u_i,u_{i+1}+\cdots+u_{j+1},u_{j+2},\ldots,u_r)M(u_{i+1},\ldots,u_j)\cr
&-\sum_{i=2}^{r-1} ic_iR'(u_{i+1},\cdots,u_r)
+\sum_{i=2}^{r-1} ic_iR'(u_2,\ldots,u_{r-i+1}).}\eqno(3.2.4)$$
In particular, in the case where the mould $R$ satisfies the strict Fay
relations, i.e.~$C_{R'}=0$, the mould $P=Darit(N)\cdot R$ satisfies the 
Fay relations
$$C_{P'}(u_1,\ldots,u_r)=\sum_{i=2}^{r-1} ic_i\Bigl(R'(u_2,u_3,\ldots,u_{r-i+1})
		-R'(u_{i+1},u_{i+2},\ldots,u_r)\Bigr).\eqno(3.2.5)$$
Thus if $R$ satisfies the strict Fay relations and furthermore $swapM$ is
strictly circ-neutral, then $P=Darit(N)\cdot R$ satisfies the strict relations:
$$\F(P')(u_1,\ldots,u_r)=0.\eqno(3.2.6)$$
}
\noindent {\bf Proof.} The proof, which is long and technical, can be found
in Appendix A.\hfill{$\diamondsuit$}
\vskip .3cm
Let $N=\bar{\be}(\tau)$ and $M=\Delta^{-1}\bar{\be}(\tau)$. It was shown in the proof of Theorem 3.1 that this mould is alternal, push-invariant
and $swap(M)+C$ is circ-neutral for the mould $C$ of(0,4).  Thus 
$M=\Delta^{-1}\bar{\be}(\tau)$ satisfies the properties required for $M$ in Theorem 3.2. In order to apply 
the result of Theorem 3.2 to compute the Fay correction for
$Darit\bigl(\bar{\be}(\tau)\bigr)\cdot \hat{T}_{01}$, we also need to know
the Fay correction for the mould $R=\hat{T}_{01}$. The following lemma shows that
in fact $\hat{T}_{01}$ satisfies the strict Fay relations.  Its proof,
which is annoyingly technical, is banished to Appendix B.
\vskip .3cm
\noindent {\bf Lemma 3.3.} {\it The mould $\hat{T}_{01}$ satisfies the 
strict Fay relations
$$\F\bigl(\hat{T}'_{01}\bigr)(u_1,\ldots,u_r)=0,\ \ \ r\ge 2.$$}

\noindent We will now apply Theorems 3.1 and 3.2 and Lemma 3.3 directly to give the explicit Fay 
correction for $\bar{\ba}(\tau)$.
\vskip .3cm
\noindent {\bf Theorem 3.4.} {\it The Lie-like elliptic associator
$\bar{\ba}(\tau)$ satisfies the Fay relations
$$\F\bigl(\bar{\ba}'(\tau)\bigr)=C_{\bar{\ba}'(\tau)}(u_1,\ldots,u_r)$$
for $r\ge 2$, where the correction mould $C_{\bar{\ba}'(\tau)}$ is given 
for even $r$ by
$$C_{\bar{\ba}'(\tau)}(u_1,\ldots,u_r)=0$$
and for odd $r$ by
$$C_{\bar{\ba}'(\tau)}(u_1,\ldots,u_r)=
\zeta(r)(u_2+\cdots+u_r)
+\sum_{{{i=3}\atop{i\ {\rm odd}}}}^{r-1} \zeta(i)\Bigl(\hat{T}'_{01}(u_2,u_3,\ldots,u_{r-i+1})-\hat{T}'_{01}(u_{i+1},u_{i+2},\ldots,u_r)\Bigr).\eqno(3.2.7)$$}
									\vskip .2cm
\noindent {\bf Proof.} By (3.2.3), the Fay correction satisfied
									by $\bar{\ba}(\tau)$ is the sum of that of $-\overline{\be}(\tau)$ and that of
									$Darit\bigl(\overline{\be}(\tau)\bigr)\cdot \hat{T}_{01}$.  For even $r$,
									the Fay correction for $-\bar{\be}(\tau)$ is zero by Theorem 3.1.
									Since $\hat{T}_{01}$ satisfies the strict Fay relations by Lemma 3.3, the
Fay correction for $Darit\bigl(\bar{\be}(\tau)\bigr)\cdot \hat{T}_{01}$
is given by
(3.2.5).  But in the case where $r$ is even, the terms in the sum in (3.2.5)
	where $i$ is even are all zero because $c_i=0$, and the terms where
	$i$ is odd are also all zero, because
	the terms $R'(u_2,\ldots,u_{r-i+1})$ and $R'(u_{i+1},\ldots,u_r)$ are
	both of odd depth $r-i$, and the mould $R=\hat{T}_{01}$ is zero in all
	odd depths because of the Bernoulli coefficient in its definition (3.2.1).

	When $r$ is odd, the Fay correction term coming from $-\bar{\be}(\tau)$ is
	given in Theorem 3.1, and the Fay correction term coming from
	$Darit\bigl(\bar{\be}(\tau)\bigr)\cdot \hat{T}_{01}$ comes directly from
	(3.2.5), giving the desired formula (3.2.7) for the total Fay
correction. This concludes the proof.\hfill{$\diamondsuit$}
	\vskip .3cm
	Up to depth 5, the Fay correction for $\bar{\ba}(\tau)$ is given by
	$$\cases{\F\bigl(\bar{\ba}'(\tau)\bigr)(u_1,u_2)=0\cr
		\F\bigl(\bar{\ba}'(\tau)\bigr)(u_1,u_2,u_3)=\zeta(3)(u_2+u_3)\cr
			\F\bigl(\bar{\ba}'(\tau)\bigr)(u_1,u_2,u_3,u_4)=0\cr
			\F\bigl(\bar{\ba}'(\tau)\bigr)(u_1,u_2,u_3,u_4,u_5)=\zeta(5)(u_2+u_3+u_4+u_5)
			+\zeta(3)\Bigl(\displaystyle{{{u_2-u_3}\over{12u_2u_3}}-{{u_4-u_5}\over{12u_4u_5}}}\Bigr).}\eqno(3.2.8)$$
			\vskip .4cm
			\noindent {\bf \S 3.3. The group-like elliptic generating series and elliptic associator.}
			\vskip .3cm
			In this section we briefly recall how to ``undo'' the linearization
			of the moulds $\bar{\be}(\tau)$ and the double linearization of the mould $\bar{\ba}(\tau)$, so as to
			explicitly compute the Fay corrections of the {\it reduced group-like elliptic
			generating series} $\bar{E}(\tau)$ and the {\it reduced group-like elliptic associator} $\bar{A}(\tau)$.
			We underline the
			point that the natural place for simple closed formulas is in the Lie-like
			situation; the formulas in the group-like situation are deduced from these
			directly, for which reason we do not attempt to give a closed formula, but simply show the process.
			\vskip .3cm
In fact, we saw in (3.2.6) of Theorem 3.2 that if in fact $swap(\Delta^{-1}N)$ is 
strictly circ-neutral and $R$ satisfies the strict Fay relations, 
then $P=Darit(N)\cdot R$ also satisfies the strict Fay relations. Thus
also 
$${\rm exp}\bigl(Darit(N)\bigr)\cdot R=R+Darit(N)\cdot R+{{1}\over{2}}
Darit(N)^2\cdot R+\cdots$$
satisfies the strict Fay relations, so in the case of moulds satisfying
the strict relations, the closed formula is immediate.
It is only the correction moulds involved in the construction of
$\bar{\be}(\tau)$ and $\bar{\ba}(\tau)$ which make
a step-by-step 
calculation of the correction moulds in the group-like case necessary.

			\vskip .3cm
			Let us start with the group-like elliptic
generating series. Recall from (2.3.3) that
$\bar{E}_\tau=\Psi(a)=exp(\psi)(a)$. Thus in mould terms, since
we saw in (2.4.2) that the derivation $\psi$ corresponds to the
derivation $Darit\bigl(\bar{\be}(\tau)\bigr)$, we have
$$\bar{E}(\tau)=exp\Bigl(Darit\bigl(\bar{\be}(\tau)\bigr)\Bigr)\cdot a=\sum_{n\ge 0} {{1}\over{n!}}Darit^n\bigl(\bar{\be}(\tau)\bigr)\cdot a$$
which since $Darit\bigl(\bar{\be}(\tau)\bigr)\cdot a=\bar{\be}(\tau)$ by 
(2.4.1), we rewrite as
$$\bar{E}(\tau)=a+\bar{\be}(\tau)+\sum_{n\ge 2} {{1}\over{n!}}Darit^{n-1}\bigl(\bar{\be}(\tau)\bigr)\cdot \bar{\be}(\tau).\eqno(4.2.1)$$

			The depth $0$ term $a$ does not contribute to
the Fay correction mould for $\bar{E}(\tau)$.  Theorem 3.1 gives the
contribution of the term $\bar{\be}(\tau)$. 
Now set $N=R=\bar{\be}(\tau)$.  Then the Fay correction
mould $C_{R'}$ is given in Theorem 3.1, and the
constant correction mould $C$ such that $swap\bigl(\Delta^{-1}\bar{\be}(\tau)\bigr)+C$ is circ-neutral is the mould defined in (0.4).  Therefore
(3.2.4) from Theorem 3.2 gives us the Fay correction mould for
the $n=2$ term $Darit\bigl(\bar{\be}(\tau)\bigr)\cdot \bar{\be}(\tau)$ of (4.2.1). Next we set $N=\bar{\be}(\tau)$ and
$R=Darit\bigl(\bar{\be}(\tau)\bigr)\cdot \bar{\be}(\tau)$; then
using the Fay correction just computed for this mould $R$, Theorem 3.2
gives the Fay correction mould for the $n=3$ term 
$Darit\bigl(\bar{\be}(\tau)\bigr)^2\cdot \bar{\be}(\tau)$.  
Proceeding step by step in this manner, we compute each depth
of the Fay correction mould.  Up to depth 4, it is given by
$$\cases{\F\bigl(\bar{E}'(\tau)\bigr)(u_1,u_2)=0\cr
\F\bigl(\bar{E}'(\tau)\bigr)(u_1,u_2,u_3)=0\cr
\F\bigl(\bar{E}'(\tau)\bigr)(u_1,u_2,u_3,u_4)=
2\zeta(3)\sum_{r\ge 3\ {\rm odd}} \zeta(r)\bigl(u_2^r-u_4^r\bigr).}\eqno(4.2.1)$$
\vskip .3cm
We now turn to the case of $\bar{A}(\tau)$.  Recall from (2.3.5) that
$\bar{A}_\tau=\Psi(e^{t_{01}})$.
					Since $\bar{\ba}(\tau)$ is doubly linearized with respect to
					$\bar{A}(\tau)$, our first step is to undo one level of linearization, by
					computing the Fay correction for the log of the elliptic associator, given by $\Psi(t_{01})$, or in mould
terms by 
					$${\rm log}\bar{A}(\tau)={\rm exp}\Bigl(Darit\bigl(
								\overline{\be}(\tau)\bigr)\Bigr)\cdot T_{01}.$$

					Since by (3.2.3) we have
					$$\bar{\ba}(\tau)=-\bar{\be}(\tau)+Darit\bigl(\bar{\be}(\tau)\bigr)\cdot \hat{T}_{01},$$
					the log elliptic associator is given by
					$$\eqalign{{\rm log}\,\bar{A}(\tau)&=exp\bigl(Darit(\bar{\be}(\tau)\bigr)\cdot T_{01}\cr
							&=\sum_{n\ge 0}{{1}\over{n!}}Darit\bigl(\bar{\be}(\tau)\bigr)^n\cdot T_{01}\cr
						&=T_{01}+Darit\bigl(\bar{\be}(\tau)\bigr)\cdot T_{01}+\sum_{n\ge 2}{{1}\over{n!}}Darit\bigl(\bar{\be}(\tau)\bigr)^n\cdot T_{01}\cr
						&=T_{01}+\bar{\ba}(\tau)+\sum_{n\ge 2}{{1}\over{n!}}Darit\bigl(\bar{\be}(\tau)\bigr)^{n-1}\cdot \bar{\ba}(\tau)\cr
						&=T_{01}+\bar{\ba}(\tau)+\sum_{n\ge 2}{{1}\over{n!}}Darit\bigl(\bar{\be}(\tau)\bigr)^{n-1}\cdot \Bigl(-\bar{\be}(\tau)+Darit\bigl(\bar{\be}(\tau)\bigr)\cdot\hat{T}_{01}\Bigr)\cr
						&=T_{01}+\bar{\ba}(\tau)
						-\sum_{n\ge 2}{{1}\over{n!}}Darit\bigl(\bar{\be}(\tau)\bigr)^{n-1}\cdot \bar{\be}(\tau)
						+\sum_{n\ge 2}{{1}\over{n!}}Darit\bigl(\bar{\be}(\tau)\bigr)^n\cdot 
						\hat{T}_{01}.  }$$
						We already know how to determine the Fay corrections of $T_{01}$ (which is
								zero) and $\bar{\ba}(\tau)$ (which is given in Theorem 3.4). The third term
is exactly the term calculated in the preceding part for the 
elliptic generating series.  The fourth term is computed in the exact
same way, using Theorem 3.2 repeatedly.

Up to depth 4, we find that the log elliptic associator satisfies
						the Fay relations up to depth 4 given by (3.4) and (3.5):
$$\cases{\F\bigl({\rm log}\,\bar{A}'(\tau)\bigr)(u_1,u_2)=0\cr
\F\bigl({\rm log}\,\bar{A}'(\tau)\bigr)(u_1,u_2,u_3)=\zeta(3)(u_2+u_3)\cr
\F\bigl({\rm log}\,\bar{A}'(\tau)\bigr)(u_1,u_2,u_3,u_4)=
-\zeta(3)\sum_{r\ge 3\ {\rm odd}} \zeta(r)(u_2^r-u_4^r)}$$
where we note that since $dar$ is an automorphism, we have
$$\bigl({\rm log}\,\bar{A}(\tau)\bigr)'={\rm log}\,\bar{A}'(\tau).$$
\vskip .3cm
We now undo the second level of linearization by computing the Fay 
correction for $\bar{A}(\tau)$ as the exponential of ${\rm log}\bar{A}
(\tau)$.  We show how to compute the Fay relations for ${\rm exp}({\cal P})$ where ${\cal P}\in ARI$ is any mould satisfying a family of Fay relations
with correction mould $C_{{\cal P}'}$. Just as for the log, since $dar$ is an automorphism we have 
$$exp(P)'=dar^{-1}
exp({\cal P})=exp({\cal P}')=1+{\cal P}'+{{1}\over{2}}({\cal P}')^2
									+{{1}\over{6}}({\cal P}')^3+\cdots$$
									and furthermore
									$$dar^{-1}({\cal P}^r)=({\cal P}')^r.$$
									Thus 
									$$C_{exp({\cal P})'}(u_1,\ldots,u_r)=\F\bigl(exp({\cal P}')\bigr)(u_1,\ldots,u_r)=
									C_{{\cal P}'}(u_1,\ldots,u_r)+\sum_{r\ge 2} {{1}\over{r!}}
									\F\bigl(({\cal P}')^r\bigr)(u_1,\ldots,u_r).$$
									Then we write ${\cal P}'=log\bigl(exp({\cal P}')\bigr)$ to get the
									final answer in terms of $exp({\cal P}')$.
									Since ${\cal P}'(u_1)=exp({\cal P}')(u_1)$ in depth 1, we find that
									in depth 2,
									$$\eqalign{C_{exp({\cal P}')}(u_1,u_2)
										&=C_{{\cal P}'}(u_1,u_2)+{{1}\over{2}}exp({\cal P}')(u_1)exp({\cal P}')(u_2)\cr
											&\qquad +{{1}\over{2}}exp({\cal P}')(-u_1)exp({\cal P}')(u_1+u_2)
											+{{1}\over{2}}exp({\cal P}')(u_2)exp({\cal P}')(-u_1-u_2).}$$
											and indeed the general formula for the Fay correction of 
											$exp({\cal P})$ is given by
											$$C_{{\cal P}'}(u_1,\ldots,u_r)=\F(dar{\cal P}'')(u_1,\ldots,u_r),\eqno(3.6)$$
											where the mould ${\cal P}''$ is defined by
											$${\cal P}''(w)=
											\sum_{n=2}^r {{1}\over{n!}}\sum_{w=a_1\ldots a_n}\prod_{j=1}^n{\cal P}'(a_j)$$
											with $w=(u_1,\ldots,u_r)$.

Thus the Fay correction mould for $\bar{A}(\tau)$ can be
computed directly from that of ${\rm log}\,\bar{A}(\tau)$ with a 
closed formula.  In depth 2, for example, we obtain the expression:
$$\F(\bar{A}'_\tau)(u_1,u_2)=
{{1}\over{2}}\bar{A}'_\tau(u_1)\bar{A}'_\tau(u_2)+
{{1}\over{2}}	\bar{A}'_\tau(-u_1)\bar{A}'_\tau(u_1+u_2)+
{{1}\over{2}}	\bar{A}'_\tau(u_2)\bar{A}'_\tau(-u_1-u_2).$$

\vfill\eject
\noindent {\bf Appendix A:  Proof of Theorem 3.2.}
\vskip .2cm
Let $N$, $M$, $C$, $R$ and $C_{R'}$ be as in the statement of the proposition,
so $M$ is bialternal with constant correction $C$ and $N=\Delta\,M$,
while $R$ satisfies the Fay relations $\F(R')=C_{R'}$.
Set $Q=R'=dar^{-1}R$; since $R$ is assumed even in depth 1, the mould $Q$ is
odd in depth 1. We will use the identity 
$$P=Darit(N)\cdot R=dar\cdot arat(M)\cdot dar^{-1}R=
dar\cdot arat(M)\cdot Q,$$
where for each $M\in ARI$, $arat(M)$ is a derivation\footnote{$^*$}{In the more 
familiar notation used by \'Ecalle and in [S1], we have
 $arat(M)\cdot Q=-arit(M)\cdot Q+lu(M,Q)$.} of $ARI$ given explicitly by
$$\bigl(arat(M)\cdot Q\bigr)(u_1,\ldots,u_r)=
\sum_{{{w=abc}\atop{ac,b\ne \emptyset}}} 
\Bigl(Q(a\rceil c)M(b) -Q(a\lceil c)M(b)\Bigr),$$
where $w=(u_1,\ldots,u_r)$, and the ``flexion'' notation is defined as follows:
in any decomposition of $w=(u_1,\ldots,u_r)$ into words $a_1\cdots a_m$ (for example
$w=abc$) with $a_l=(u_i,\ldots,u_j)$, $a_{l+1}=(u_{j+1},\ldots,u_k)$, we set
$$\cases{a_l\rceil=(u_i,\ldots,u_{j-1},u_j+u_{j+1}+\cdots+u_l)\cr
\lceil a_{l+1}=(u_i+\cdots+u_j+u_{j+1},u_{j+2},\ldots,u_k).}$$
Starting from (A.3) below, we will also use the flexion notation for moulds in
$\overline{ARI}$, which is defined as follows: 
in any decomposition of $w=(v_1,\ldots,v_r)$ into words $a_1\cdots a_m$ 
with $a_l=(v_i,\ldots,v_j)$, $a_{l+1}=(v_{j+1},\ldots,v_k)$, we set
$$\cases{a_l\rfloor=(v_i-v_{j+1},\ldots,v_j-v_{j+1})\cr
\lfloor a_{l+1}=(v_{j+1}-v_j,\ldots,v_k-v_j).}$$
\vskip .1cm
The goal of this proof is to compute the Fay correction mould
$C_{P'}$ for the Fay relations
$$C_{P'}:=\F(P')=\F\bigl(arat(M)\cdot Q\bigr).$$
To do this, we will resort to the swap version of the Fay relations given
in (2.3):
$$swap(P')\bigl(sh\bigl((v_1), (v_2,\ldots,v_r)\bigr)\bigr)
+push^{-1}\cdot swap(P')(v_1,\ldots,v_r)=push^{-1}\cdot swap(C_{P'})(v_1,\ldots,v_r).\eqno(A.1)$$
Let us simplify the terminology by setting
$$\cases{\M=swap(M)\cr
\Q=swap(Q)\cr
\P=swap(P')=swap(arat(M)\cdot Q)\cr
C_{\bf P}=push^{-1}\cdot swap(C_{P'}).}$$
Then (A.1) translates to
$$C_{\P}(v_1,\ldots,v_r):=\P\Bigl(sh\bigl((v_1),(v_2,\ldots,v_r)\bigr)\Bigr)+push^{-1}\P(v_1,\ldots,v_r).\eqno(A.2)$$
\vskip .3cm
Letting $w=(v_1,\ldots,v_r)$ and using the known commutation relations 
of $swap$ with $arit$ and with $lu$ (cf.~\S 2.4 of [S2])
together with the fact that $M$ is push-invariant, we find that
$$\P(w)=swap\Bigl(arat\bigl(M)\bigr)\cdot Q\Bigr)(w)
=\sum_{{{w=abc}\atop{a,b\ne\emptyset}}} \Q(ac)\M(\lfloor b)
-\sum_{{{w=abc}\atop{ac,b\ne\emptyset}}} \Q(ac)\M(b\rfloor)
+\sum_{{{w=ab}\atop{a,b\ne \emptyset}}} \Q(a\rfloor)\M(b).\eqno(A.3)$$
We can now proceed to the computation of the Fay correction mould $C_{P'}$
via the computation of $C_{\bf P}$, for which we compute and sum the two terms
from the right-hand side of (A.3).
\vskip .3cm
We do this by adapting a method due to A.~Salerno (cf. [SS, appendix]), 
consisting in separating the terms running the $v_1$ through $a$, $b$ and $c$
as follows.  We write $w'=v_2\cdots v_r$, so $w=v_1w'$.
\vskip .2cm
$$ \P\Bigl(sh\bigl((v_1),(v_2,\ldots,v_r)\bigr)\Bigr)=
\sum_{{{w'=abc}\atop{b\ne\emptyset}}} \Q\bigl(sh(v_1,a)c\bigr)\M(\lfloor b)
+\sum_{{{w'=abc}\atop{a\ne\emptyset}}} \Q(ac)\M\bigl(\lfloor sh(v_1,b)\bigr)$$
$$+\sum_{{{w'=abc}\atop{a,b\ne\emptyset}}} \Q\bigl(a\,sh(v_1,c)\bigr)\M(\lfloor b)
-\sum_{{{w'=abc}\atop{b\ne\emptyset}}} \Q\bigl(sh(v_1,a)c\bigr)\M(b\rfloor)
-\sum_{{{w'=abc}\atop{ac\ne\emptyset}}} \Q(ac)\M\bigl(sh(v_1,b)\rfloor\bigr)$$
$$-\sum_{{{w'=abc}\atop{b\ne\emptyset}}} \Q\bigl(a\,sh(v_1,c)\bigr)\M(b\rfloor)
+\sum_{{{w'=ab}\atop{b\ne\emptyset}}} \Q(sh(v_1,a)\rfloor)\M(b)
+\sum_{{{w'=ab}\atop{a\ne\emptyset}}} \Q(a\rfloor)\M(sh(v_1,b)).\eqno(A.4)$$
We split up the second and fifth sums of (A.4) as follows:
$$\sum_{{{w'=abc}\atop{a\ne\emptyset}}} \Q(ac)\M\bigl(\lfloor sh(v_1,b)\bigr)
-\sum_{{{w'=abc}\atop{ac\ne\emptyset}}} \Q(ac)\M\bigl(sh(v_1,b)\rfloor\bigr)=$$
$$\sum_{{{w'=abc}\atop{a\ne\emptyset}}} \Q(ac)\M\bigl(\lfloor sh(v_1,b)\bigr)
-\sum_{{{w'=abc}\atop{a\ne\emptyset}}} \Q(ac)\M\bigl(sh(v_1,b)\rfloor\bigr)
-\sum_{{{w'=bc}\atop{b,c\ne\emptyset}}} \Q(c)\M\bigl(sh(v_1,b)\rfloor\bigr)
-\Q(w')\M(v_1-v_2).$$ 
Since $M+C$ is alternal, where $C=(c_i)_{i\ge 0}$ is the constant correction mould for $M$, we can write
$$\M\bigl(sh(v_1,b)\bigr)=\M\bigl(\lfloor sh(v_1,b)\bigr)=\M\bigl(sh(v_1,b)\rfloor\bigr)=-(|b|+1)c_{|b|+1}$$
for $|b|>1$, so we see that all the terms for $b\ne \emptyset$ in the 
first two sums cancel out, and it becomes equal to
$$\sum_{{{w'=ac}\atop{a\ne\emptyset}}} \Q(ac)\M(v_1-a_l)
-\sum_{{{w'=ac}\atop{a\ne\emptyset}}} \Q(ac)\M(v_1-c_f)
+\sum_{{{w'=bc}\atop{b,c\ne\emptyset}}} \Q(c)(|b|+1)c_{|b|+1}
-\Q(w')\M(v_1-v_2).$$ 

Writing $a=(v_2,\ldots,v_i)$ and $c=(v_{i+1},\ldots,v_r)$,
and setting $v_{r+1}=0$, this can be rewritten as
$$\sum_{i=2}^{r} \Q(v_2,\ldots,v_r)\M(v_1-v_i)
-\sum_{i=2}^r \Q(v_2,\ldots,v_r)\M(v_1-v_{i+1})
-\Q(v_2,\ldots,v_r)\M(v_1-v_2)
+\sum_{{{w'=bc}\atop{b,c\ne\emptyset}}} (|b|+1)c_{|b|+1}\Q(c),$$
in which most terms cancel out, 
so that finally the second and fifth term in (A.4) sum to just
$$-\Q(v_2,\ldots,v_r)\M(v_1)
+\sum_{{{w'=bc}\atop{b,c\ne\emptyset}}} (|b|+1)c_{|b|+1}\Q(c).$$
Replacing this in (A.4), we obtain

$$\P\Bigl(sh\bigl((v_1),(v_2,\ldots,v_r)\bigr)\Bigr)=
\sum_{{{w'=abc}\atop{b\ne\emptyset}}} \Q\bigl(sh(v_1,a)c\bigr)\M(\lfloor b)
+\sum_{{{w'=abc}\atop{a,b\ne\emptyset}}} \Q\bigl(a\,sh(v_1,c)\bigr)\M(\lfloor b)$$
$$-\sum_{{{w'=abc}\atop{b\ne\emptyset}}} \Q\bigl(sh(v_1,a)c\bigr)\M(b\rfloor)
-\sum_{{{w'=abc}\atop{b\ne\emptyset}}} \Q\bigl(a\,sh(v_1,c)\bigr)\M(b\rfloor)
+\sum_{{{w'=ab}\atop{b\ne\emptyset}}} \Q(sh(v_1,a)\rfloor)\M(b)$$
$$+\sum_{{{w'=ab}\atop{a\ne\emptyset}}} \Q(a\rfloor)\M(sh(v_1,b))
-\Q(v_2,\ldots,v_r)\M(v_1)
+\sum_{{{w'=bc}\atop{b,c\ne\emptyset}}} (|b|+1)c_{|b|+1}\Q(c).\eqno(A.5)$$

Let us now add and subtract the term $a=\emptyset$ to the second sum of
(A.5), in order to make the first four sums all over the same indices:
$$\sum_{{{w'=abc}\atop{a,b\ne\emptyset}}} \Q\bigl(a\,sh(v_1,c)\bigr)\M(\lfloor b)
= \sum_{{{w'=abc}\atop{b\ne\emptyset}}} \Q\bigl(a\,sh(v_1,c)\bigr)\M(\lfloor b)-\sum_{{{w'=bc}\atop{b\ne \emptyset}}} \Q\bigl(sh(v_1,c)\bigr)\M(b).$$
We will also modify the form of the sixth term in (A.5), noting that
because of the $a\rfloor$ term, the shuffle term in $M$ is not an independent
factor: for example, if $r=4$, in the decomposition $w'=ab=(v_2)(v_3,v_4)$,
the sixth term is 
$$\Q(v_2-v_1)\M(v_1,v_3,v_4)+\Q(v_2-v_3)\M(v_3,v_1,v_4)+\Q(v_2-v_3)\M(v_3,v_4,v_1),$$
since the term $a\rfloor$ is with respect to $v_1$ in the first shuffle 
term $(v_1,v_3,v_4)$, but with respect to $v_3$ for all the other shuffle
terms of $sh\bigl((v_1),(v_3,v_4)\bigr)$. To take advantage of the alternality
of $\M$, we rewrite the sixth term as
$$\sum_{{{w'=ab}\atop{a\ne\emptyset}}} \Q(a\rfloor)\M(sh(v_1,b))=
\Q(v_2-v_1,\ldots,v_r-v_1)\M(v_1)
+\sum_{{{w'=ab}\atop{a,b\ne\emptyset}}} 
\Q(a-b_1)\M(sh(v_1,b))$$
$$\qquad\qquad\qquad\qquad +\sum_{{{w'=ab}\atop{a,b\ne\emptyset}}} 
\Q(a-v_1)M(v_1b)
-\sum_{{{w'=ab}\atop{a,b\ne\emptyset}}} 
\Q(a-b_1)\M(v_1b)$$
$$=\Q(v_2-v_1,\ldots,v_r-v_1)\M(v_1)
+\sum_{{{w'=ab}\atop{a,b\ne\emptyset}}} 
\Q(a-v_1)M(v_1b)
-\sum_{{{w'=ab}\atop{a,b\ne\emptyset}}} 
\Q(a-b_1)\M(v_1b)
-\sum_{{{w'=ab}\atop{a,b\ne\emptyset}}} 
\Q(a-b_1)(|b|+1)c_{|b|+1}$$
Plugging these two modifications into (A.5), we obtain
$$\P\Bigl(sh\bigl((v_1),(v_2,\ldots,v_r)\bigr)\Bigr)=
\sum_{{{w'=abc}\atop{b\ne\emptyset}}} \Q\bigl(sh(v_1,a)c\bigr)\M(\lfloor b)
+\sum_{{{w'=abc}\atop{b\ne\emptyset}}} \Q\bigl(a\,sh(v_1,c)\bigr)\M(\lfloor b)
-\sum_{{{w'=abc}\atop{b\ne\emptyset}}} \Q\bigl(sh(v_1,a)c\bigr)\M(b\rfloor)$$
$$-\sum_{{{w'=abc}\atop{b\ne\emptyset}}} \Q\bigl(a\,sh(v_1,c)\bigr)\M(b\rfloor)
-\sum_{{{w'=bc}\atop{b\ne\emptyset}}}\Q\bigl(sh(v_1,c)\bigr)\M(b)
+\sum_{{{w'=ab}\atop{b\ne\emptyset}}} \Q(sh(v_1,a)\rfloor)\M(b)$$
$$+\Q(v_2-v_1,\ldots,v_r-v_1)\M(v_1)
+\sum_{{{w'=ab}\atop{a,b\ne\emptyset}}} \Q(a-v_1)M(v_1b)
-\sum_{{{w'=ab}\atop{a,b\ne\emptyset}}} \Q(a-b_1)\M(v_1b)$$
$$-\sum_{{{w'=ab}\atop{a,b\ne\emptyset}}} \Q(a-b_1)(|b|+1)c_{|b|+1}
-\Q(v_2,\ldots,v_r)\M(v_1)
+\sum_{{{w'=bc}\atop{b,c\ne\emptyset}}} (|b|+1)c_{|b|+1}\Q(c).\eqno(A.6)$$
The seventh term is the $b=\emptyset$ case of the eighth term and the
eleventh term is the $b=\emptyset$ case of the ninth term, so we can
simplify this as
$$\P\Bigl(sh\bigl((v_1),(v_2,\ldots,v_r)\bigr)\Bigr)=
\sum_{{{w'=abc}\atop{b\ne\emptyset}}} \Q\bigl(sh(v_1,a)c\bigr)\M(\lfloor b)
+\sum_{{{w'=abc}\atop{b\ne\emptyset}}} \Q\bigl(a\,sh(v_1,c)\bigr)\M(\lfloor b)$$
$$-\sum_{{{w'=abc}\atop{b\ne\emptyset}}} \Q\bigl(sh(v_1,a)c\bigr)\M(b\rfloor)
-\sum_{{{w'=abc}\atop{b\ne\emptyset}}} \Q\bigl(a\,sh(v_1,c)\bigr)\M(b\rfloor)
-\sum_{{{w'=bc}\atop{b\ne\emptyset}}}\Q\bigl(sh(v_1,c)\bigr)\M(b)$$
$$+\sum_{{{w'=ab}\atop{b\ne\emptyset}}} \Q(sh(v_1,a)\rfloor)\M(b)
+\sum_{{{w'=ab}\atop{a\ne\emptyset}}} \Q(a-v_1)M(v_1b)
-\sum_{{{w'=ab}\atop{a\ne\emptyset}}} \Q(a-b_1)\M(v_1b)$$
$$-\sum_{{{w'=ab}\atop{a,b\ne\emptyset}}} \Q(a-b_1)(|b|+1)c_{|b|+1}
+\sum_{{{w'=bc}\atop{b,c\ne\emptyset}}} (|b|+1)c_{|b|+1}\Q(c).\eqno(A.7)$$

We will now recompute the sum of the first four terms in (A.7), with
the decompositions
$$w'=abc=(v_1,\ldots,v_i)(v_{i+1},\ldots,v_j)(v_{j+1},\ldots,v_r).$$
Set $a=a'v_i$ and $c=v_{j+1}c'$,  
and observe that in the first term $\Q\bigl(sh(v_1,a)c\bigr)\M(\lfloor b)$,
the flexion $\lfloor b$ is given by 
$$\lfloor b=\cases{(v_{i+1}-v_i,\ldots,v_j-v_i)&for every term in
the shuffle $sh(v_1,a')v_i$\cr 
(v_{i+1}-v_1,\ldots,v_j-v_1)&for the term $av_1$.}$$ In the second term
$\Q\bigl(a\,sh(v_1,c)\bigr)\M(\lfloor b)$, however, the left flexion is always
$\lfloor b=(v_{i+1}-v_i,\ldots,v_j-v_i)$.  Similarly, in the third term
$\Q\bigl(sh(v_1,a)c)\M(b\rfloor)$ the right flexion is always
$b\rfloor=(v_{i+1}-v_{j+1},\ldots,v_j-v_{j+1})$, but in the fourth term
$\Q\bigl(a\,sh(v_1,c)\bigr)\M(b\rfloor)$, it is given by
$$b\rfloor=\cases{(v_{i+1}-v_{j+1},\ldots,v_j-v_{j+1})&for every term in 
$v_{j+1}sh(v_1,c')$\cr
(v_{i+1}-v_1,\ldots,v_j-v_1)&for $v_1c$.}$$
Observing that writing $a=a'v_i$ and $c=v_{j+1}c'$, we have
$$\Q\bigl(sh(v_1,a)c\bigr)=\Q\bigl(sh(v_1,a')v_ic)+\Q(av_1c).$$
and
$$\Q\bigl(a\,sh(v_1,c)\bigr)=\Q(av_1c)+\Q\bigl(av_{j+1}sh(v_1,c')\bigr),$$
we recombine the first four terms in (A.7) as
$$\eqalign{\P\Bigl(sh\bigl((v_1),(v_2,\ldots,v_r)\bigr)\Bigr)
&=\sum_{{{w'=abc}\atop{b\ne\emptyset}}} \Q\bigl(sh(v_1,ac)\bigr)\M(b-a_l)
+\sum_{{{w'=abc}\atop{b\ne\emptyset}}} \Q\bigl(av_1c\bigr)\M(b-v_1)\cr
&-\sum_{{{w'=abc}\atop{b\ne\emptyset}}} \Q\bigl(sh(v_1,ac)\bigr)\M(b-c_f)
-\sum_{{{w'=abc}\atop{b\ne\emptyset}}} \Q\bigl(av_1c\bigr)\M(b-v_1)\cr
&-\sum_{{{w'=bc}\atop{b\ne\emptyset}}}\Q\bigl(sh(v_1,c)\bigr)\M(b)
+\sum_{{{w'=ab}\atop{b\ne\emptyset}}} \Q(sh(v_1,a)\rfloor)\M(b)\cr
&+\sum_{{{w'=ab}\atop{a\ne\emptyset}}} \Q(a-v_1)M(v_1b)
-\sum_{{{w'=ab}\atop{a\ne\emptyset}}} \Q(a-b_1)\M(v_1b)\cr
&-\sum_{{{w'=ab}\atop{a,b\ne\emptyset}}} \Q(a-b_1)(|b|+1)c_{|b|+1}
+\sum_{{{w'=bc}\atop{b,c\ne\emptyset}}} (|b|+1)c_{|b|+1}\Q(c)}\eqno(A.8)$$
where we use the notation 
$$w'=abc=(v_1,\ldots,v_i)(v_{i+1},\ldots,v_j)(v_{j+1},\ldots,v_r)$$
and
$$a_l=\cases{a_{last}=v_i&if $a\ne\emptyset$\cr
a_l=0&if $a=\emptyset$},\qquad
c_f=\cases{c_{first}=v_{j+1}&if $c\ne\emptyset$\cr
c_f=0&if $c=\emptyset$.}$$
The second and fourth terms in (A.8) cancel out, leaving
$$\eqalign{\P\Bigl(sh\bigl((v_1),(v_2,\ldots,v_r)\bigr)\Bigr)
&=\sum_{{{w'=abc}\atop{b,ac\ne\emptyset}}} \Q\bigl(sh(v_1,ac)\bigr)\M(b-a_l)
-\sum_{{{w'=abc}\atop{b,ac\ne\emptyset}}} \Q\bigl(sh(v_1,ac)\bigr)\M(b-c_f)\cr
&-\sum_{{{w'=bc}\atop{b\ne\emptyset}}}\Q\bigl(sh(v_1,c)\bigr)\M(b)
+\sum_{{{w'=ab}\atop{b\ne\emptyset}}} \Q(sh(v_1,a)\rfloor)\M(b)\cr
&+\sum_{{{w'=ab}\atop{a\ne\emptyset}}} \Q(a-v_1)M(v_1b)
-\sum_{{{w'=ab}\atop{a\ne\emptyset}}} \Q(a-b_1)\M(v_1b)\cr
&-\sum_{{{w'=ab}\atop{a,b\ne\emptyset}}} \Q(a-b_1)(|b|+1)c_{|b|+1}
+\sum_{{{w'=bc}\atop{b,c\ne\emptyset}}} (|b|+1)c_{|b|+1}\Q(c)}\eqno(A.9)$$
(where we note that the terms with $ac=\emptyset$ cancel out in the first 
two sums, so we can assume $ac\ne\emptyset$ in these sums).

We will now finish our calculation by using the swapped version of the
Fay relation satisfied by $\Q$:
$$\Q(sh(a,b))=-push^{-1}\Q(ab)+C_\Q(ab)$$
for any non-empty sequences $a$, $b$, to simplify the first four
terms of (A.9).

$$\eqalign{
&\sum_{{{w'=abc}\atop{b\ne\emptyset}}} \Q\bigl(sh(v_1,ac)\bigr)\M(b-a_l)
-\sum_{{{w'=abc}\atop{b\ne\emptyset}}} \Q\bigl(sh(v_1,ac)\bigr)\M(b-c_f)
=\sum_{{{w'=abc}\atop{ac,b\ne \emptyset}}} push^{-1}\Q(v_1ac)\M(b-c_f)\cr
&-\sum_{{{w'=abc}\atop{ac,b\ne \emptyset}}} push^{-1}\Q(v_1ac)\M(b-a_l)
+\sum_{{{w'=abc}\atop{ac,b\ne \emptyset}}} C_\Q(v_1ac)\M(b-a_l)
-\sum_{{{w'=abc}\atop{ac,b\ne \emptyset}}} C_\Q(v_1ac)\M(b-c_f).}$$
Furthermore, we use the swapped Fay to simplify the third term in (A.9) as
$$\eqalign{-\sum_{{{w'=bc}\atop{b\ne\emptyset}}}\Q\bigl(sh(v_1,c)\bigr)\M(b)
&=-\sum_{{{w'=bc}\atop{b,c\ne\emptyset}}}\Q\bigl(sh(v_1,c)\bigr)\M(b)
-\Q(v_1)\M(v_2,\ldots,v_r)\cr
&=\sum_{{{w'=bc}\atop{b\ne\emptyset}}}push^{-1}\Q(v_1c)\M(b)
-\sum_{{{w'=bc}\atop{b\ne\emptyset}}}C_\Q(v_1c)\M(b)
-\Q(v_1)\M(v_2,\ldots,v_r),}$$
and the fourth term in (A.9) as
$$\eqalign{\sum_{{{w'=ab}\atop{b\ne\emptyset}}} \Q(sh(v_1,a)\rfloor)\M(b)
&=\sum_{{{w'=ab}\atop{a,b\ne\emptyset}}} \Q(sh(v_1,a)\rfloor)\M(b)
+\Q(v_1-v_2)\M(v_2,\ldots,v_r)\cr
&=-\sum_{{{w'=ab}\atop{a,b\ne\emptyset}}} push^{-1}\Q(v_1a\rfloor)\M(b)
+\sum_{{{w'=ab}\atop{b\ne\emptyset}}} C_\Q(v_1a\rfloor)\M(b)
+\Q(v_1-v_2)\M(v_2,\ldots,v_r).}$$
Replacing the first four terms in (A.9) in this way, we obtain
$$\eqalign{\P\Bigl(sh\bigl((v_1),(v_2,\ldots,v_r)\bigr)\Bigr)
&=\sum_{{{w'=abc}\atop{ac,b\ne \emptyset}}} push^{-1}\Q(v_1ac)\M(b-c_f)
-\sum_{{{w'=abc}\atop{ac,b\ne \emptyset}}} push^{-1}\Q(v_1ac)\M(b-a_l)\cr
&+\sum_{{{w'=bc}\atop{b,c\ne\emptyset}}}push^{-1}\Q(v_1c)\M(b)
-\Q(v_1)\M(v_2,\ldots,v_r)\cr
&-\sum_{{{w'=ab}\atop{a,b\ne\emptyset}}} push^{-1}\Q(v_1a\rfloor)\M(b)
+\Q(v_1-v_2)\M(v_2,\ldots,v_r)\cr
&+\sum_{{{w'=ab}\atop{a\ne\emptyset}}} \Q(a-v_1)M(v_1b)
-\sum_{{{w'=ab}\atop{a\ne\emptyset}}} \Q(a-b_1)\M(v_1b)\cr
&-\sum_{{{w'=bc}\atop{b\ne\emptyset}}}C_\Q(v_1c)\M(b)
+\sum_{{{w'=ab}\atop{b\ne\emptyset}}} C_\Q(v_1a\rfloor)\M(b)\cr
&+\sum_{{{w'=abc}\atop{ac,b\ne \emptyset}}} C_\Q(v_1ac)\M(b-a_l)
-\sum_{{{w'=abc}\atop{ac,b\ne \emptyset}}} C_\Q(v_1ac)\M(b-c_f)\cr
&-\sum_{{{w'=ab}\atop{a,b\ne\emptyset}}} \Q(a-b_1)(|b|+1)c_{|b|+1}
+\sum_{{{w'=bc}\atop{b,c\ne\emptyset}}} (|b|+1)c_{|b|+1}\Q(c)}\eqno(A.10)$$

Let us write (A.10) as
$$\P\Bigl(sh\bigl((v_1),(v_2,\ldots,v_r)\bigr)\Bigr)={\cal A}+{\cal C},$$
where ${\cal C}$ denotes the main part given by the first eight terms
and ${\cal C}$ denotes the correction part given in the last six terms.
In order to prove (A.2), we will show in (A.11)-(A.17) below that
$$push({\cal A})(v_1,\ldots,v_r)=-P(v_1,\ldots,v_r),\ \ {\rm i.e.}\ \ 
({\cal A})(v_1,\ldots,v_r)=-push^{-1}P(v_1,\ldots,v_r),$$
and thus ${\cal C}$ is the desired correction mould $C_{\P}$.

To apply the push-operator to ${\cal A}$, we first
break the first two terms of ${\cal A}$ into two sums.
$$\eqalign{{\cal A}(v_1,\ldots,v_r)&=\sum_{{{w'=abc}\atop{b,c\ne \emptyset}}} push^{-1}\Q(v_1ac)\M(b-c_f)
+\sum_{{{w'=ab}\atop{a,b\ne \emptyset}}} push^{-1}\Q(v_1a)\M(b)\cr
&-\sum_{{{w'=abc}\atop{a,b\ne \emptyset}}} push^{-1}\Q(v_1ac)\M(b-a_l)
-\sum_{{{w'=bc}\atop{b,c\ne \emptyset}}} push^{-1}\Q(v_1c)\M(b)\cr
&+\sum_{{{w'=bc}\atop{b,c\ne\emptyset}}}push^{-1}\Q(v_1c)\M(b)
-\Q(v_1)\M(v_2,\ldots,v_r)\cr
&-\sum_{{{w'=ab}\atop{a,b\ne\emptyset}}} push^{-1}\Q(v_1a\rfloor)\M(b)
+\Q(v_1-v_2)\M(v_2,\ldots,v_r)\cr
&+\sum_{{{w'=ab}\atop{a\ne\emptyset}}} \Q(a-v_1)M(v_1b)
-\sum_{{{w'=ab}\atop{a\ne\emptyset}}} \Q(a-b_1)\M(v_1b).}\eqno(A.11)$$
The fourth and fifth terms cancel.  Writing this out in indices with
$abc=(v_2,\ldots,v_i)(v_{i+1},\ldots,v_j)(v_{j+1},\ldots,v_r)$, we find
$$\eqalign{{\cal A}(v_1,\ldots,v_r)&=\sum_{1\le i<j<r} push^{-1}\Q(v_1,v_2,\ldots,v_i,v_{j+1},\ldots,v_r)
\M(v_{i+1}-v_{j+1},\ldots,v_j-v_{j+1})\cr
&+\sum_{1<i<r} push^{-1}\Q(v_1,v_2,\ldots,v_i)
\M(v_{i+1},\ldots,v_r)\cr
&-\sum_{1<i<j\le r} push^{-1}\Q(v_1,v_2,\ldots,v_i,v_{j+1},\ldots,v_r)
\M(v_{i+1}-v_i,\ldots,v_j-v_i)\cr
&-\sum_{1<i<r} push^{-1}\Q(v_1-v_{i+1},\ldots,v_i-v_{i+1})\M(v_{i+1},\ldots,v_r)\cr
&+\sum_{1<i\le r}\Q(v_2-v_1,\ldots,v_i-v_1)\M(v_1,v_{i+1},\ldots,v_r)\cr
&-\sum_{1<i\le r}\Q(v_2-v_{i+1},\ldots,v_i-v_{i+1})\M(v_1,v_{i+1},\ldots,v_r)\cr
&-\Q(v_1)\M(v_2,\ldots,v_r)+\Q(v_1-v_2)\M(v_2,\ldots,v_r).}\eqno(A.12)$$
Now we compute $push({\cal A})$ by applying the $push$ operator to (A.12): it maps
$v_1\mapsto -v_r$ and $v_k\mapsto v_{k-1}-v_r$ for $2\le k\le r$, so we obtain
$$\eqalign{push&({\cal A})(v_1,\ldots,v_r)=\sum_{1\le i<j<r} 
push^{-1}\Q(-v_r,v_1-v_r,\ldots,v_{i-1}-v_r,v_j-v_r,\ldots,v_{r-1}-v_r)
\M(v_i-v_j,\ldots,v_{j-1}-v_j)\cr
&+\sum_{1<i<r} 
push^{-1}\Q(-v_r,v_1-v_r,\ldots,v_{i-1}-v_r)
\M(v_i-v_r,\ldots,v_{r-1}-v_r)\cr
&-\sum_{1<i<j\le r} 
push^{-1}\Q(-v_r,v_1-v_r,\ldots,v_{i-1}-v_r,v_j-v_r,\ldots,v_{r-1}-v_r)
\M(v_i-v_{i-1},\ldots,v_{j-1}-v_{i-1})\cr
&-\sum_{1<i<r} 
push^{-1}\Q(-v_i,v_1-v_i,\ldots,v_{i-1}-v_i)\M(v_i-v_r,\ldots,v_{r-1}-v_r)\cr
&+\sum_{1<i\le r}
\Q(v_1,\ldots,v_{i-1})\M(-v_r,v_i-v_r,\ldots,v_{r-1}-v_r)\cr
&-\sum_{1<i\le r}
\Q(v_1-v_i,\ldots,v_{i-1}-v_i)\M(-v_r,v_i-v_r,\ldots,v_{r-1}-v_r)\cr
&-\Q(-v_r)\M(v_1-v_r,\ldots,v_{r-1}-v_r)
+\Q(-v_1)\M(v_1-v_r,\ldots,v_{r-1}-v_r).}\eqno(A.13)$$
Now we unpush the $\Q$s (recall that $Q$ is odd in depth 1),
and also some of the $\M$'s (given that $\M$ being bialternal is push-invariant):
$$\eqalign{push&({\cal A})(v_1,\ldots,v_r)=\sum_{1\le i<j<r} 
\Q(v_1,\ldots,v_{i-1},v_j,\ldots,v_r)
\M(v_i-v_j,\ldots,v_{j-1}-v_j)\cr
&+\sum_{1<i<r} 
\Q(v_1,\ldots,v_{i-1},v_r)
\M(v_i-v_r,\ldots,v_{r-1}-v_r)\cr
&-\sum_{1<i<j\le r} 
\Q(v_1,\ldots,v_{i-1},v_j,\ldots,v_r)
\M(v_i-v_{i-1},\ldots,v_{j-1}-v_{i-1})\cr
&-\sum_{1<i<r} 
\Q(v_1,\ldots,v_i)
\M(v_{i+1}-v_i,\ldots,v_r-v_i)\cr
&+\sum_{1<i\le r}
\Q(v_1,\ldots,v_{i-1})\M(v_i,\ldots,v_r)\cr
&-\sum_{1<i\le r}
\Q(v_1-v_i,\ldots,v_{i-1}-v_i)\M(v_i,\ldots,v_r)\cr
&+\Q(v_r)\M(v_1-v_r,\ldots,v_{r-1}-v_r)
-\Q(v_1)\M(v_2-v_1,\ldots,v_r-v_1).}\eqno(A.14)$$
Now we rewrite this in terms of $w=abc=(v_1,\ldots,v_{i-1})(v_i,\ldots,v_{j-1})
(v_j,\ldots,v_r)$:
$$\eqalign{push({\cal A})(v_1,\ldots,v_r)&=\sum_{{{w=abc}\atop{b\ne\emptyset,|c|\ge 2}}} \Q(ac)\M(b\rfloor)
+\sum_{{{w=abc}\atop{a,b\ne\emptyset,|c|=1}}} \Q(ac)\M(b\rfloor)
-\sum_{{{w=abc}\atop{a,b,c\ne \emptyset}}} \Q(ac)\M(\lfloor b)\cr
&-\sum_{{{w=ab}\atop{|a|\ge 2,b\ne\emptyset}}} \Q(a)M(\lfloor b)
+\sum_{{{w=ab}\atop{a,b\ne\emptyset}}} \Q(a)\M(b)
-\sum_{{{w=ab}\atop{a,b\ne\emptyset}}} \Q(a\rfloor)\M(b)\cr
&+\Q(v_r)\M(v_1-v_r,\ldots,v_{r-1}-v_r)
-\Q(v_1)\M(v_2-v_1,\ldots,v_r-v_1).}\eqno(A.15)$$
Term 7 is the $a=\emptyset$ case of term 2, and term 8 is the $|a|=1$ case of term 4:
$$\eqalign{push({\cal A})(v_1,\ldots,v_r)&=\sum_{{{w=abc}\atop{b\ne\emptyset,|c|\ge 2}}} \Q(ac)\M(b\rfloor)
+\sum_{{{w=abc}\atop{b\ne\emptyset,|c|=1}}} \Q(ac)\M(b\rfloor)
-\sum_{{{w=abc}\atop{a,b,c\ne \emptyset}}} \Q(ac)\M(\lfloor b)\cr
&-\sum_{{{w=ab}\atop{a,b\ne\emptyset}}} \Q(a)M(\lfloor b)
+\sum_{{{w=ab}\atop{a,b\ne\emptyset}}} \Q(a)\M(b)
-\sum_{{{w=ab}\atop{a,b\ne\emptyset}}} \Q(a\rfloor)\M(b).}\eqno(A.16)$$
Finally, the first, second and fifth terms combine, and the third and fourth combine, 
to form
$$push({\cal A})(v_1,\ldots,v_r)=\sum_{{{w=abc}\atop{ac,b\ne\emptyset}}} \Q(ac)\M(b\rfloor)
-\sum_{{{w=abc}\atop{a,b\ne \emptyset}}} \Q(ac)\M(\lfloor b)
-\sum_{{{w=ab}\atop{a,b\ne\emptyset}}} \Q(a\rfloor)\M(b).\eqno(A.17)$$
Comparing this with (A.3), we see that $push({\cal A})(v_1,\ldots,v_r)$ is 
equal to $-\P(v_1,\ldots,v_r)$. Thus since we rephrased (A.10) as
$$\P\Bigl(sh\bigl((v_1),(v_2,\ldots,v_r)\bigr)\Bigr)={\cal A}(v_1,\ldots,v_r)
+{\cal C}(v_1,\ldots,v_r)$$
and (A.17) shows that
$${\cal A}(v_1,\ldots,v_r)=-push^{-1}\P(v_1,\ldots,v_r),$$
comparison with (A.2) shows that  $C_\P$ is given by ${\cal C}$, i.e.~we have
$$\eqalign{C_\P&(v_1,\ldots,v_r)
=-\sum_{{{w'=bc}\atop{b\ne\emptyset}}}C_\Q(v_1c)\M(b)
+\sum_{{{w'=ab}\atop{b\ne\emptyset}}} C_\Q(v_1a\rfloor)\M(b)
+\sum_{{{w'=abc}\atop{ac,b\ne \emptyset}}} C_\Q(v_1ac)\M(b-a_l)\cr
&-\sum_{{{w'=abc}\atop{ac,b\ne \emptyset}}} C_\Q(v_1ac)\M(b-c_f)
-\sum_{{{w'=ab}\atop{a,b\ne\emptyset}}} (|b|+1)c_{|b|+1}\Q(a-b_1)
+\sum_{{{w'=bc}\atop{b,c\ne\emptyset}}} (|b|+1)c_{|b|+1}\Q(c).}\eqno(A.18)$$
To compute the desired correction $C_{P'}$, recall that
$$C_\P=push^{-1}\cdot swap(C_{P'}), \ \ i.e.\ \ 
C_{P'}=swap\cdot push(C_{\P}).$$
To compute $C_{P'}$ explicitly, we first write (A.18) with indices, then push
it, and finally swap the result.
$$\eqalign{C_\P&(v_1,\ldots,v_r)
=-\sum_{1<i\le r}C_\Q(v_1,v_{i+1},\ldots,v_r)\M(v_2,\ldots,v_i)\cr
&+\sum_{1\le i<r} C_\Q(v_1-v_{i+1},v_2-v_{i+1},\ldots,v_i-v_{i+1})\M(v_{i+1},\ldots,v_r)\cr
&+\sum_{1<i<j\le r} C_\Q(v_1,\ldots,v_i,v_{j+1},\ldots,v_r)\M(v_{i+1}-v_i,\ldots,v_j-v_i)
+\sum_{1<i<r} C_\Q(v_1,v_{i+1},\ldots,v_r)\M(v_2,\ldots,v_i)\cr
&-\sum_{1<i<j<r} C_\Q(v_1,\ldots,v_i,v_{j+1},\ldots,v_r)\M(v_{i+1}-v_{j+1},
\ldots,v_j-v_{j+1})
-\sum_{1<i<r} C_\Q(v_1,\ldots,v_i)\M(v_{i+1},\ldots,v_r)\cr
&-\sum_{1<i<r} (r-i+1)c_{r-i+1}\Q(v_2-v_{i+1},\ldots,v_i-v_{i+1})
+\sum_{1<i<r} iC_i\Q(v_{i+1},\ldots,v_r).}\eqno(A.19)$$
The first and fourth term combine to form 
$$-\sum_{1<i\le r}C_\Q(v_1,v_{i+1},\ldots,v_r)\M(v_2,\ldots,v_i)
+\sum_{1<i<r} C_\Q(v_1,v_{i+1},\ldots,v_r)\M(v_2,\ldots,v_i)=-C_\Q(v_1)\M(v_2,\ldots,v_r)$$
which then combines with the sixth sum, becoming
$$\eqalign{C_\P(v_1,\ldots,v_r)
&=\sum_{1\le i<r} C_\Q(v_1-v_{i+1},v_2-v_{i+1},\ldots,v_i-v_{i+1})\M(v_{i+1},\ldots,v_r)\cr
&+\sum_{1<i<j\le r} C_\Q(v_1,\ldots,v_i,v_{j+1},\ldots,v_r)\M(v_{i+1}-v_i,\ldots,v_j-v_i)\cr
&-\sum_{1<i<j<r} C_\Q(v_1,\ldots,v_i,v_{j+1},\ldots,v_r)\M(v_{i+1}-v_{j+1},\ldots,v_j-v_{j+1})\cr
&-\sum_{1\le i<r} C_\Q(v_1,\ldots,v_i)\M(v_{i+1},\ldots,v_r)\cr 
&-\sum_{1<i<r} (r-i+1)c_{r-i+1}\Q(v_2-v_{i+1},\ldots,v_i-v_{i+1})
+\sum_{1<i<r} iC_i\Q(v_{i+1},\ldots,v_r).}\eqno(A.20)$$
So we have
$$\eqalign{push&(C_\P)(v_1,\ldots,v_r)
=\sum_{1\le i<r} C_\Q(-v_i,v_1-v_i,\ldots,v_{i-1}-v_i)\M(v_i-v_r,\ldots,v_{r-1}-v_r)\cr
&+\sum_{1<i<j\le r} C_\Q(-v_r,v_1-v_r,\ldots,v_{i-1}-v_r,v_j-v_r,\ldots,v_{r-1}-v_r)
\M(v_i-v_{i-1},\ldots,v_{j-1}-v_{i-1})\cr
&-\sum_{1<i<j<r} C_\Q(-v_r,v_1-v_r,\ldots,v_{i-1}-v_r,v_j-v_r,\ldots,v_{r-1}-v_r)
\M(v_i-v_j,\ldots,v_{j-1}-v_j)\cr
&-\sum_{1\le i<r} C_\Q(-v_r,v_1-v_r,\ldots,v_{i-1}-v_r)\M(v_i-v_r,\ldots,v_{r-1}-v_r)\cr 
&-\sum_{1<i<r} (r-i+1)c_{r-i+1}\Q(v_1-v_i,\ldots,v_{i-1}-v_i)
+\sum_{1<i<r} iC_i\Q(v_i-v_r,\ldots,v_{r-1}-v_r).}\eqno(A.21)$$
Rewriting this in terms of $push(C_\Q)$ gives
$$\eqalign{push&(C_\P)(v_1,\ldots,v_r)
=\sum_{1\le i<r} pushC_\Q(v_1,\ldots,v_{i-1},v_i)\M(v_i-v_r,\ldots,v_{r-1}-v_r)\cr
&+\sum_{1<i<j\le r} pushC_\Q(v_1,\ldots,v_{i-1},v_j,\ldots,v_{r-1},v_r)
\M(v_i-v_{i-1},\ldots,v_{j-1}-v_{i-1})\cr
&-\sum_{1<i<j<r} pushC_\Q(v_1,\ldots,v_{i-1},v_j,\ldots,v_{r-1},v_r)
\M(v_i-v_j,\ldots,v_{j-1}-v_j)\cr
&-\sum_{1\le i<r} pushC_\Q(v_1,\ldots,v_{i-1},v_r)\M(v_i-v_r,\ldots,v_{r-1}-v_r)\cr 
&-\sum_{1<i<r} (r-i+1)c_{r-i+1}\Q(v_1-v_i,\ldots,v_{i-1}-v_i)
+\sum_{1<i<r} iC_i\Q(v_i-v_r,\ldots,v_{r-1}-v_r).}\eqno(A.22)$$
Now we compute $C_{P'}=swap\cdot push(C_\P)$. We use the notation $u_i+\cdots+u_j=u_{i,j}$ for $i\le j$.
$$\eqalign{C_{P'}&(u_1,\ldots,u_r)=swap\cdot push(C_\P)(u_1,\ldots,u_r)\cr
&=\sum_{1\le i<r} pushC_\Q(u_{1,r},\ldots,u_{1,r-i+2},u_{1,r-i+1})
\M(u_{2,r-i+1},u_{2,r-i},\ldots,u_2)\cr
&+\sum_{1<i<j\le r} pushC_\Q(u_{1,r},\ldots,u_{1,r-i+2},u_{1,r-j+1},\ldots,u_{1,2},u_1)
\M(-u_{r-i+2},\ldots,-u_{r-j+3,r-i+2})\cr
&-\sum_{1<i<j<r} pushC_\Q(u_{1,r},\ldots,u_{1,r-i+2},u_{1,r-j+1},\ldots,u_{1,2},u_1)
\M(u_{r-j+2,r-i+1},\ldots,u_{r-j+2})\cr
&-\sum_{1\le i<r} pushC_\Q(u_{1,r},\ldots,u_{1,r-i+2},u_1)\M(u_{2,r-i+1},\ldots,u_2)\cr 
&-\sum_{1<i<r} (r-i+1)c_{r-i+1}\Q(u_{r-i+2,r},\ldots,u_{r-i+2})
+\sum_{1<i<r} iC_i\Q(u_{2,r-i+1},\ldots,u_2).}\eqno(A.23)$$
Now we use the fact that
$$swap(C)(u_1,\ldots,u_r)=C(u_{1,r},u_{1,r-1},\ldots,u_{1,2},u_1).$$
$$\eqalign{C_{P'}&(u_1,\ldots,u_r)=swap\cdot push(C_\P)(u_1,\ldots,u_r)\cr
&=\sum_{1\le i<r} swap\cdot pushC_\Q(u_{1,r-i+1},u_{r-i+2},\ldots,u_r)
swap\M(u_2,\ldots,u_{r-i+1})\cr
&+\sum_{1<i<j\le r} swap\cdot pushC_\Q(u_1,\ldots,u_{r-j+1},u_{r-j+2,r-i+2},u_{r-i+3},\ldots,u_r)
swap\M(u_{r-j+3},\ldots,u_{r-i+2})\cr
&-\sum_{1<i<j<r} swap\cdot pushC_\Q(u_1,\ldots,u_{r-j+1},u_{r-j+2,r-i+2},u_{r-i+3},\ldots,u_r)
swap\M(u_{r-j+2},\ldots,u_{r-i+1})\cr
&-\sum_{1\le i<r} swap\cdot pushC_\Q(u_1,u_{2,r-i+2},u_{r-i+3},\ldots,u_r)
swap\M(u_2,\ldots,u_{r-i+1})\cr 
&-\sum_{1<i<r} (r-i+1)c_{r-i+1}swap\Q(u_{r-i+2},\cdots,u_r)
+\sum_{1<i<r} iC_iswap\Q(u_2,\ldots,u_{r-i+1})}\eqno(A.24)$$
where the $\M$-factor in the second term is computed as
$$\M(-u_{r-i+2},\ldots,-u_{r-j+3,r-i+2})
=swap(M)(-u_{r-j+3,r-i+2},u_{r-j+3},\ldots,u_{r-i+1}),$$
which is equal to $swap(M)(u_{r-j+3},\ldots,u_{r-i+2}$ thanks to the push-invariance of $swap(M)$.
Now, recall from (1.2.4) that
$$swap\bigl(\F(Q))(v_1,\ldots,v_r)=swap(Q)(v_1,\ldots,v_r)+swap(Q)\Bigl(sh\bigl((v_1),
(v_2,\ldots,v_r)\bigr)\Bigr),$$
so since $\F(Q)=C_Q$ and
$$C_{\Q}(v_1,\ldots,v_r):=\Q\Bigl(sh\bigl((v_1),(v_2,\ldots,v_r)\bigr)\Bigr)+push^{-1}\Q(v_1,\ldots,v_r),$$
we have $C_Q=swap\cdot push(C_\Q)$.  Replacing this in (A.24) as well as $Q=swap\Q$, $M=swap\M$, we obtain
$$\eqalign{C_{P'}(u_1,\ldots,u_r)
&=\sum_{1\le i<r} C_Q(u_{1,r-i+1},u_{r-i+2},\ldots,u_r) M(u_2,\ldots,u_{r-i+1})\cr
&+\sum_{1<i<j\le r} C_Q(u_1,\ldots,u_{r-j+1},u_{r-j+2,r-i+2},u_{r-i+3},\ldots,u_r)
M(u_{r-j+3},\ldots,u_{r-i+2})\cr
&-\sum_{1<i<j<r} C_Q(u_1,\ldots,u_{r-j+1},u_{r-j+2,r-i+2},u_{r-i+3},\ldots,u_r)
M(u_{r-j+2},\ldots,u_{r-i+1})\cr
&-\sum_{1\le i<r} C_Q(u_1,u_{2,r-i+2},u_{r-i+3},\ldots,u_r) M(u_2,\ldots,u_{r-i+1})\cr 
&-\sum_{i=2}^{r-1} (r-i+1)c_{r-i+1}Q(u_{r-i+2},\cdots,u_r)
+\sum_{i=2}^{r-1} iC_iQ(u_2,\ldots,u_{r-i+1}).}\eqno(A.25)$$
The first sum corresponds to the term $C_Q(a\rceil c)M(b)$ for 
$$abc=(u_1)(u_2,\ldots,u_{r-i+1})(u_{r-i+2},\ldots,u_r),$$
where the condition $i<r$ means that $b\ne\emptyset$,
and the second sum corresponds to the same term $C_Q(a\rceil c)M(b)$ for
$$abc=(u_1,\ldots,u_{r-j+2})(u_{r-j+3},\ldots,u_{r-i+2})(u_{r-i+3},\ldots,u_r)$$
with $j<r$, i.e. $|a|\ge 2$ and $b\ne\emptyset$; thus these two terms combine as the sum over
$w=abc$ with $a,b\ne\emptyset$ of the terms $C_Q(a\rceil c)M(b)$.  Similarly, the fourth
sum corresponds to the term $C_Q(a\lceil c)M(b)$ for
$$abc=(u_1)(u_2,\ldots,u_{r-i+1})(u_{r-i+2},\ldots,u_r)$$
so $|a|=1$ and again the condition $i<r$ means that $b\ne\emptyset$, while the third sum
corresponds to the sum of the terms $C_Q(a\lceil c)M(b)$ over the decompositions
$$abc=(u_1,\ldots,u_{r-j+1})(u_{r-j+2},\ldots,u_{r-i+1})(u_{r-i+2},\ldots,u_r),$$
in which the condition $j<r$ implies that $|a|\ge 2$ and the condition $i<j$ means that $b\ne\emptyset$,
so again these two terms combine to form the sum of $C_Q(a\lceil c)M(b)$ over the decompositions
$w=abc$ with $a,b\ne\emptyset$.  In terms of indices, this means that using the standard form
$$abc=(u_1,\ldots,u_i)(u_{i+1},\ldots,u_j)(u_{j+1},\ldots,u_r),$$
the terms of (A.25) can be reindexed as
$$\eqalign{C_{P'}(u_1,\ldots,u_r)
&=\sum_{1<i<j\le r} C_Q(u_1,\ldots,u_{i-1},u_i+\cdots+u_j,u_{j+1},\ldots,u_r)M(u_{i+1},\ldots,u_j)\cr
&-\sum_{1<i<j\le r} C_Q(u_1,\ldots,u_i,u_{i+1}+\cdots+u_{j+1},u_{j+2},\ldots,u_r)M(u_{i+1},\ldots,u_j)\cr
&-\sum_{i=2}^{r-1} iC_iQ(u_{i+1},\cdots,u_r)
+\sum_{i=2}^{r-1} iC_iQ(u_2,\ldots,u_{r-i+1}).}\eqno(A.26)$$
This concludes the proof of Theorem 3.2.
\vskip 1.5cm
\noindent {\bf Appendix B: Proof of Lemma 3.3.}
\vskip 1cm
Since $-a+{{1}\over{2}}ad(b)(a)$ is the start of
the power series expansion of $t_{01}$ in $a,b$,
$\hat{t}_{01}$ is simply the depth $r\ge 2$ part of $t_{01}$,
whose associated mould is extremely simple;
indeed, for each even depth $r\ge 2$, we have 
$$\hat{T}_{01}(u_1,\ldots,u_r)={{B_r}\over{r!}}T_r(u_1,\ldots,u_r),$$
where $T_r$ is the mould concentrated in depth $r$ defined by
$$T_r(u_1,\ldots,u_r):=
\sum_{i=1}^r (-1)^{i-1} \Bigl({{r-1}\atop{i-1}}\Bigr)u_i.\eqno(B.1)$$
It suffices to check that each $T_r$ satisfies the strict Fay 
	relations for even $r$, i.e.~that for all even $r\ge 2$,
	$$T'_r(u_1,\ldots,u_r)+\sum_{j=1}^{r-1} T'_r(u_2,\ldots,u_j,-\overline{u}_j,\overline{u}_{j+1},u_{j+2},\ldots,u_r)
	+T'_r(u_2,\ldots,u_r,-u_{1,r}) =0\eqno(B.2)$$
	which we write out explicitly as
	$${{\sum_{i=1}^r (-1)^{i-1} \Bigl({{r-1}\atop{i-1}}\Bigr)u_i}\over
		{u_1\cdots u_r}}=
{{\sum_{i=1}^{r-1} (-1)^{i-1}\Bigl({{r-1}\atop{i-1}}\Bigr)u_{i+1}-(-1)^{r-1}u_{1,r}}\over{u_2\cdots u_ru_{1,r}}}$$
$$+\sum_{j=1}^{r-1}
{{\sum_{i=1}^{j-1}(-1)^{i-1}\Bigl({{r-1}\atop{i-1}}\Bigr)u_{i+1}
	-(-1)^{j-1}\Bigl({{r-1}\atop{j-1}}\Bigr)u_{1,j}
	+(-1)^j\Bigl({{r-1}\atop{j}}\Bigr)u_{1,j+1}
	+\sum_{i=j+2}^r(-1)^{i-1}\Bigl({{r-1}\atop{i-1}}\Bigr)u_i}
\over{u_2\cdots u_ju_{1,j}u_{1,j+1}u_{j+2}\cdots u_r}}.\eqno(B.3)$$
Let us put this over the common denominator
$$D=u_1\cdots u_ru_{1,2}\cdots u_{1,r};$$
the equality we need to prove is then given by
$$u_{1,2}\cdots u_{1,r}\sum_{i=1}^r (-1)^{i-1}\Bigl({{r-1}\atop{i-1}}\Bigr)u_i$$
$$=\sum_{j=1}^{r-1} D_j
\biggl(\sum_{i=1}^{j-1}(-1)^{i-1}\Bigl({{r-1}\atop{i-1}}\Bigr)u_{i+1}
		-(-1)^{j-1}\Bigl({{r-1}\atop{j-1}}\Bigr)u_{1,j}
		+(-1)^j\Bigl({{r-1}\atop{j}}\Bigr)u_{1,j+1}
		+\sum_{i=j+2}^r(-1)^{i-1}\Bigl({{r-1}\atop{i-1}}\Bigr)u_i\biggr).\eqno(B.4)
$$
where 
$$D_j=\cases{u_1u_{j+1}u_{1,2}\cdots u_{1,j-1}u_{1,j+2}\cdots u_{1,r}&$1\le j\le r-1$\cr
	u_1u_{1,2}\cdots u_{1,r-1}&$j=r$.}$$
	To prove (B.4), we first show that the degree $r$ polynomial on the
	right-hand side is divisible by $u_{1,k}$ for $k=2,\ldots,r$.  The
	factor $u_{1,k}$ appears in $D_j$ for all $j\ne k-1,k$, so if we set
	$u_{1,k}=0$, i.e. $u_k=-u_{k+1}-\cdots-u_r$ in the right-hand side of (B.4),
	only the two terms 
	$$D_{k-1}\biggl(\sum_{i=1}^{k-2}(-1)^{i-1}\Bigl({{r-1}\atop{i-1}}\Bigr)u_{i+1}
			-(-1)^{k-2}\Bigl({{r-1}\atop{k-2}}\Bigr)u_{1,k-1}
			+(-1)^{k-1}\Bigl({{r-1}\atop{k-1}}\Bigr)u_{1,k}
			+\sum_{i=k+1}^r (-1)^{i-1}\Bigl({{r-1}\atop{i-1}}\Bigr)u_i\biggr)$$
	$$+D_k\biggl(\sum_{i=1}^{k-1}(-1)^{i-1}\Bigl({{r-1}\atop{i-1}}\Bigr)u_{i+1}
			-(-1)^{k-1}\Bigl({{r-1}\atop{k-1}}\Bigr)u_{1,k}
			+(-1)^k\Bigl({{r-1}\atop{k}}\Bigr)u_{1,k+1}
			+\sum_{i=k+2}^r (-1)^{i-1}\Bigl({{r-1}\atop{i-1}}\Bigr)u_i\biggr).\eqno(B.5)$$
	We write
	$$D_{k-1}=Du_ku_{1,k+1}\ \ \ {\rm and}\ \ \ 
	D_k=Du_{k+1}u_{1,k-1}$$
	with 
	$$D=u_1u_{1,2}\cdots u_{1,k-2}u_{1,k+2}\cdots u_{1,r}.$$
	Observe that if $u_{1,k}=0$ then $u_{1,k-1}=-u_k$ and $u_{1,k+1}=u_{k+1}$.
	We use this to rewrite (B.5) as
	$$Du_ku_{k+1}\biggl(\sum_{i=1}^{k-2}(-1)^{i-1}\Bigl({{r-1}\atop{i-1}}\Bigr)u_{i+1}
			+(-1)^{k-2}\Bigl({{r-1}\atop{k-2}}\Bigr)u_k
			+\sum_{i=k+1}^r (-1)^{i-1}\Bigl({{r-1}\atop{i-1}}\Bigr)u_i\biggr)$$
	$$-Du_ku_{k+1}\biggl(\sum_{i=1}^{k-1}(-1)^{i-1}\Bigl({{r-1}\atop{i-1}}\Bigr)u_{i+1}
			+(-1)^k\Bigl({{r-1}\atop{k}}\Bigr)u_{k+1}
			+\sum_{i=k+2}^r (-1)^{i-1}\Bigl({{r-1}\atop{i-1}}\Bigr)u_i\biggr),$$
	or better
	$$Du_ku_{k+1}\biggl(\sum_{i=1}^{k-2}(-1)^{i-1}\Bigl({{r-1}\atop{i-1}}\Bigr)u_{i+1}
			+(-1)^{k-2}\Bigl({{r-1}\atop{k-2}}\Bigr)u_k
			+(-1)^k\Bigl({{r-1}\atop{k}}\Bigr)u_{k+1}
			+\sum_{i=k+2}^r (-1)^{i-1}\Bigl({{r-1}\atop{i-1}}\Bigr)u_i\biggr)$$
	$$-Du_ku_{k+1}\biggl(\sum_{i=1}^{k-2}(-1)^{i-1}\Bigl({{r-1}\atop{i-1}}\Bigr)u_{i+1}
			+(-1)^{k-2}\Bigl({{r-1}\atop{k-2}}\Bigr)u_k
			+(-1)^k\Bigl({{r-1}\atop{k}}\Bigr)u_{k+1}
			+\sum_{i=k+2}^r (-1)^{i-1}\Bigl({{r-1}\atop{i-1}}\Bigr)u_i\biggr)$$
	which is identically equal to zero.  This shows that the right-hand side of
	(B.4) is divisible by $u_{1,k}$ for $k=2,\ldots,r$.  To conclude the
	proof of the equality (B.4), we need to show that the right-hand side
	also vanishes when we set
	$$\sum_{i=1}^r (-1)^{i-1}\Bigl({{r-1}\atop{i-1}}\Bigr)u_i=0,$$
	an equality which can be written equivalently as
	$$\sum_{i=j+2}^r(-1)^{i-1}\Bigl({{r-1}\atop{i-1}}\Bigr)u_i=
	\sum_{i=1}^{j+1} (-1)^i\Bigl({{r-1}\atop{i-1}}\Bigr)u_i$$
	for $1\le j\le r$.  Let us make this substitution in the right-hand side of
	(B.4) (and reindex the first sum from $2$ to $j$ instead of $1$ to $j-1$):
		$$=\sum_{j=1}^r D_j   
		\biggl(\sum_{i=2}^{j}(-1)^{i-2}\Bigl({{r-1}\atop{i-2}}\Bigr)u_i
				-(-1)^{j-1}\Bigl({{r-1}\atop{j-1}}\Bigr)u_{1,j}
				+(-1)^j\Bigl({{r-1}\atop{j}}\Bigr)u_{1,j+1}
				+\sum_{i=1}^{j+1} (-1)^i\Bigl({{r-1}\atop{i-1}}\Bigr)u_i\biggr)$$
		$$=\sum_{j=1}^r D_j  \Biggl( 
				(-1)^{j}\Bigl({{r-1}\atop{j-1}}\Bigr)u_{1,j}
				+(-1)^j\Bigl({{r-1}\atop{j}}\Bigr)u_{1,j+1}
				-u_1-(-1)^j\Bigl({{r-1}\atop{j}}\Bigr)u_{j+1}
				+\sum_{i=2}^{j}(-1)^{i}\biggl(\Bigl({{r-1}\atop{i-2}}\Bigr)
					+\Bigl({{r-1}\atop{i-1}}\Bigr)\biggr)u_i\Biggr)$$
		$$=\sum_{j=1}^r D_j  \Biggl(-u_1 
				+(-1)^{j}\Bigl({{r-1}\atop{j-1}}\Bigr)u_{1,j}
				+(-1)^j\Bigl({{r-1}\atop{j}}\Bigr)u_{1,j}
				+(-1)^i\sum_{i=2}^{j}\biggl(\Bigl({{r-1}\atop{i-2}}\Bigr)
					+\Bigl({{r-1}\atop{i-1}}\Bigr)\biggr)u_i\Biggr)$$
		$$=\sum_{j=1}^r D_j  \Biggl(-u_1 
				+(-1)^{j}\biggl(\Bigl({{r-1}\atop{j-1}}\Bigr)
					+\Bigl({{r-1}\atop{j}}\Bigr)\biggr)u_{1,j}
				+(-1)^i\sum_{i=2}^{j}\biggl(\Bigl({{r-1}\atop{i-2}}\Bigr)
					+\Bigl({{r-1}\atop{i-1}}\Bigr)\biggr)u_i\Biggr)$$
		$$=\sum_{j=1}^r D_j  \Biggl(-u_1 
				+(-1)^{j}\biggl({{(r-1)!}\over{(j-1)!(r-j)!}}
					+{{(r-1!)}\over{j!(r-j-1)!}}\biggr)
				u_{1,j}
				+(-1)^i\sum_{i=2}^{j}
				\biggl({{(r-1)!}\over{(i-2)!(r-i+1)!}}
					+\Bigl({{(r-1)!}\over{(i-1)!(r-i)!}}\biggr)
					u_i\Biggr)$$
				$$=\sum_{j=1}^r D_j  \Biggl(-u_1 
					+(-1)^{j}\biggl(j\times {{(r-1)!}\over{j!(r-j)!}}
						+(r-j)\times {{(r-1!)}\over{j!(r-j)!}}\biggr)
					u_{1,j}$$
					$$\qquad\qquad +\sum_{i=2}^{j}(-1)^{i}
					\biggl((i-1)\times {{(r-1)!}\over{(i-1)!(r-i+1)!}}
						+(r-i+1)\times \Bigl({{(r-1)!}\over{(i-1)!(r-i+1)!}}\biggr)
						u_i\Biggr)$$
					$$=\sum_{j=1}^r D_j  \Biggl(-u_1 
						+(-1)^{j}\Bigl({{r}\atop{j}}\Bigr)u_{1,j}
						+\sum_{i=2}^{j}(-1)^{i}
						\Bigl({{r}\atop{i-1}}\Bigr) u_i\Biggr)$$
					$$=\sum_{j=1}^r D_j  \Biggl(-u_1 
						+\sum_{i=1}^j (-1)^{j}\Bigl({{r}\atop{j}}\Bigr)u_i
						+\sum_{i=2}^{j}(-1)^{i}
						\Bigl({{r}\atop{i-1}}\Bigr) u_i\Biggr)$$
					$$=\sum_{j=1}^r D_j  \Biggl(
						\sum_{i=1}^j \biggl((-1)^{j}\Bigl({{r}\atop{j}}\Bigr)
							+(-1)^{i}\Bigl({{r}\atop{i-1}}\Bigr)\biggr) u_i\Biggr)$$
					$$=\sum_{j=1}^r 
					u_1u_{j+1}u_{1,2}\cdots u_{1,j-1}u_{1,j+2}\cdots u_{1,r}\Biggl(
						\sum_{i=1}^j \biggl((-1)^{j}\Bigl({{r}\atop{j}}\Bigr)
							+(-1)^{i}\Bigl({{r}\atop{i-1}}\Bigr)\biggr) u_i\Biggr).\eqno(B.6)$$
	Let us call this polynomial $P$.
	\vskip .2cm
	\noindent {\bf Claim.} $u_{1,k}$ is a factor of $P$ for $k=1,\ldots,r$.

	For $k=1$ this is obvious since $u_1$ appears as a factor in every term.
	Fix a $k$ with $1<k<r$; then there are only two terms in $P$ where $u_{1,k}$ 
	does not appear, namely the term for $j=k-1$ and the term for $j=k$:
$$u_1u_{k}u_{1,2}\cdots u_{1,k-2}u_{1,k+1}\cdots u_{1,r}\Biggl(
		\sum_{i=1}^{k-1} \biggl((-1)^{k-1}\Bigl({{r}\atop{k-1}}\Bigr)
			+(-1)^{i}\Bigl({{r}\atop{i-1}}\Bigr)\biggr) u_i\Biggr)$$
$$+u_1u_{k+1}u_{1,2}\cdots u_{1,k-1}u_{1,k+2}\cdots u_{1,r}\Biggl(
		\sum_{i=1}^k \biggl((-1)^{k}\Bigl({{r}\atop{k}}\Bigr)
			+(-1)^{i}\Bigl({{r}\atop{i-1}}\Bigr)\biggr) u_i\Biggr).$$
	Add and subtract a term between these two lines as follows:
$$u_1u_{k}u_{1,2}\cdots u_{1,k-2}u_{1,k+1}\cdots u_{1,r}\Biggl(
		\sum_{i=1}^{k-1} \biggl((-1)^{k-1}\Bigl({{r}\atop{k-1}}\Bigr)
			+(-1)^{i}\Bigl({{r}\atop{i-1}}\Bigr)\biggr) u_i\Biggr)$$
$$+u_1u_{1,k-1}u_{1,2}\cdots u_{1,k-2}u_{1,k+1}\cdots u_{1,r}\Biggl(
		\sum_{i=1}^{k-1} \biggl((-1)^{k-1}\Bigl({{r}\atop{k-1}}\Bigr)
			+(-1)^{i}\Bigl({{r}\atop{i-1}}\Bigr)\biggr) u_i\Biggr)$$
$$-u_1u_{1,k-1}u_{1,2}\cdots u_{1,k-2}u_{1,k+1}\cdots u_{1,r}\Biggl(
		\sum_{i=1}^{k-1} \biggl((-1)^{k-1}\Bigl({{r}\atop{k-1}}\Bigr)
			+(-1)^{i}\Bigl({{r}\atop{i-1}}\Bigr)\biggr) u_i\Biggr)$$
$$+u_1u_{k+1}u_{1,2}\cdots u_{1,k-1}u_{1,k+2}\cdots u_{1,r}\Biggl(
		\sum_{i=1}^k \biggl((-1)^{k}\Bigl({{r}\atop{k}}\Bigr)
			+(-1)^{i}\Bigl({{r}\atop{i-1}}\Bigr)\biggr) u_i\Biggr).$$
$$=u_1u_{1,k}u_{1,2}\cdots u_{1,k-2}u_{1,k+1}\cdots u_{1,r}\Biggl(
		\sum_{i=1}^{k-1} \biggl((-1)^{k-1}\Bigl({{r}\atop{k-1}}\Bigr)
			+(-1)^{i}\Bigl({{r}\atop{i-1}}\Bigr)\biggr) u_i\Biggr)$$
	$$-u_1u_{1,2}\cdots u_{1,k-2}u_{1,k-1}(u_{1,k}+u_{k+1})\cdots u_{1,r}\Biggl(
			\sum_{i=1}^{k-1} \biggl((-1)^{k-1}\Bigl({{r}\atop{k-1}}\Bigr)
				+(-1)^{i}\Bigl({{r}\atop{i-1}}\Bigr)\biggr) u_i\Biggr)$$
$$+u_1u_{k+1}u_{1,2}\cdots u_{1,k-1}u_{1,k+2}\cdots u_{1,r}\Biggl(
		\sum_{i=1}^k \biggl((-1)^{k}\Bigl({{r}\atop{k}}\Bigr)
			+(-1)^{i}\Bigl({{r}\atop{i-1}}\Bigr)\biggr) u_i\Biggr).$$
$$=u_1u_{1,k}u_{1,2}\cdots u_{1,k-2}u_{1,k+1}\cdots u_{1,r}\Biggl(
		\sum_{i=1}^{k-1} \biggl((-1)^{k-1}\Bigl({{r}\atop{k-1}}\Bigr)
			+(-1)^{i}\Bigl({{r}\atop{i-1}}\Bigr)\biggr) u_i\Biggr)$$
$$-u_1u_{1,2}\cdots u_{1,k-2}u_{1,k-1}u_{1,k}\cdots u_{1,r}\Biggl(
		\sum_{i=1}^{k-1} \biggl((-1)^{k-1}\Bigl({{r}\atop{k-1}}\Bigr)
			+(-1)^{i}\Bigl({{r}\atop{i-1}}\Bigr)\biggr) u_i\Biggr)$$
$$-u_1u_{k+1}u_{1,2}\cdots u_{1,k-2}u_{1,k-1}\cdots u_{1,r}\Biggl(
		\sum_{i=1}^{k-1} \biggl((-1)^{k-1}\Bigl({{r}\atop{k-1}}\Bigr)
			+(-1)^{i}\Bigl({{r}\atop{i-1}}\Bigr)\biggr) u_i\Biggr)$$
$$+u_1u_{k+1}u_{1,2}\cdots u_{1,k-1}u_{1,k+2}\cdots u_{1,r}\Biggl(
		\sum_{i=1}^k \biggl((-1)^{k}\Bigl({{r}\atop{k}}\Bigr)
			+(-1)^{i}\Bigl({{r}\atop{i-1}}\Bigr)\biggr) u_i\Biggr).$$
$$=u_1u_{1,k}u_{1,2}\cdots u_{1,k-2}u_{1,k+1}\cdots u_{1,r}\Biggl(
		\sum_{i=1}^{k-1} \biggl((-1)^{k-1}\Bigl({{r}\atop{k-1}}\Bigr)
			+(-1)^{i}\Bigl({{r}\atop{i-1}}\Bigr)\biggr) u_i\Biggr)$$
$$-u_1u_{1,2}\cdots u_{1,k-2}u_{1,k-1}u_{1,k}\cdots u_{1,r}\Biggl(
		\sum_{i=1}^{k-1} \biggl((-1)^{k-1}\Bigl({{r}\atop{k-1}}\Bigr)
			+(-1)^{i}\Bigl({{r}\atop{i-1}}\Bigr)\biggr) u_i\Biggr)$$
$$+u_1u_{k+1}u_{1,2}\cdots u_{1,k-1}u_{1,k+2}\cdots u_{1,r}\Biggl(
		\sum_{i=1}^k \biggl((-1)^{k}\Bigl({{r}\atop{k}}\Bigr)
			+(-1)^{i}\Bigl({{r}\atop{i-1}}\Bigr)\biggr) u_i-
		\sum_{i=1}^{k-1} \biggl((-1)^{k-1}\Bigl({{r}\atop{k-1}}\Bigr)
			+(-1)^{i}\Bigl({{r}\atop{i-1}}\Bigr)\biggr) u_i\Biggr).$$
	The factor $u_{1,k}$ appears in the first two lines of this, but in the
	third line, we will show that it appears in the large bracketed factor.
	For this, we compute the coefficient of $u_i$ in the large bracketed
	factor.  For $i>k$, $u_i$ does not occur.  The coefficient of $u_k$ is
$$\eqalign{(-1)^k\Bigl({{r}\atop{k}}\Bigr)+(-1)^k\Bigl({{r}\atop{k-1}}\Bigr)
	&=(-1)^k\biggl({{r!}\over{k!(r-k)!}}+{{r!}\over{(k-1)!(r-k+1)!}}\biggr)\cr
		&=(-1)^k\biggl((r-k+1)\times {{r!}\over{k!(r-k+1)!}}+k\times {{r!}\over{k!(r-k+1)!}}\biggr)\cr
		&=(-1)^k\biggl({{(r+1)!}\over{k!(r-k+1)!}}\biggr)\cr
		&=(-1)^k\Bigl({{r+1}\atop{k}}\Bigr).}$$
	For $1\le i<k$, the coefficient is
$$(-1)^k\Bigl({{r}\atop{k}}\Bigr)+(-1)^i\Bigl({{r}\atop{i-1}}\Bigr)
	-(-1)^{k-1}\Bigl({{r}\atop{k-1}}\Bigr)-(-1)^i\Bigl({{r}\atop{i-1}}\Bigr)$$
$$(-1)^k\Bigl({{r+1}\atop{k}}\Bigr)+(-1)^i\Bigl({{r}\atop{i-1}}\Bigr)
	-(-1)^i\Bigl({{r}\atop{i-1}}\Bigr)$$
	$$(-1)^k\Bigl({{r+1}\atop{k}}\Bigr).$$
	Thus the large bracketed factor in the last line is nothing other than
	$$(-1)^k\Bigl({{r+1}\atop{k}}\Bigr)u_{1,k}$$
	and thus every term contains the factor $u_{1,k}$; thus $u_{1,k}$ is a
	factor of $P$, proving the claim.

	Since $P$ is of degree $r$ and $u_{1,k}$ is a factor of $P$ for $1\le k\le r$,
	we have $P=cu_{1,1}u_{1,2}\cdots u_{1,r}$.  Since $u_{1,j}=u_1+\cdots+u_j$, 
	when we expand $P$ in the variables $u_1,\ldots,u_r$, the monomial $u_1^r$ 
	appears with coefficient $c$.  In order to compute the monomial $cu_1^r$, 
	we consider the expression of $P$ given in (B.6) and set $u_2=\cdots=u_r=0$.
	Because $u_{j+1}$ appears in every term of $P$ except for $j=r$, 
	using the fact that $u_{1,j}$ becomes equal to $u_1$ when $u_2=\cdots=u_r=0$,
	we only need to consider the term $j=r$ and $i=1$: we find
	$$P(u_1,0,\ldots,0)=u_1^{r-1}\times \Bigl((-1)^r+(-1)\Bigr)u_1.$$
	Given that $r$ is even, this is equal to 0 and thus $P=0$, which concludes
	the unfortunately complicated proof of Lemma 3.3.\hfill{$\diamondsuit$}
	\vskip 1cm
	\noindent {\bf References}
	\vskip .4cm
	\noindent{[AT]} A.~Alekseev, C.~Torossian, The Kashiwara-Vergne conjecture
	and Drinfeld's associators, {\it Annals of Math.} {\bf 175:2} (2012), 415-463.
	\vskip .2cm
	\noindent{[BMS]} J.~Broedel, N.~Matthes, O.~Schlotterer, Relations between elliptic multiple zeta values and a special derivation algebra, {\it J. Phys. A} {\bf 49} (2016) 15, 155203.
	\vskip .2cm
\noindent{[Ec]} J.~\'Ecalle (with computational assistancy from S.~Carr), The flexion structure and dimorphy: flexion units, singulators, generators, and the enumeration of multizeta irreducibles, in {\it Asymptotics in dynamics, geometry and PDEs}; Pisa, 2011.
	\vskip .2cm
	\noindent{[En1]} B.~Enriquez, Elliptic associators, {\it Selecta Math.} (N.S.) {\bf 20:2} (2014), 491-584.
	\vskip .2cm
	\noindent{[En2]} B.~Enriquez, Analogues elliptiques des nombres multiz\^etas, {\it Bull.
		Soc. Math. France} {\bf 144:3} (2016), 395-427.
		\vskip .2cm
		\noindent{[HM} R.~Hain and M.~Matsumoto, Universal mixed elliptic motives,
{\it J. Inst. Math. Jussieu} {\bf 19(3)} (2020), 663-766.
\vskip .2cm
\noindent{[LMS]} P.~Lochak, N.~Matthes, L.~Schneps, Elliptic multizetas and the elliptic double shuffle relations, {\it IMRN} {\bf 2021(1)}, 695-763.
\vskip .2cm
\noindent{[R]} G.~Racinet, S\'erie g\'en\'eratrices non commutatives de polyz\^etas
et associateurs de Drinfel'd, thesis, 2000.
\vskip .2cm
\noindent{[RS]} E.~Raphael, L.~Schneps, On the elliptic Kashiwara-Vergne Lie algebra, ArXiv:1809.09340, to appear in {\it J. Lie Algebra}.
\vskip .2cm
\noindent{[S1]} L.~Schneps, Elliptic double shuffle, Grothendieck-Teichm\"uller and mould theory, {\it Annales Math. Qu\'ebec} {\bf 44(2)} (2020), 261-289.
\vskip .2cm
\noindent{[S2]} L.~Schneps, ARI, GARI, Zig and Zag, arXiv:1507.01535, 2015.
\vskip .2cm
\noindent{[SS]} A.~Salerno, L.~Schneps, Mould theorey and the double
shuffle Lie algebra structure, in {\it Periods in Quantum Field
Theory and Arithmetic}, Proceedings in Mathematics and Statistics,
Springer 2019.

\bye
\vfill\eject
END OF ARTICLE HERE
\vfill\eject
\noindent We have 
$$P=arat(M)\cdot Q=-arit(M)\cdot Q+lu(M,Q),$$
so we compute
$$\eqalign{swap(P)&=swap\bigl(arat(M)\cdot Q\bigr)\cr
	&=-swap\bigl(arit(M)\cdot Q\bigr)+swap\bigl(mu(M,Q)\bigr)-swap\bigl(mu(Q,M)\bigr)\cr
		&=-axit\bigl(swapM,-push(swapM)\bigr)\cdot swap(Q)-mu\bigl(swap(Q),swapM\bigr)\cr
		&\qquad+swap\bigl(mu(Q,M)\bigr)+swap\bigl(mu(M,Q)\bigr)-swap\bigl(mu(Q,M)\bigr)\cr
		&=-axit\bigl(swapM,-push(swapM)\bigr)\cdot swap(Q)-mu\bigl(swap(Q),swapM\bigr) +swap\bigl(mu(M,Q)\bigr)\cr
		&=-arit\bigl(swapM\bigr)\cdot swap(Q)-mu\bigl(swap(Q),swapM\bigr)
		+swap\bigl(mu(M,Q)\bigr),}\eqno(A.3)$$
		where the last equality follows since $M$ is bialternal and thus 
		$push$-invariant, and $axit(P,-P)=arit(P)$.
		In order to compute the left-hand side of (A.1), we will use (A.3) and
	separately compute the three terms
$$\cases{\Bigl(arit\bigl(swapM\bigr)\cdot swap(Q)\Bigr)
	\bigl(sh\bigl((v_1),(v_2,\ldots,v_r)\bigr)\bigr)\cr
		\Bigl(mu\bigl(swap(Q),swapM\bigr)\Bigr)
		\bigl(sh\bigl((v_1),(v_2,\ldots,v_r)\bigr)\bigr)\cr
		\Bigl(swap\bigl(mu(M,Q)\bigr)\Bigr)
		\bigl(sh\bigl((v_1),(v_2,\ldots,v_r)\bigr)\bigr).}\eqno(A.3)$$

		INTERLUDE

$$\sum_{w=abc,a\ne 0} \Q(ac)\M(\lfloor b)-\sum_{w=abc,c\ne 0} \Q(ac)M(b\rfloor)
	-\sum_{w=ab} \Q(a)\M(b)+\sum_{w=ab} \Q(a\rfloor)\M(b)$$
	The two middle terms join:
$$\sum_{w=abc,a\ne 0} \Q(ac)\M(\lfloor b)-\sum_{w=abc} \Q(ac)M(b\rfloor)
	+\sum_{w=ab} \Q(a\rfloor)\M(b)$$
	Simplify by adding and subtracting the first term part with $a=0$:
$$\sum_{w=abc} \Q(ac)\M(\lfloor b)-\sum_{w=abc} \Q(ac)M(b\rfloor)
	+\sum_{w=ab} \Q(a\rfloor)\M(b)-\sum_{w=abc,a=0}\Q(ac)M(\lfloor b)$$
$$=\sum_{w=abc} \Q(ac)\M(\lfloor b)-\sum_{w=abc} \Q(ac)M(b\rfloor)
	+\sum_{w=ab} \Q(a\rfloor)\M(b)-\sum_{w=bc}\Q(c)\M(b)$$

	END INTERLUDE

	\noindent {\bf Left-hand side: first term of (A.3).} We compute
	$$\Bigl(arit\bigl(swapM\bigr)\cdot swap(Q)\Bigr)\bigl(sh((v_1),(v_2,\ldots,v_r))\bigr)$$
	by adapting a method due to A.~Salerno (cf. [SS, appendix]), as follows.
	The definition of the $arit$-operator on moulds in $\overline{ARI}$ is given
	by the following formula.
	$$\Bigl(arit\bigl(swapM\bigr)\cdot swap(Q)\Bigr)(w)=\sum_{w=abc,c\ne 0}
	Q(ac)M(b\rfloor) - \sum_{w=abc, a\ne 0} Q(ac)M(\lfloor b).$$
	We thus have
	$$\Bigl(arit(swapM)\cdot swap(Q)\Bigr)
\bigl(sh((v_1),(v_2,\ldots,v_r))\bigr) 
	\qquad\qquad\qquad\qquad\qquad \qquad\qquad\qquad\qquad\qquad\qquad\qquad\qquad  $$

	a        b        c
	v1v2v3   v4       0
	v1v2     v3v4     0
	v1v2     v3       v4 
	v1       v2v3v4   0
	v1       v2v3     v4
	v1       v2       v3v4
	0       v1v2v3    v4
	0       v1v2      v3v4
	0       v1        v2v3v4
	v2v1v3   v4       0
	v2v1     v3v4     0
	v2v1     v3       v4 
	v2       v1v3v4   0
	v2       v1v3     v4
	v2       v1       v3v4
	0       v2v1v3    v4
	0       v2v1      v3v4
	0       v2        v1v3v4
	v2v3v1   v4       0
	v2v3     v1v4     0
	v2v3     v1       v4 
	v2       v3v1v4   0
	v2       v3v1     v4
	v2       v3       v1v4
	0       v2v3v1    v4
	0       v2v3      v1v4
	0       v2        v3v1v4
	v2v3v4   v1       0
	v2v3     v4v1     0
	v2v3     v4       v1 
	v2       v3v4v1   0
	v2       v3v4     v1
	v2       v3       v4v1
	0       v2v3v4    v1
	0       v2v3      v4v1
	0       v2        v3v4v1

	(a)(b)(c)           sh(v1,a)(b)(c)
	abc=(v2v3)(v4)(0)   (v1v2v3)(v4)(0)   (v2v1v3)(v4)(0)   (v2v3v1)(v4)(0)  
	abc=(v2)(v3v4)(0)   (v1v2)(v3v4)(0)   (v2v1)(v3v4)(0)
	abc=(v2)(v3)(v4)    (v1v2)(v3)(v4)    (v2v1)(v3)(v4)
	abc=(0)(v2v3v4)(0)  (v1)(v2v3v4)(0)
	abc=(0)(v2v3)(v4)   (v1)(v2v3)(v4)
abc=(0)(v2)(v3v4)   (v1)(v2)(v3v4)

	(a)sh(v1,b)(c)
	abc=(v2v3v4)(0)(0)  (v2v3v4)(v1)(0)
	abc=(v2v3)(v4)(0)   (v2v3)(v1v4)(0)   (v2v3)(v4v1)(0)
	abc=(v2v3)(0)(v4)   (v2v3)(v1)(v4)
	abc=(v2)(v3v4)(0)   (v2)(v1v3v4)(0)  (v2)(v3v1v4)(0)  (v2)(v3v4v1)(0)
	abc=(v2)(v3)(v4)    (v2)(v1v3)(v4)   (v2)(v3v1)(v4)
	abc=(v2)(0)(v3v4)   (v2)(v1)(v3v4)
	abc=(0)(v2v3)(v4)   (0)(v1v2v3)(v4)  (0)(v2v1v3)(v4)   (0)(v2v3v1)(v4)
	abc=(0)(v2)(v3v4)   (0)(v1v2)(v3v4)  (0)(v2v1)(v3v4)
abc=(0)(0)(v2v3v4)  (0)(v1)(v2v3v4)

	(a)(b)sh(v1,c)
	abc=(v2v3)(v4)(0)   (v2v3)(v4)(v1)
	abc=(v2)(v3v4)(0)   (v2)(v3v4)(v1)
	abc=(v2)(v3)(v4)    (v2)(v3)(v1v4)    (v2)(v3)(v4v1)
	abc=(0)(v2v3v4)(0)  (0)(v2v3v4)(v1)
	abc=(0)(v2v3)(v4)   (0)(v2v3)(v1v4)   (0)(v2v3)(v4v1)
abc=(0)(v2)(v3v4)   (0)(v2)(v1v3v4)   (0)(v2)(v3v1v4)   (0)(v2)(v3v4v1)

$$=\sum_{v_2\cdots v_r=abc, b,c\ne \emptyset} swap(Q)\bigl((sh(v_1,a)c\bigr)swapM(b\rfloor)
		+\sum_{v_2\cdots v_r=abc, c\ne \emptyset} swap(Q)(ac)swapM\bigl(sh(v_1,b)\rfloor\bigr)$$
		$$+\sum_{v_2\cdots v_r=abc,b\ne \emptyset} swap(Q)\bigl(a\,sh(v_1,c)\bigr)swapM(b\rfloor)
		-\sum_{v_2\cdots v_r=abc,b\ne \emptyset} swap(Q)\bigl(sh(v_1,a)c\bigr)swapM(\lfloor b)$$
		$$-\sum_{v_2\cdots v_r=abc, a\ne \emptyset} swap(Q)(ac)swapM\bigl(\lfloor sh(v_1,b)\bigr)
		-\sum_{v_2\cdots v_r=abc, a,b\ne \emptyset} swap(Q)\bigl(a\,sh(v_1,c)\bigr)swapM(\lfloor b).\eqno(A.4)$$
		Because $swapM$ is alternal, we can eliminate the second and fifth terms.
		Indeed, this follows from alternality when $b\ne\emptyset$, and in
		the case $b=\emptyset$, the second and fifth terms reduce to
		$$\sum_{v_2\cdots v_r=abc, c\ne \emptyset} swap(Q)(ac)swapM(sh(v_1,b)\rfloor)
		-\sum_{v_2\cdots v_r=abc, a\ne \emptyset} swap(Q)(ac)swapM(\lfloor sh(v_1,b))$$
		$$=\sum_{v_2\cdots v_r=ac, c\ne \emptyset} swap(Q)(ac)swapM(v_1\rfloor)
		-\sum_{v_2\cdots v_r=ac, a\ne \emptyset} swap(Q)(ac)swapM(\lfloor v_1)$$
		$$=\sum_{i=1}^{r-1} swap(Q)(v_2,\ldots,v_r)swapM(v_1-v_{i+1})
		-\sum_{i=2}^r swap(Q)(v_2,\ldots,v_r)swapM(v_1-v_i),\eqno(A.5)$$
		which is zero.

		Thus we can drop the second and fifth terms from (A.4), and it reduces to
		$$(A.4)=\sum_{v_2\cdots v_r=abc, b,c\ne \emptyset} swap(Q)\bigl(sh(v_1,a)c\bigr)swapM(b\rfloor)
		+\sum_{v_2\cdots v_r=abc,b\ne \emptyset} swap(Q)\bigl(a\,sh(v_1,c)\bigr)swapM(b\rfloor)$$
		$$-\sum_{v_2\cdots v_r=abc,b\ne \emptyset} swap(Q)\bigl(sh(v_1,a)c\bigr)swapM(\lfloor b)
	-\sum_{v_2\cdots v_r=abc, a,b\ne \emptyset} swap(Q)\bigl(a\,sh(v_1,c)\bigr)swapM(\lfloor b),\eqno(A.6)$$
	which we rewrite as

	$$=\sum_{v_2\cdots v_r=abc, b,c\ne \emptyset} 
	\Biggl(swap(Q)(sh(v_1,ac))swapM(b-c_f)+swap(Q)(av_1c)swapM(b-v_1)\Biggr)$$
$$-\sum_{v_2\cdots v_r=abc,a,b\ne \emptyset} \Biggl(swap(Q)(sh(v_1,ac))swapM(b-a_{last})
		+swap(Q)(av_1c)swapM(b-v_1)\Biggr)$$
	$$+\sum_{i=1}^{r-1} swap(Q)(v_2,\ldots,v_i,v_1)swapM(v_{i+1}-v_1,\ldots,v_r-v_1)$$
	$$-\sum_{i=2}^r swap(Q)(v_1,v_{i+1},\ldots,v_r)swapM(v_2-v_1,\ldots,v_i-v_1)$$
	which after some agreeable cancellation becomes
$$=\sum_{v_2\cdots v_r=abc, a,b,c\ne \emptyset} swap(Q)(sh(v_1,ac))\Biggl(swapM(b-c_f)
		-swapM(b-a_{last})\Biggr)$$
	$$+\sum_{j=2}^{r-1}
	swap(Q)(sh((v_1),(v_{j+1},\ldots,v_r))swapM(v_2-v_{j+1},\ldots,v_j-v_{j+1})$$
			$$-\sum_{i=2}^{r-1}swap(Q)(sh((v_1),(v_2,\ldots,v_i))swapM(v_{i+1}-v_i,\ldots,v_r-v_i),$$
				or more succinctly,
				$$=\sum_{1\le i<j\le r-1} swap(Q)(sh((v_1),(v_2,\ldots,v_i,v_{j+1},\ldots,v_r)))
				swapM(v_{i+1}-v_{j+1},\ldots,v_j-v_{j+1})$$
				$$-\sum_{2\le i<j\le r} swap(Q)(sh((v_1),(v_2,\ldots,v_i,v_{j+1},\ldots,v_r)))
				swapM(v_{i+1}-v_i,\ldots,v_j-v_i).$$

		Now we use the assumption that $dar\,Q$ satisfies the strict relations, and 
				apply (A.1) to simplify the shuffle terms in $swap(Q)$: 
				$$=-\sum_{1\le i<j\le r-1} 
				swap(Q)(v_2-v_1,\ldots,v_i-v_1,v_{j+1}-v_1,\ldots,v_r-v_1,-v_1)
				swapM(v_{i+1}-v_{j+1},\ldots,v_j-v_{j+1})$$
				$$+\sum_{2\le i<j\le r} 
				swap(Q)(v_2-v_1,\ldots,v_i-v_1,v_{j+1}-v_1,\ldots,v_r-v_1,-v_1)
				swapM(v_{i+1}-v_i,\ldots,v_j-v_i).$$
				We make two further modifications that are useful for the proof of the desired equality
				(A.1).  Since $M$ is assumed bialternal, it is $push$-invariant, so
				we change nothing by applying $push$ to the $swapM$ terms in the 
				second line:
				$$=-\sum_{1\le i<j\le r-1} 
swap(Q)(v_2-v_1,\ldots,v_i-v_1,v_{j+1}-v_1,\ldots,v_r-v_1,-v_1)
	swapM(v_{i+1}-v_{j+1},\ldots,v_j-v_{j+1})$$
$$+\sum_{2\le i<j\le r} swap(Q)(v_2-v_1,\ldots,v_i-v_1,v_{j+1}-v_1,\ldots,v_r-v_1,-v_1)
	swapM(v_i-v_j,v_{i+1}-v_j,\ldots,v_{j-1}-v_j).$$
	Finally, we reindex the second sum by $i\mapsto i-1$ and $j\mapsto j-1$:
	$$=-\sum_{1\le i<j\le r-1} 
swap(Q)(v_2-v_1,\ldots,v_i-v_1,v_{j+1}-v_1,\ldots,v_r-v_1,-v_1)
	swapM(v_{i+1}-v_{j+1},\ldots,v_j-v_{j+1})$$
$$+\sum_{1\le i<j\le r-1} swap(Q)(v_2-v_1,\ldots,v_{i+1}-v_1,v_{j+2}-v_1,\ldots,v_r-v_1,-v_1)
	swapM(v_{i+1}-v_{j+1},\ldots,v_{j}-v_{j+1}),$$
	which allows us to combine the two sums as
	$$=\sum_{1\le i<j\le r-1} 
	swapM(v_{i+1}-v_{j+1},\ldots,v_j-v_{j+1})\times$$
	$$\Biggl(swap(Q)(v_2-v_1,\ldots,v_{i+1}-v_1,v_{j+2}-v_1,\ldots,v_r-v_1,-v_1)-
			swap(Q)(v_2-v_1,\ldots,v_i-v_1,v_{j+1}-v_1,\ldots,v_r-v_1,-v_1)\Biggr).$$
	Finally, we separate out the term where $i=1$, to obtain the
	following expression for the first term of (A.3):
		$$\bigl(arit(swapM)\cdot swap(Q)\bigr)\bigl(sh((v_1),(v_2,\ldots,v_r))\bigr)
		   =\sum_{j=2}^{r-1} swapM(v_2-v_{j+1},\ldots,v_j-v_{j+1})\times$$
		   $$\Biggl(swap(Q)(v_2-v_1,v_{j+2}-v_1,\ldots,v_r-v_1,-v_1)-
				   swap(Q)(v_{j+1}-v_1,\ldots,v_r-v_1,-v_1)\Biggr)$$
		   $$+\sum_{2\le i<j\le r-1} 
		   swapM(v_{i+1}-v_{j+1},\ldots,v_j-v_{j+1})\times$$
		   $$\Biggl(swap(Q)(v_2-v_1,\ldots,v_{i+1}-v_1,v_{j+2}-v_1,\ldots,v_r-v_1,-v_1)-
				   swap(Q)(v_2-v_1,\ldots,v_i-v_1,v_{j+1}-v_1,\ldots,v_r-v_1,-v_1)\Biggr).\eqno(A.7)$$
		   \vskip .2cm
		   \noindent {\bf Left-hand side: second term of (A.3).} 
	For the second term of (A.3), we have
$$\Bigl(mu\bigl(swap(Q),swapM\bigr)\Bigr)
	\bigl(sh\bigl((v_1),(v_2,\ldots,v_r)\bigr)\bigr)\qquad\qquad\qquad\qquad\qquad
	\qquad\qquad\qquad$$
	$$=\sum_{v_2\cdots v_r=ab} \Biggl(swap(Q)\bigl(sh(v_1,a)\bigr)swapM(b)+
			swap(Q)(a)swapM\bigl(sh(v_1,b)\bigr)\Biggr).\eqno(A.8)$$

	The second term of the summand disappears when $a=0$ since $swap(Q)\in \overline{ARI}$.
	Because $M$ is bialternal, the second term also disappears for all $b$ except 
	$b=\emptyset$.  Thus with $a=(v_2,\ldots,v_i)$ and $b=(v_{i+1},\ldots,v_r)$,
	we can rewrite (A.8) as
$$=\sum_{i=2}^{r-1} swap(Q)\bigl(sh((v_1),(v_2,\ldots,v_i))\bigr)
	swapM(v_{i+1},\ldots,v_r)$$
	$$+swap(Q)(v_1)swapM(v_2,\ldots,v_r)+swap(Q)(v_2,\ldots,v_r)swapM(v_1).$$

	Since $Q$ satisfies the strict Fay relations, we can use (A.1) in the sum to obtain
$$\Bigl(mu\bigl(swap(Q),swapM\bigr)\Bigr)
	\bigl(sh\bigl((v_1),(v_2,\ldots,v_r)\bigr)\bigr)=
	\qquad \qquad \qquad \qquad \qquad \qquad \qquad \qquad \qquad \qquad
	\qquad\qquad$$
	$$-\sum_{i=2}^{r-1}swap(Q)(v_2-v_1,\ldots,v_i-v_1,-v_1)swapM(v_{i+1},\ldots,v_r)\qquad\qquad\qquad\qquad\qquad$$
	$$+ swap(Q)(v_1)swapM(v_2,\ldots,v_r)+swap(Q)(v_2,\ldots,v_r)swapM(v_1).\eqno(A.9)$$
	\vskip .1cm
	\noindent {\bf Left-hand side: third term of (A.3).} 
	We first write
$$mu(M,Q)\bigr)(u_1,\ldots,u_r)=\sum_{i=1}^{r-1} M(u_1,\ldots,u_i)
	Q(u_{i+1},\ldots,u_r),$$
	so
	$$swap\bigl(mu(M,Q)\bigr)(v_1,\ldots,v_r)=
	\sum_{i=1}^{r-1} M(v_r,v_{r-1}-v_r,\ldots,v_{r-i+1}-v_{r-i+2})Q(v_{r-i}-v_{r-i+1},\ldots, v_1-v_2)$$
$$=\sum_{i=1}^{r-1} swapM(v_{r-i+1},v_{r-i+2},\ldots,v_{r-1},v_r)
	swap(Q)(v_1-v_{r-i+1},v_2-v_{r-i+1},\ldots,v_{r-i}-v_{r-i+1})$$

	Thus,
	$$swap\bigl(mu(M,Q)\bigr)\bigl(sh((v_1),(v_2,\ldots,v_r))\bigr)
	\qquad \qquad \qquad \qquad \qquad \qquad \qquad \qquad \qquad \qquad $$
$$=\sum_{i=1}^{r-1} \Biggl(swapM(v_{r-i+1},\ldots,v_r)
		swap(Q)\bigl(sh((v_1-v_{r-i+1}),(v_2-v_{r-i+1},\ldots,v_{r-i}-v_{r-i+1}))\bigr)$$
		$$+swapM(v_1,v_{r-i+2},\ldots,v_r)swapQ(v_2-v_1,\ldots,v_{r-i+1}-v_1)\Biggr)$$
$$+\sum_{i=2}^{r-1} swapM(v_{r-i+2}\cdot sh((v_1),(v_{r-i+3},\ldots,v_r))
		swap(Q)(v_2-v_{r-i+2},\ldots,v_{r-i+1}-v_{r-i+2}).\eqno(A.10)$$
		Since $Q$ satisfies the strict Fay relations, we can use the equality
		$$swap(Q)\bigl(sh((v_1-v_{r-i+1}),(v_2-v_{r-i+1},\ldots,v_{r-i}-v_{r-i+1}))\bigr)=
		-swap(Q)\bigl(v_2-v_1,\ldots,v_{r-i}-v_1,v_{r-i+1}-v_1)$$
		in the first line.
		Furthermore, since $swapM$ is alternal, we can use the fact that
		$$swapM(v_{r-i+2}\cdot sh((v_1),(v_{r-i+3},\ldots,v_r)\bigr)=
			-swapM(v_1,v_{r-i+2},\ldots,v_r)$$
			in the third line.  This gives
			$$swap\bigl(mu(M,Q)\bigr)\bigl(sh((v_1),(v_2,\ldots,v_r))\bigr)
			\qquad \qquad \qquad \qquad \qquad \qquad \qquad \qquad \qquad \qquad $$
			$$=\sum_{i=1}^{r-1} \Biggl(-swapM(v_{r-i+1},\ldots,v_r)
				swap(Q)(v_2-v_1,\ldots,v_{r-i+1}-v_1)$$
				$$+swapM(v_1,v_{r-i+2},\ldots,v_r)swapQ(v_2-v_1,\ldots,v_{r-i+1}-v_1)\Biggr)$$
			$$-\sum_{i=2}^{r-1} swapM(v_1,v_{r-i+2},\ldots,v_r)
			swap(Q)(v_2-v_{r-i+2},\ldots,v_{r-i+1}-v_{r-i+2})\Biggr)$$
		$$=\sum_{i=2}^{r-1} swapM(v_1,v_{r-i+2},\ldots,v_r)\Biggl(swapQ(v_2-v_1,\ldots,v_{r-i+1}-v_1)
			-swap(Q)(v_2-v_{r-i+2},\ldots,v_{r-i+1}-v_{r-i+2})\Biggr)$$
		$$-\sum_{i=1}^{r-1} swapM(v_{r-i+1},\ldots,v_r)
		swap(Q)(v_2-v_1,\ldots,v_{r-i+1}-v_1)$$
$$+swapM(v_1)swap(Q)(v_2-v_1,\ldots,v_r-v_1).$$
Finally, we reindex the second sum with $i\mapsto r-i$ and separate out the
term for $i=1$:
$$=\sum_{i=2}^{r-1} swapM(v_1,v_{r-i+2},\ldots,v_r)\Biggl(swapQ(v_2-v_1,\ldots,v_{r-i+1}-v_1)
		-swap(Q)(v_2-v_{r-i+2},\ldots,v_{r-i+1}-v_{r-i+2})\Biggr)$$
$$-\sum_{i=2}^{r-1} swapM(v_{i+1},\ldots,v_r)
	swap(Q)(v_2-v_1,\ldots,v_{i+1}-v_1)$$
$$+swapM(v_1)swap(Q)(v_2-v_1,\ldots,v_r-v_1)
	-swapM(v_2,\ldots,v_r)swap(Q)(v_2-v_1).\eqno(A.11)$$

	\vskip .1cm
	\noindent {\bf Left-hand side of (A.1): total.}  The three terms of (A.3) calculated
	in (A.7), (A.9) and (A.11) combine as in (A.3) (with negative signs for the
			first two terms) to yield an explicit expression for the left-hand side of 
	the desired equality (A.1).  
	$$swap\bigl(arat(M)\cdot Q\bigr)\bigl(sh((v_1),(v_2,\ldots,v_r))\bigr)=
	\qquad\qquad \qquad\qquad \qquad\qquad \qquad\qquad$$
	$$\sum_{2\le i<j\le r-1} 
	swapM(v_{i+1}-v_{j+1},\ldots,v_j-v_{j+1})\times$$
$$\Biggl(swap(Q)(v_2-v_1,\ldots,v_i-v_1,v_{j+1}-v_1,\ldots,v_r-v_1,-v_1)
		-swap(Q)(v_2-v_1,\ldots,v_{i+1}-v_1,v_{j+2}-v_1,\ldots,v_r-v_1,-v_1)\Biggr)$$
	$$+\sum_{j=2}^{r-1} swapM(v_2-v_{j+1},\ldots,v_j-v_{j+1})\times$$
	$$\Biggl(swap(Q)(v_{j+1}-v_1,\ldots,v_r-v_1,-v_1)-swap(Q)(v_2-v_1,v_{j+2}-v_1,\ldots,v_r-v_1,-v_1)\Biggr)$$
	$$+\sum_{i=2}^{r-1}swap(Q)(v_2-v_1,\ldots,v_i-v_1,-v_1)swapM(v_{i+1},\ldots,v_r)\qquad\qquad\qquad\qquad\qquad$$
	$$- swap(Q)(v_1)swapM(v_2,\ldots,v_r)-swap(Q)(v_2,\ldots,v_r)swapM(v_1)$$
$$+\sum_{i=2}^{r-1} swapM(v_1,v_{r-i+2},\ldots,v_r)\Biggl(swapQ(v_2-v_1,\ldots,v_{r-i+1}-v_1)
		-swap(Q)(v_2-v_{r-i+2},\ldots,v_{r-i+1}-v_{r-i+2})\Biggr)$$
$$-\sum_{i=2}^{r-1} swapM(v_{i+1},\ldots,v_r)
	swap(Q)(v_2-v_1,\ldots,v_{i+1}-v_1)$$
$$+swapM(v_1)swap(Q)(v_2-v_1,\ldots,v_r-v_1)
	-swapM(v_2,\ldots,v_r)swap(Q)(v_2-v_1).\eqno(A.12)$$
	The third and fifth sums combine to form
	$$=\sum_{2\le i<j\le r-1} 
	swapM(v_{i+1}-v_{j+1},\ldots,v_j-v_{j+1})\times$$
$$\Biggl(swap(Q)(v_2-v_1,\ldots,v_i-v_1,v_{j+1}-v_1,\ldots,v_r-v_1,-v_1)
		-swap(Q)(v_2-v_1,\ldots,v_{i+1}-v_1,v_{j+2}-v_1,\ldots,v_r-v_1,-v_1)\Biggr)$$
	$$+\sum_{j=2}^{r-1} swapM(v_2-v_{j+1},\ldots,v_j-v_{j+1})\times$$
	$$\Biggl(swap(Q)(v_{j+1}-v_1,\ldots,v_r-v_1,-v_1)-swap(Q)(v_2-v_1,v_{j+2}-v_1,\ldots,v_r-v_1,-v_1)\Biggr)$$
	$$+\sum_{i=2}^{r-1}swapM(v_{i+1},\ldots,v_r)\times$$
	$$\Biggl(swap(Q)(v_2-v_1,\ldots,v_i-v_1,-v_1)-swap(Q)(v_2-v_1,\ldots,v_{i+1}-v_1)\Biggr)$$
	$$- swap(Q)(v_1)swapM(v_2,\ldots,v_r)-swap(Q)(v_2,\ldots,v_r)swapM(v_1)$$
$$+\sum_{i=2}^{r-1} swapM(v_1,v_{r-i+2},\ldots,v_r)\Biggl(swapQ(v_2-v_1,\ldots,v_{r-i+1}-v_1)
		-swap(Q)(v_2-v_{r-i+2},\ldots,v_{r-i+1}-v_{r-i+2})\Biggr)$$
$$+swapM(v_1)swap(Q)(v_2-v_1,\ldots,v_r-v_1)
	-swapM(v_2,\ldots,v_r)swap(Q)(v_2-v_1).\eqno(A.13)$$
	The newly formed third sum can now be considered as the $j=r$ term of the first sum, with the
	convention $v_{r+1}=0$.
	$$=\sum_{2\le i<j\le r} 
	swapM(v_{i+1}-v_{j+1},\ldots,v_j-v_{j+1})\times$$
$$\Biggl(swap(Q)(v_2-v_1,\ldots,v_i-v_1,v_{j+1}-v_1,\ldots,v_r-v_1,-v_1)
		-swap(Q)(v_2-v_1,\ldots,v_{i+1}-v_1,v_{j+2}-v_1,\ldots,v_r-v_1,-v_1)\Biggr)$$
	$$+\sum_{j=2}^{r-1} swapM(v_2-v_{j+1},\ldots,v_j-v_{j+1})\times$$
	$$\Biggl(swap(Q)(v_{j+1}-v_1,\ldots,v_r-v_1,-v_1)-swap(Q)(v_2-v_1,v_{j+2}-v_1,\ldots,v_r-v_1,-v_1)\Biggr)$$
	$$- swap(Q)(v_1)swapM(v_2,\ldots,v_r)-swap(Q)(v_2,\ldots,v_r)swapM(v_1)$$
$$+\sum_{i=2}^{r-1} swapM(v_1,v_{r-i+2},\ldots,v_r)\Biggl(swapQ(v_2-v_1,\ldots,v_{r-i+1}-v_1)
		-swap(Q)(v_2-v_{r-i+2},\ldots,v_{r-i+1}-v_{r-i+2})\Biggr)$$
$$+swapM(v_1)swap(Q)(v_2-v_1,\ldots,v_r-v_1)
	-swapM(v_2,\ldots,v_r)swap(Q)(v_2-v_1).\eqno(A.14)$$

	The terms
	$$swapM(v_1)swap(Q)(v_2-v_1,\ldots,v_r-v_1) -swap(Q)(v_2,\ldots,v_r)swapM(v_1)$$
	can be considered as the term $i=1$ term of the third sum,
	whereas using $-swap(Q)(v_1)=swap(Q)(-v_1)$ thanks to the oddness of $swap(Q)$ in depth $1$, 
	the terms
	$$swap(Q)(-v_1)swapM(v_2,\ldots,v_r) -swapM(v_2,\ldots,v_r)swap(Q)(v_2-v_1)$$
	can be considered as the $j=r$ term of the second sum:
	$$=\sum_{2\le i<j\le r} 
	swapM(v_{i+1}-v_{j+1},\ldots,v_j-v_{j+1})\times$$
$$\Biggl(swap(Q)(v_2-v_1,\ldots,v_i-v_1,v_{j+1}-v_1,\ldots,v_r-v_1,-v_1)
		-swap(Q)(v_2-v_1,\ldots,v_{i+1}-v_1,v_{j+2}-v_1,\ldots,v_r-v_1,-v_1)\Biggr)$$
	$$+\sum_{j=2}^r swapM(v_2-v_{j+1},\ldots,v_j-v_{j+1})\times$$
	$$\Biggl(swap(Q)(v_{j+1}-v_1,\ldots,v_r-v_1,-v_1)-swap(Q)(v_2-v_1,v_{j+2}-v_1,\ldots,v_r-v_1,-v_1)\Biggr)$$
$$+\sum_{i=1}^{r-1} swapM(v_1,v_{r-i+2},\ldots,v_r)\Biggl(swapQ(v_2-v_1,\ldots,v_{r-i+1}-v_1)
		-swap(Q)(v_2-v_{r-i+2},\ldots,v_{r-i+1}-v_{r-i+2})\Biggr)\eqno(A.15)$$

	Finally, since $M$ is bialternal, $swapM$ is $push$-invariant, so we can apply the $push$ to the $swapM$
	terms in the third sum without changing the value, obtaining
	$$=\sum_{2\le i<j\le r} 
	swapM(v_{i+1}-v_{j+1},\ldots,v_j-v_{j+1})\times$$
$$\Biggl(swap(Q)(v_2-v_1,\ldots,v_i-v_1,v_{j+1}-v_1,\ldots,v_r-v_1,-v_1)
		-swap(Q)(v_2-v_1,\ldots,v_{i+1}-v_1,v_{j+2}-v_1,\ldots,v_r-v_1,-v_1)\Biggr)$$
	$$+\sum_{j=2}^r swapM(v_2-v_{j+1},\ldots,v_j-v_{j+1})\times$$
	$$\Biggl(swap(Q)(v_{j+1}-v_1,\ldots,v_r-v_1,-v_1)-swap(Q)(v_2-v_1,v_{j+2}-v_1,\ldots,v_r-v_1,-v_1)\Biggr)$$
$$+\sum_{i=1}^{r-1} swapM(v_{r-i+2}-v_1,\ldots,v_r-v_1,-v_1)\Biggl(swapQ(v_2-v_1,\ldots,v_{r-i+1}-v_1)
		-swap(Q)(v_2-v_{r-i+2},\ldots,v_{r-i+1}-v_{r-i+2})\Biggr)\eqno(A.16)$$

	This concludes our computation of the left-hand side of the desired equality (A.1).
	We will now compute the right-hand side of (A.1), in four steps.
	\vskip .2cm
	\noindent {\bf Step 1. $\Bigl(arat(M)\cdot Q\Bigr)(u_1,\ldots,u_r)$.} For $w=(u_1,\ldots,u_r)$, we have 
	$$\Bigl(arat(M)\cdot Q\Bigr)(u_1,\ldots,u_r)=\sum_{w=abc} 
	Q(a\rceil c)M(b)-Q(a\lceil c)M(b).$$
	Writing $a=(u_1,\ldots,u_i)$, $b=(u_{i+1},\ldots,u_j)$, $c=(u_{j+1},\ldots,u_r)$, this becomes
	$$\Bigl(arat(M)\cdot Q\Bigr)(u_1,\ldots,u_r)=
	\sum_{1\le i<j\le r-1} M(u_{i+1},\ldots,u_j)\times$$
$$\Biggl(Q(u_1,\ldots,u_{i-1},u_i+u_{i+1}+\cdots+u_j,u_{j+1}, \ldots,u_r)
		-Q(u_1,\ldots,u_i,u_{i+1}+\cdots+u_{j+1},u_{j+2},\ldots, u_r)\Biggr)$$
$$+\sum_{1\le i\le r-1} M(u_1,\ldots,u_i)
	\Biggl(Q(u_{i+1},\ldots,u_r)-Q(u_1+\cdots+u_{i+1},u_{i+2},\ldots, u_r)\Biggr)$$
$$+\sum_{1\le i\le r-1} M(u_{i+1},\ldots,u_r)
	\Biggl(Q(u_1,\ldots,u_{i-1},u_i+u_{i+1}+\cdots+u_r) -Q(u_1,\ldots,u_i)\Biggr),$$
	where the first sum contains the decompositions $w=abc$ in which $a,b,c\ne \emptyset$, the
	second sum contains the terms with $a=\emptyset$ and the third those with $c=\emptyset$.
	\vskip .2cm
	\noindent {\bf Step 2. $swap\Bigl(arat(M)\cdot Q\Bigr)(v_1,\ldots,v_r)$.} We apply the swap to obtain
	$$swap\Bigl(arat(M)\cdot Q\Bigr)(v_1,\ldots,v_r)=
	\sum_{1\le i<j\le r-1} M(v_{r-i}-v_{r-i+1},\ldots,v_{r-j+1}-v_{r-j+2})\times$$
	$$\Biggl(
			Q(v_r,v_{r-1}-v_r,\ldots,v_{r-i+2}-v_{r-i+3},v_{r-j+1}-v_{r-i+2},v_{r-j}-v_{r-j+1},\ldots, v_1-v_2)
			$$
			$$\qquad\qquad 
			-Q(v_r,v_{r-1}-v_r,\ldots,v_{r-i+1}-v_{r-i+2},v_{r-j}-v_{r-i+1},v_{r-j-1}-v_{r-j}, \ldots,v_1-v_2)
			\Biggr)$$
$$+\sum_{1\le i\le r-1} M(v_r,v_{r-1}-v_r,\ldots,v_{r-i+1}-v_{r-i+2})
	\Biggl(Q(v_{r-i}-v_{r-i+1}, \ldots,v_1-v_2)-
			Q(v_{r-i},v_{r-i-1}-v_{r-i},\ldots, v_1-v_2)\Biggr)$$
$$+\sum_{1\le i\le r-1} M(v_{r-i}-v_{r-i+1},\ldots,v_1-v_2)
	\Biggl(Q(v_r,v_{r-1}-v_r,\ldots,v_{r-i+2}-v_{r-i+3},v_1-v_{r-i+2})$$
			$$\qquad\qquad\qquad\qquad\qquad\qquad\qquad\qquad -Q(v_r,v_{r-1}-v_r,\ldots,v_{r-i+1}-v_{r-i+2})\Biggr).$$
	\noindent {\bf Step 3. $swap\bigl(arat(M)\cdot Q\bigr)(v_1,\ldots,v_r)$ in terms
		of $swapM$ and $swap(Q)$.}
		$$=\sum_{1\le i<j\le r-1} swapM(v_{r-j+1}-v_{r-i+1},\ldots,v_{r-i}-v_{r-i+1})\times$$
		$$\Biggl(swap(Q)(v_1,\ldots,v_{r-j+1},v_{r-i+2},\ldots,v_r)-
				swap(Q)(v_1,\ldots,v_{r-j},v_{r-i+1},\ldots,v_r)\Biggr)$$
$$+\sum_{1\le i\le r-1} swapM(v_{r-i+1},v_{r-i+2},\ldots,v_r)
	\Biggl(swap(Q)(v_1-v_{r-i+1},\ldots,v_{r-i}-v_{r-i+1})-
			swap(Q)(v_1,\ldots,v_{r-i})\Biggr)$$
	$$+\sum_{1\le i\le r-1} swapM(v_1-v_{r-i+1},v_2-v_{r-i+1},\ldots,v_{r-i}-v_{r-i+1})\times$$
	$$\qquad\qquad\qquad\qquad\qquad \Biggl(swap(Q)(v_1,v_{r-i+2},\ldots,v_r))-swap(Q)(v_{r-i+1},\ldots,v_r)\Biggr).$$
	\vskip .3cm
	\noindent {\bf Step 4. The inverse $push$.}  We have
	$$swap\Bigl(arat(M)\cdot Q\Bigr)(v_2-v_1,\ldots,v_r-v_1,-v_1)=$$
	$$=\sum_{1\le i<j\le r-1} swapM(v_{r-j+2}-v_{r-i+2},\ldots,v_{r-i+1}-v_{r-i+2})\times$$
$$\Biggl(swap(Q)(v_2-v_1,\ldots,v_{r-j+2}-v_1,v_{r-i+3}-v_1,\ldots,v_r-v_1,-v_1)
		\qquad\qquad\qquad\qquad$$
		$$-swap(Q)(v_2-v_1,\ldots,v_{r-j+1}-v_1,v_{r-i+2}-v_1,\ldots,v_r-v_1,-v_1)\Biggr)$$
	$$+\sum_{1\le i\le r-1} swapM(v_{r-i+2}-v_1,\ldots,v_r-v_1,-v_1)\times$$
	$$\Biggl(swap(Q)(v_2-v_{r-i+2},\ldots,v_{r-i+1}-v_{r-i+2})-
			swap(Q)(v_2-v_1,\ldots,v_{r-i+1}-v_1)\Biggr)$$
	$$+\sum_{1\le i\le r-1} swapM(v_2-v_{r-i+2},\ldots,v_{r-i+1}-v_{r-i+2})\times$$
	$$\Biggl(swap(Q)(v_2-v_1,v_{r-i+3}-v_1,\ldots,v_r-v_1,-v_1))-
	swap(Q)(v_{r-i+2}-v_1,\ldots,v_r-v_1,-v_1)\Biggr).$$
	We reindex the double sum by $j\mapsto r-i+1$, $i\mapsto r-j+1$, and we reindex the
	final sum by $j=r-i+1$:
	$$=\sum_{2\le i<j\le r} swapM(v_{i+1}-v_{j+1},\ldots,v_j-v_{j+1})\times$$
$$\Biggl(swap(Q)(v_2-v_1,\ldots,v_{i+1}-v_1,v_{j+2}-v_1,\ldots,v_r-v_1,-v_1)
		\qquad\qquad\qquad\qquad$$
		$$-swap(Q)(v_2-v_1,\ldots,v_i-v_1,v_{j+1}-v_1,\ldots,v_r-v_1,-v_1)\Biggr)$$
	$$+\sum_{i=1}^{r-1} swapM(v_{r-i+2}-v_1,\ldots,v_r-v_1,-v_1)\times$$
	$$\Biggl(swap(Q)(v_2-v_{r-i+2},\ldots,v_{r-i+1}-v_{r-i+2})-
			swap(Q)(v_2-v_1,\ldots,v_{r-i+1}-v_1)\Biggr)$$
	$$+\sum_{j=2}^r swapM(v_2-v_{j+1},\ldots,v_{j}-v_{j+1})\times$$
	$$\Biggl(swap(Q)(v_2-v_1,v_{j+2}-v_1,\ldots,v_r-v_1,-v_1))-
	swap(Q)(v_{j+1}-v_1,\ldots,v_r-v_1,-v_1)\Biggr).\eqno(A.17)$$
	This is thus the expression for $push^{-1}.swap\bigl(arat(M)\cdot Q\bigr)$.  Since the right-hand
	side of the desired equality (A.1) is $-push^{-1}.swap\bigl(arat(M)\cdot Q\bigr)$, we conclude the
	proof by noting that the negative of (A.17) is equal to the left-hand side as given 
	in (A.16).\hfill{$\diamondsuit$}
	\vskip .3cm
	\noindent {\bf Porism to Lemma 2.9.} {\it Let $M$ be a bialternal mould with
		constant correction $C=(c_r)_{r\ge 2}$, and let $Q\in ARI$ be odd in depth 1 and satisfy the
			strict Fay relations. Then $arat(M)\cdot Q$ satisfies the Fay relations for the constants
			$(c_r)_{r\ge 2}$.}
			\vskip .3cm
			If $M$ is  non-strictly bialternal, then we need to make some modifications in the computations at the point where
			we claim that the shuffle terms in swap $M$ are zero.  In particular, we must rewrite (A.6) by recalculating the
			second and fifth terms of (A.4) in the case where 
			$$swapM\bigl(sh(v_1,b)\bigr)=-rc_r$$
			for $|b|=r$.

	The second and fifth terms of (A.4) are given by
$$=\sum_{v_2\cdots v_r=abc, c\ne \emptyset} swap(Q)(ac)swapM(sh(v_1,b)\rfloor)
	-\sum_{v_2\cdots v_r=abc, a\ne \emptyset} swap(Q)(ac)swapM(\lfloor sh(v_1,b))$$
	which we separate into cases as
$$+\sum_{v_2\cdots v_r=abc, a,b,c\ne \emptyset} swap(Q)(ac)swapM(sh(v_1,b)\rfloor)
	-\sum_{v_2\cdots v_r=abc, a,b,c=\emptyset} swap(Q)(ac)swapM(\lfloor sh(v_1,b))$$
$$+\sum_{v_2\cdots v_r=abc, a=\emptyset,b,c\ne \emptyset} swap(Q)(ac)swapM(sh(v_1,b)\rfloor)
	-\sum_{v_2\cdots v_r=abc, a,b\ne \emptyset,c=\emptyset} swap(Q)(ac)swapM(\lfloor sh(v_1,b))$$
$$+\sum_{v_2\cdots v_r=abc, b=\emptyset,a,c\ne \emptyset} swap(Q)(ac)swapM(sh(v_1,b)\rfloor)
	-\sum_{v_2\cdots v_r=abc, a,c\ne \emptyset,b=\emptyset} swap(Q)(ac)swapM(\lfloor sh(v_1,b))$$
$$+\sum_{v_2\cdots v_r=abc, a,b=\emptyset,c\ne \emptyset} swap(Q)(ac)swapM(sh(v_1,b)\rfloor)
	-\sum_{v_2\cdots v_r=abc, a\ne \emptyset,b,c=\emptyset} swap(Q)(ac)swapM(\lfloor sh(v_1,b)).$$
	In the first two terms, we have
	$$swapM(sh(v_1,b)\rfloor)=swapM(\lfloor sh(v_1,b))=-c_{|b|+1},$$
	so these two terms sum to zero.  We rewrite the remaining terms as
$$\sum_{i=2}^{r-1} swap(Q)(c)swapM\bigl(sh((v_1),(v_2,\ldots,v_i))\rfloor\bigr)
	-\sum_{i=2}^{r-1} swap(Q)(v_2,\ldots,v_i)swapM\bigl(\lfloor sh((v_1),(v_{i+1},\ldots,v_r)))$$
$$+\sum_{i=2}^{r-1} swap(Q)(v_2,\ldots,v_r)swapM(v_1-v_{i+1})
	-\sum_{i=2}^{r-1} swap(Q)(v_2,\ldots,v_r)swapM(v_1-v_i)$$
$$+swap(Q)(v_2,\ldots,v_r)swapM(v_1-v_2)
	-swap(Q)(v_2,\ldots,v_r)swapM(v_1-v_r).$$
	The last four terms sum to zero, leaving
$$\sum_{i=2}^{r-1} swap(Q)(v_{i+1},\ldots,v_r)swapM\bigl(sh((v_1),(v_2,\ldots,v_i)\rfloor\bigr)
		-\sum_{i=2}^{r-1} swap(Q)(v_2,\ldots,v_i)swapM\bigl(\lfloor sh((v_1),(v_{i+1},\ldots,v_r)))$$
		$$=\sum_{i=2}^{r-1} \Bigl((r-i+1)c_{r-i+1}swap(Q)(v_1,\ldots,v_i)-ic_iswap(Q)(v_{i+1},\ldots,v_r)\Bigr).\eqno(A.18)$$
		This is the term that must be added to (A.7) to compute the first term of (A.3) when $M$ is not strictly bialternal.

		For the second term of (A.3), we start with (A.8):
		$$=\sum_{v_2\cdots v_r=ab} \Biggl(swap(Q)\bigl(sh(v_1,a)\bigr)swapM(b)+
			swap(Q)(a)swapM\bigl(sh(v_1,b)\bigr)\Biggr)$$
		$$=swap(Q)(v_1)swapM(v_2,\ldots,v_r)+\sum_{i=2}^{r-1} swap(Q)\bigl(sh((v_1),(v_2,\ldots,v_i))\bigr)swapM(v_{i+1},\ldots,v_r)$$
		$$\qquad\qquad +\sum_{i=2}^{r-1} swap(Q)(v_2,\ldots,v_i)swapM\bigl(sh((v_1),(v_{i+1},\ldots,v_r))\bigr)
		+swap(Q)(v_2,\ldots,v_r)swapM(v_1),$$
		but instead of eliminating the shuffle term in $swapM$ when $b\ne \emptyset$, we replace it by 
		$-c_{|b|+1}$, obtaining
		$$=swap(Q)(v_1)swapM(v_2,\ldots,v_r)+\sum_{i=2}^{r-1} swap(Q)\bigl(sh((v_1),(v_2,\ldots,v_i))\bigr)swapM(v_{i+1},\ldots,v_r)\qquad\qquad\qquad$$
		$$\qquad\qquad -\sum_{i=2}^{r-1} c_{r-i+1}swap(Q)(v_2,\ldots,v_i)+
		swap(Q)(v_2,\ldots,v_r)swapM(v_1).$$
		We now use the strict Fay relation on $swap(Q)$ to substitute 
		$$swap(Q)\bigl(sh((v_1),(v_2,\ldots,v_i))\bigr)=-swap(Q)(v_2-v_1,\ldots,v_i-v_1,-v_1),$$
		obtaining
		$$=swap(Q)(v_1)swapM(v_2,\ldots,v_r)-\sum_{i=2}^{r-1} swap(Q)(v_2-v_1,\ldots,v_i-v_1,-v_1)swapM(v_{i+1},\ldots,v_r)\qquad\qquad\qquad$$
		$$\qquad\qquad -\sum_{i=2}^{r-1} c_{r-i+1}swap(Q)(v_2,\ldots,v_i)+
swap(Q)(v_2,\ldots,v_r)swapM(v_1).$$
Comparing with (A.9), we see that we need to add the correction term
$$-\sum_{i=2}^{r-1} c_{r-i+1}swap(Q)(v_2,\ldots,v_i)\eqno(A.19)$$
to the second term of (A.3).

We now consider the third term of (A.3), as computed in (A.11).
We start with (A.10) with the substitution
$$swap(Q)\bigl(sh((v_1-v_{r-i+1}),(v_2-v_{r-i+1},\ldots,v_{r-i}-v_{r-i+1}))\bigr)=
-swap(Q)(v_2-v_1,\ldots,v_{r-i}-v_1,r_{r-i+1}-v_1),$$
	giving
$$\sum_{i=1}^{r-1} \Biggl(-swapM(v_{r-i+1},\ldots,v_r)
		swap(Q)(v_2-v_1,\ldots,v_{r-i}-v_1,r_{r-i+1}-v_1)$$
		$$+swapM(v_1,v_{r-i+2},\ldots,v_r)swapQ(v_2-v_1,\ldots,v_{r-i+1}-v_1)\Biggr)$$
$$+\sum_{i=2}^{r-1} swapM(v_{r-i+2}\cdot sh((v_1),(v_{r-i+3},\ldots,v_r))
		swap(Q)(v_2-v_{r-i+2},\ldots,v_{r-i+1}-v_{r-i+2}),$$
		and then we substitute
		$$swapM(v_{r-i+2}\cdot sh((v_1),(v_{r-i+3},\ldots,v_r))=-swapM(v_1,v_{r-i+2},\ldots,v_r)-c_i$$ 
			to obtain
			$$\sum_{i=1}^{r-1} \Biggl(-swapM(v_{r-i+1},\ldots,v_r)
				swap(Q)(v_2-v_1,\ldots,v_{r-i}-v_1,r_{r-i+1}-v_1)$$
				$$+swapM(v_1,v_{r-i+2},\ldots,v_r)swapQ(v_2-v_1,\ldots,v_{r-i+1}-v_1)\Biggr)$$
			$$-\sum_{i=2}^{r-1} swapM(v_1,v_{r-i+2},\ldots,v_r) swap(Q)(v_2-v_{r-i+2},\ldots,v_{r-i+1}-v_{r-i+2})$$
			$$-\sum_{i=2}^{r-1}c_iswap(Q)(v_2-v_{r-i+2},\ldots,v_{r-i+1}-v_{r-i+2}),$$
			so that 
			$$-\sum_{i=2}^{r-1}c_iswap(Q)(v_2-v_{r-i+2},\ldots,v_{r-i+1}-v_{r-i+2})\eqno(A.20)$$
			is the term that must be added to (A.11) in the non-strict case.

			To obtain the complete correction, we add up (A.18), (A.19) and (A.20): 
			$$=\sum_{i=2}^{r-1} (r-i+1)c_{r-i+1}swap(Q)(v_2,\ldots,v_i) -\sum_{i=2}^{r-1} ic_iswap(Q)(v_{i+1},\ldots,v_r)$$
			$$-\sum_{i=2}^{r-1} c_{r-i+1}swap(Q)(v_2,\ldots,v_i) 
			-\sum_{i=2}^{r-1}c_iswap(Q)(v_2-v_{r-i+2},\ldots,v_{r-i+1}-v_{r-i+2}).$$
			Reindexing the first and third sums by $i\mapsto r-i+1$, this gives
			$$=\sum_{i=2}^{r-1} c_i\Bigl(swap(Q)(v_1,\ldots,v_{r-i+1}) -swap(Q)(v_{i+1},\ldots,v_r)$$
				$$-swap(Q)(v_2,\ldots,v_{r-i+1})-swap(Q)(v_2-v_{r-i+2},\ldots,v_{r-i+1}-v_{r-i+2})\Bigr).$$
			\vfill\eject
			We now work out the general expressions for terms 1, 2 and 3 of (A.3) when $Q$ satisfies Fay for a family
		$(\epsilon_r)_{r\ge 2}$ and $M$ has a correction term $C=(c_r)_{r\ge 2}$.
		\vskip .2cm
		\noindent {\bf Term 1 of (A.3).}
		The full expression of this term is given in (A.4).

	We saw that the first, third, fourth and sixth terms of (A.4) add up to
$$\sum_{1\le i<j\le r-1} swap(Q)(sh((v_1),(v_2,\ldots,v_i,v_{j+1},\ldots,v_r)))
	swapM(v_{i+1}-v_{j+1},\ldots,v_j-v_{j+1})$$
$$-\sum_{2\le i<j\le r} swap(Q)(sh((v_1),(v_2,\ldots,v_i,v_{j+1},\ldots,v_r)))
	swapM(v_{i+1}-v_i,\ldots,v_j-v_i).$$

	We saw in (A.18) that the second and fifth terms add up to 
	$$\sum_{i=2}^{r-1} \Bigl((r-i+1)c_{r-i+1}swap(Q)(v_2,\ldots,v_i)-ic_iswap(Q)(v_{i+1},\ldots,v_r)\Bigr)$$

	Therefore (A.4) adds up to
$$(A.4)=\sum_{1\le i<j\le r-1} swap(Q)(sh((v_1),(v_2,\ldots,v_i,v_{j+1},\ldots,v_r)))
	swapM(v_{i+1}-v_{j+1},\ldots,v_j-v_{j+1})$$
$$-\sum_{2\le i<j\le r} swap(Q)(sh((v_1),(v_2,\ldots,v_i,v_{j+1},\ldots,v_r)))
	swapM(v_{i+1}-v_i,\ldots,v_j-v_i)$$
	$$+\sum_{i=2}^{r-1} \Bigl((r-i+1)c_{r-i+1}swap(Q)(v_2,\ldots,v_i)-ic_iswap(Q)(v_{i+1},\ldots,v_r)\Bigr)$$
	which, separating out the $i=1$ term from the first double sum, becomes
$$=\sum_{j=2}^{r-1} swap(Q)(sh((v_1),(v_{j+1},\ldots,v_r)))
	swapM(v_2-v_{j+1},\ldots,v_j-v_{j+1})$$
$$+\sum_{2\le i<j\le r-1} swap(Q)(sh((v_1),(v_2,\ldots,v_i,v_{j+1},\ldots,v_r)))
	swapM(v_{i+1}-v_{j+1},\ldots,v_j-v_{j+1})$$
$$-\sum_{2\le i<j\le r} swap(Q)(sh((v_1),(v_2,\ldots,v_i,v_{j+1},\ldots,v_r)))
	swapM(v_{i+1}-v_i,\ldots,v_j-v_i)$$
	$$+\sum_{i=2}^{r-1} \Bigl((r-i+1)c_{r-i+1}swap(Q)(v_2,\ldots,v_i)-ic_iswap(Q)(v_{i+1},\ldots,v_r)\Bigr).$$
	Now we will use the assumption that $Q$ satisfies the Fay relation for the family $(\epsilon_r)_{r\ge 2}$ in the
	form:
	$$swap(Q)\bigl(sh\bigl((v_1), (v_2,\ldots,v_r)\bigr)\bigr)=
	-swap(Q)(v_2-v_1,\ldots,v_r-v_1,-v_1)-r\epsilon_rv_2$$
	to write:
	$$\cases{swap(Q)\bigl(sh((v_1),(v_{j+1},\ldots,v_r))\bigr)=
		-swap(Q)(v_{j+1}-v_1,\ldots,v_r-v_1,-v_1)-(r-j+1)\epsilon_{r-j+1}v_{j+1}&
			$i=1$\cr
			swap(Q)\bigl(sh((v_1),(v_2,\ldots,v_i,v_{j+1},\ldots,v_r))\bigr)=\cr
			\qquad -swap(Q)(v_2-v_1,\ldots,v_i-v_1,v_{j+1}-v_1,\ldots,v_r-v_1,-v_1)
			-(r-j+i)\epsilon_{r-j+i}v_2&$i\ge 2$.}$$
			Using this, we rewrite:
			$$(A.4)=-\sum_{j=2}^{r-1} 
			\Bigl(swap(Q)(v_{j+1}-v_1,\ldots,v_r-v_1,-v_1)+(r-j+1)\epsilon_{r-j+1}v_{j+1}\Bigr)swapM(v_2-v_{j+1},\ldots,v_j-v_{j+1})$$
			$$-\sum_{2\le i<j\le r-1} 
\Bigl(swap(Q)(v_2-v_1,\ldots,v_i-v_1,v_{j+1}-v_1,\ldots,v_r-v_1,-v_1)
		+(r-j+i)\epsilon_{r-j+i}v_2\Bigr)
	swapM(v_{i+1}-v_{j+1},\ldots,v_j-v_{j+1})$$
	$$-\sum_{2\le i<j\le r} 
\Bigl(swap(Q)(v_2-v_1,\ldots,v_i-v_1,v_{j+1}-v_1,\ldots,v_r-v_1,-v_1)
		+(r-j+i)\epsilon_{r-j+i}v_2\Bigr)
	swapM(v_{i+1}-v_i,\ldots,v_j-v_i)$$
	$$+\sum_{i=2}^{r-1} \Bigl((r-i+1)c_{r-i+1}swap(Q)(v_2,\ldots,v_i)-ic_iswap(Q)(v_{i+1},\ldots,v_r)\Bigr),$$
	and then separating out the correction terms:
	$$(A.4)=-\sum_{j=2}^{r-1} 
	swap(Q)(v_{j+1}-v_1,\ldots,v_r-v_1,-v_1)swapM(v_2-v_{j+1},\ldots,v_j-v_{j+1})$$
	$$-\sum_{2\le i<j\le r-1} 
swap(Q)(v_2-v_1,\ldots,v_i-v_1,v_{j+1}-v_1,\ldots,v_r-v_1,-v_1)
	swapM(v_{i+1}-v_{j+1},\ldots,v_j-v_{j+1})$$
	$$-\sum_{2\le i<j\le r} 
swap(Q)(v_2-v_1,\ldots,v_i-v_1,v_{j+1}-v_1,\ldots,v_r-v_1,-v_1)
	swapM(v_{i+1}-v_i,\ldots,v_j-v_i)$$
	$$+\sum_{i=2}^{r-1} \Bigl((r-i+1)c_{r-i+1}swap(Q)(v_2,\ldots,v_i)-ic_iswap(Q)(v_{i+1},\ldots,v_r)\Bigr),$$
	$$-\sum_{j=2}^{r-1} 
	(r-j+1)\epsilon_{r-j+1}v_{j+1}\,swapM(v_2-v_{j+1},\ldots,v_j-v_{j+1})$$
	$$-\sum_{2\le i<j\le r-1} 
	(r-j+i)\epsilon_{r-j+i}v_2\, swapM(v_{i+1}-v_{j+1},\ldots,v_j-v_{j+1})$$
	$$-\sum_{2\le i<j\le r} 
	(r-j+i)\epsilon_{r-j+i}v_2\, swapM(v_{i+1}-v_i,\ldots,v_j-v_i)$$

	\vskip .3cm
	\noindent {\bf Term 2 of (A.3).}
	We need to compute 
	$$(A.8)=\sum_{v_2\cdots v_r=ab} \Biggl(swap(Q)\bigl(sh(v_1,a)\bigr)swapM(b)+
			swap(Q)(a)swapM\bigl(sh(v_1,b)\bigr)\Biggr)$$
	$$=\sum_{i=2}^{r-1} swap(Q)\bigl(sh((v_1),(v_2,\ldots,v_i)))\bigr)swapM(v_{i+1},\ldots,v_r)$$
	$$\qquad\qquad\qquad +\sum_{i=2}^{r-1} swap(Q)(v_2,\ldots,v_i)swapM\bigl(sh((v_1),(v_{i+1},\ldots,v_r))\bigr)$$
	$$+swap(Q)(v_1)swapM(v_2,\ldots,v_r)+swap(Q)(v_2,\ldots,v_r)swapM(v_1)$$
	$$=-\sum_{i=2}^{r-1} swap(Q)(v_2-v_1,\ldots,v_i-v_1,-v_1)swapM(v_{i+1},\ldots,v_r)$$
	$$+swap(Q)(v_1)swapM(v_2,\ldots,v_r)+swap(Q)(v_2,\ldots,v_r)swapM(v_1)$$
$$-\sum_{i=2}^{r-1} i\epsilon_iv_2\,swapM(v_{i+1},\ldots,v_r)
	-\sum_{i=2}^{r-1} (r-i+1)c_{r-i+1} swap(Q)(v_2,\ldots,v_i)$$
	\vskip .1cm
	\noindent {\bf Term 3 of (A.3).}
	We need to compute 
$$(A.10)=\sum_{i=1}^{r-1} swapM(v_{r-i+1},\ldots,v_r)
	swap(Q)\bigl(sh((v_1-v_{r-i+1}),(v_2-v_{r-i+1},\ldots,v_{r-i}-v_{r-i+1}))\bigr)$$
	$$+\sum_{i=1}^{r-1} swapM(v_1,v_{r-i+2},\ldots,v_r)swapQ(v_2-v_1,\ldots,v_{r-i+1}-v_1)$$
$$+\sum_{i=2}^{r-1} swapM(v_{r-i+2}\cdot sh((v_1),(v_{r-i+3},\ldots,v_r))
		swap(Q)(v_2-v_{r-i+2},\ldots,v_{r-i+1}-v_{r-i+2})$$
		$$=-\sum_{i=1}^{r-1}swapM(v_{r-i+1},\ldots,v_r)swap(Q)(v_2-v_1,\ldots,v_{r-i+1}-v_1)$$
		$$+\sum_{i=1}^{r-1} swapM(v_1,v_{r-i+2},\ldots,v_r)swapQ(v_2-v_1,\ldots,v_{r-i+1}-v_1)$$
		$$-\sum_{i=2}^{r-1}swapM(v_1,v_{r-i+2},\ldots,v_r)swap(Q)(v_2-v_{r-i+2},\ldots,
			v_{r-i+1}-v_{r-i+2})$$
		$$-\sum_{i=1}^{r-1} (r-i)\epsilon_{r-i}(v_2-v_{r-i+1})swapM(v_{r-i+1},\ldots,v_r)$$
		$$-\sum_{i=2}^{r-1}ic_iswap(Q)(v_2-v_{r-i+2},\ldots,v_{r-i+1}-v_{r-i+2})$$

		\noindent {\bf Total.} In the end, the left-hand side of the Fay relation for $P=arat(M)\cdot Q$:
		$$swap(P)\bigl(sh((v_1),(v_2,\ldots,v_r))\bigr)+swap(P)(v_2-v_1,\ldots,v_r-v_1,-v_1)$$
		is equal to
		\centerline{Term 1}
		$$-\sum_{i=2}^{r-1} \Bigl((r-i+1)c_{r-i+1}swap(Q)(v_2,\ldots,v_i)-ic_iswap(Q)(v_{i+1},\ldots,v_r)\Bigr),$$
\centerline{Term 2}
$$+\sum_{i=2}^{r-1} (r-i+1)c_{r-i+1} swap(Q)(v_2,\ldots,v_i)$$
\centerline{Term 3}
$$-\sum_{i=2}^{r-1}ic_iswap(Q)(v_2-v_{r-i+2},\ldots,v_{r-i+1}-v_{r-i+2})$$
In the special case $Q=durM$, this becomes:
$$\sum_{i=2}^{r-1} (r-i+1)c_{r-i+1}v_1\,swapM(v_1,\ldots,v_i)$$
$$-\sum_{i=2}^{r-1} ic_iv_{i+1}\,swapM(v_{i+1},\ldots,v_r)$$
$$-\sum_{j=2}^{r-1} 
(r-j+1)\epsilon_{r-j+1}v_{j+1}\,swapM(v_2-v_{j+1},\ldots,v_j-v_{j+1})$$
$$-\sum_{2\le i<j\le r-1} 
(r-j+i)\epsilon_{r-j+i}v_2\, swapM(v_{i+1}-v_{j+1},\ldots,v_j-v_{j+1})$$
$$-\sum_{2\le i<j\le r} 
(r-j+i)\epsilon_{r-j+i}v_2\, swapM(v_{i+1}-v_i,\ldots,v_j-v_i)$$
$$-\sum_{i=2}^{r-1} i\epsilon_iv_2\,swapM(v_{i+1},\ldots,v_r)$$
$$-\sum_{i=2}^{r-1} (r-i+1)c_{r-i+1}v_2\, swapM(v_2,\ldots,v_i)$$
$$-\sum_{i=1}^{r-1} (r-i)\epsilon_{r-i}v_2\,swapM(v_{r-i+1},\ldots,v_r)$$
$$+\sum_{i=1}^{r-1} (r-i)\epsilon_{r-i}v_{r-i+1}\,swapM(v_{r-i+1},\ldots,v_r)$$
$$-\sum_{i=2}^{r-1}ic_iv_2\,swapM(v_2-v_{r-i+2},\ldots,v_{r-i+1}-v_{r-i+2})$$
$$+\sum_{i=2}^{r-1}ic_iv_{r-i+2}\,swapM(v_2-v_{r-i+2},\ldots,v_{r-i+1}-v_{r-i+2})$$
Putting lines 6 and 8 together at the bottom and renumbering line 8 with
$r-i\mapsto i$:
$$\sum_{i=2}^{r-1} (r-i+1)c_{r-i+1}v_1\,swapM(v_1,\ldots,v_i)$$
$$-\sum_{i=2}^{r-1} ic_iv_{i+1}\,swapM(v_{i+1},\ldots,v_r)$$
$$-\sum_{j=2}^{r-1} (r-j+1)\epsilon_{r-j+1}v_{j+1}\,swapM(v_2-v_{j+1},\ldots,v_j-v_{j+1})$$
$$-\sum_{2\le i<j\le r-1} (r-j+i)\epsilon_{r-j+i}v_2\, swapM(v_{i+1}-v_{j+1},\ldots,v_j-v_{j+1})$$
$$-\sum_{2\le i<j\le r} (r-j+i)\epsilon_{r-j+i}v_2\, swapM(v_{i+1}-v_i,\ldots,v_j-v_i)$$
$$-\sum_{i=2}^{r-1} (r-i+1)c_{r-i+1}v_2\, swapM(v_2,\ldots,v_i)$$
$$+\sum_{i=1}^{r-1} (r-i)\epsilon_{r-i}v_{r-i+1}\,swapM(v_{r-i+1},\ldots,v_r)$$
$$-\sum_{i=2}^{r-1}ic_iv_2\,swapM(v_2-v_{r-i+2},\ldots,v_{r-i+1}-v_{r-i+2})$$
$$+\sum_{i=2}^{r-1}ic_iv_{r-i+2}\,swapM(v_2-v_{r-i+2},\ldots,v_{r-i+1}-v_{r-i+2})$$
$$-\sum_{i=2}^{r-1} i\epsilon_iv_2\,swapM(v_{i+1},\ldots,v_r)$$
$$-\sum_{i=1}^{r-1} i\epsilon_{i}v_2\,swapM(v_{i+1},\ldots,v_r)$$

\vfill\eject
\centerline{\bf Some extras}
\vskip .5cm
The mould that Matthes worked with is $\hat{A'}_\tau$ which is defined 
as follows.  In this part we don't work modulo $2\pi i$.
\vskip .1cm
\noindent $\bullet$ We start with 
$$A_\tau(a,b)=\Phi_{KZ}(t_{01},t_{12})e^{2\pi it_{01}}\Phi_{KZ}(t_{01},t_{12})^{-1}$$
$$A'_\tau(a,b)=A_\tau(b,a)$$
$$\hat{A}'_\tau(a,b)=e^{-\pi i[a,b]}A'_\tau(a,b)$$
This $\hat{A'}_\tau$ is the power series that satisfies the Fay relations given by 
Matthes:
$$\cases{\F\bigl(\hat{A}'(\tau)\bigr)(u_1,u_2)=3\zeta(2)\cr
	\F\bigl(\hat{A}'(\tau)\bigr)(u_1,u_2,u_3)=\zeta(2)\bigl(dar^{-1}\hat{A}'(\tau)(u_2)
			+dar^{-1}\hat{A}'(\tau)(u_3)\bigr).}$$
	The Fay relation in depth 4 satisfied by $\hat{A}'(\tau)$ is very close to his but
not quite the same (why?)

	Since 
	$$A'_\tau(a,b)=e^{\pi i[a,b]}\hat{A}'_\tau(a,b),$$
	and
	$$ma(e^{\pi i[a,b]})(u_1,\ldots,u_r)=(-1)^r{{(\pi i)^r}\over{r!}}u_1\cdots u_r,$$
	we have
	$$\eqalign{A'(\tau)(u_1,\ldots,u_r)&=mu\Bigl(ma(e^{\pi i[a,b]}),\hat{A}'(\tau)\Bigr)(u_1,\ldots,u_r)\cr
		&=\sum_{i=0}^r ma(e^{\pi i[a,b]}(u_1,\ldots,u_i)\hat{A}'(\tau)(u_{i+1},\ldots,u_r)\cr
				&=\sum_{i=0}^r (-1)^i{{(\pi i)^i}\over{i!}}u_1\ldots u_i\cdot \hat{A}'(\tau)(u_{i+1},\ldots,u_r)\cr
				.}$$

				Since we're working mod $2\pi i$, we can drop the $\hat{}$ and just use
				$A_\tau(a,b)$ and $A'_\tau(a,b)=A_\tau(b,a)$.  The relation between them is as 
				follows.  We have
				$$A_\tau(a,b)=Dgarit(E_\tau)\cdot t_{01}=t_{01}(E_\tau,F_\tau),$$ 
				where $F_\tau$ is the partner of $E_\tau$. Therefore
				$$A'_\tau(a,b)=t_{01}\bigl(E_\tau(b,a),F_\tau(b,a)\bigr)
				=Dgarit(E_\tau(b,a))\cdot t_{01}.$$

				t01=-a+(1/2)[b,a]+(1/12).[b,[b,a]]
				Dgarit(E).t01 = -E+(1/2)[E,F]+(1/12).[F,[F,E]]+...

				Let E'=E(b,a).  Then what is F', the partner of E'?

				Knowing that [E(a,b),F(a,b)]=[a,b], we have
				[E(b,a),F(b,a)]=[b,a]  so
				[E(b,a),-F(b,a)]=[a,b].
				Thus the partner of E'=E(b,a) is F'=-F(b,a)

		Dgarit(E').t01=-E'+(1/2)[E',F']+(1/12).[F',[F',E']]+...
	=-E(b,a)+(1/2)[E(b,a),-F(b,a)]+(1/12).[F(b,a),[F(b,a),E(b,a)]]+...
	=[-E(a,b)+(1/2)[E(a,b),-F(a,b)]+(1/12).[F(a,b),[F(a,b),E(a,b)]]+...](b,a)
	=[-E-(1/2)[E(a,b),F(a,b)]+(1/12).[F,[F,E)]]+...](b,a)
	=[Dgarit(E).t01-[a,b]](b,a) 
=[a,b] + [Dgarit(E).t01](b,a)

	Therefore we have
	Dgarit(E).t01=Dgarit(E').t01-[a,b]

	The best thing is to go back to the Lie version and study
	$log_{Dari}(E')$ and $Darit(log_{Dari}(E'))\cdot t01$.

	\vfill\eject
	In this final section, we show how to use the result of Theorem 2.5 to 
	determine the Fay relations satisfied by the group-like elliptic associator
	$\bar{A}(\tau)$.  The key point is that if the Lie-like version satisfies
	a family of Fay relations, then the group-like version does as well.  The
	actual computation of the right-hand side correction terms for the group-like
	family can be done step by step.  We end with the depth 2, 3 and 4 cases,
	comparing them to the Fay relations first computed by Matthes.
	\vskip .2cm
	\noindent {\bf Theorem 2.6.} {\it Let $M$ be a a bialternal mould with
		constant correction $C$, and $Q$ a mould satisfying a family of Fay
			relations.  Then $Dgarit\bigl(exp_{Dari}(\Delta M)\bigr)\cdot \exp(Q)$
			satisfies a family of Fay relations.}
			\vskip .2cm
			\noindent {\bf Proof.} The proof follows easily from Proposition 2.3. Firstly,
			we have
			$$Dgarit(\Delta M)\cdot Q=\sum_{n\ge 0} {{1}\over{n!}}Darit(\Delta M)^n\cdot Q.$$
			By Proposition 2.3, $Q_1:=Darit(\Delta M)\cdot Q$ then also satisfies a family
			of Fay relations. But then $Darit(\Delta M)^2\cdot Q=Darit(\Delta M)\cdot Q_1$
			also satisfies a family of Fay relations and so forth.
			\vskip.3cm

			\vfill\eject
			I'm keeping separately here the complete proof of the Fay relations
			$$\F\bigl(arat(M)\cdot Q\bigr)$$
			in the special case where $M$ is strictly bialternal and $Q$ satisfies the
			strict Fay relations.
			\noindent {\bf A.1. Computation of the first term 
				$swap(P')\Bigl(sh\bigl((v_1),(v_2,\ldots,v_r)\bigr)\Bigr)$ of (A.1).} 
				\vskip .1cm
				\noindent We have 
				$$P=arat(M)\cdot Q=-arit(M)\cdot Q+lu(M,Q),$$
				so we compute
				$$\eqalign{swap(P)&=swap\bigl(arat(M)\cdot Q\bigr)\cr
					&=-swap\bigl(arit(M)\cdot Q\bigr)+swap\bigl(mu(M,Q)\bigr)-swap\bigl(mu(Q,M)\bigr)\cr
						&=-axit\bigl(swapM,-push(swapM)\bigr)\cdot swap(Q)-mu\bigl(swap(Q),swapM\bigr)\cr
						&\qquad+swap\bigl(mu(Q,M)\bigr)+swap\bigl(mu(M,Q)\bigr)-swap\bigl(mu(Q,M)\bigr)\cr
						&=-axit\bigl(swapM,-push(swapM)\bigr)\cdot swap(Q)-mu\bigl(swap(Q),swapM\bigr) +swap\bigl(mu(M,Q)\bigr)\cr
						&=-arit\bigl(swapM\bigr)\cdot swap(Q)-mu\bigl(swap(Q),swapM\bigr)
						+swap\bigl(mu(M,Q)\bigr),}\eqno(A.3)$$
						where the last equality follows since $M$ is bialternal and thus 
						$push$-invariant, and $axit(P,-P)=arit(P)$.
						In order to compute the leftmost term of (A.1), we will use (A.3) and
	separately compute the three terms
$$\cases{\Bigl(arit\bigl(swapM\bigr)\cdot swap(Q)\Bigr)
	\bigl(sh\bigl((v_1),(v_2,\ldots,v_r)\bigr)\bigr)\cr
		\Bigl(mu\bigl(swap(Q),swapM\bigr)\Bigr)
		\bigl(sh\bigl((v_1),(v_2,\ldots,v_r)\bigr)\bigr)\cr
		\Bigl(swap\bigl(mu(M,Q)\bigr)\Bigr)
		\bigl(sh\bigl((v_1),(v_2,\ldots,v_r)\bigr)\bigr).}\eqno(A.3)$$
		\noindent {\bf First term of (A.3).} We compute
		$$\Bigl(arit\bigl(swapM\bigr)\cdot swap(Q)\Bigr)\bigl(sh((v_1),(v_2,\ldots,v_r))\bigr)$$
		by adapting a method due to A.~Salerno (cf. [SS, appendix]), as follows.
		The definition of the $arit$-operator on moulds in $\overline{ARI}$ is given
		by the following formula.
		$$\Bigl(arit\bigl(swapM\bigr)\cdot swap(Q)\Bigr)(w)=\sum_{w=abc,c\ne 0}
		Q(ac)M(b\rfloor) - \sum_{w=abc, a\ne 0} Q(ac)M(\lfloor b).$$
	We thus have
	$$\Bigl(arit(swapM)\cdot swap(Q)\Bigr)
\bigl(sh((v_1),(v_2,\ldots,v_r))\bigr) 
	\qquad\qquad\qquad\qquad\qquad \qquad\qquad\qquad\qquad\qquad\qquad\qquad\qquad  $$
$$=\sum_{v_2\cdots v_r=abc, b,c\ne \emptyset} swap(Q)\bigl((sh(v_1,a)c\bigr)swapM(b\rfloor)
		+\sum_{v_2\cdots v_r=abc, c\ne \emptyset} swap(Q)(ac)swapM\bigl(sh(v_1,b)\rfloor\bigr)$$
		$$+\sum_{v_2\cdots v_r=abc,b\ne \emptyset} swap(Q)\bigl(a\,sh(v_1,c)\bigr)swapM(b\rfloor)
		-\sum_{v_2\cdots v_r=abc,b\ne \emptyset} swap(Q)\bigl(sh(v_1,a)c\bigr)swapM(\lfloor b)$$
		$$-\sum_{v_2\cdots v_r=abc, a\ne \emptyset} swap(Q)(ac)swapM\bigl(\lfloor sh(v_1,b)\bigr)
		-\sum_{v_2\cdots v_r=abc, a,b\ne \emptyset} swap(Q)\bigl(a\,sh(v_1,c)\bigr)swapM(\lfloor b).\eqno(A.4)$$
		Because $swapM$ is alternal, we can eliminate the second and fifth terms.
		Indeed, this follows from alternality when $b\ne\emptyset$, and in
		the case $b=\emptyset$, the second and fifth terms reduce to
		$$\sum_{v_2\cdots v_r=abc, c\ne \emptyset} swap(Q)(ac)swapM(sh(v_1,b)\rfloor)
		-\sum_{v_2\cdots v_r=abc, a\ne \emptyset} swap(Q)(ac)swapM(\lfloor sh(v_1,b))$$
		$$=\sum_{v_2\cdots v_r=ac, c\ne \emptyset} swap(Q)(ac)swapM(v_1\rfloor)
		-\sum_{v_2\cdots v_r=ac, a\ne \emptyset} swap(Q)(ac)swapM(\lfloor v_1)$$
		$$=\sum_{i=1}^{r-1} swap(Q)(v_2,\ldots,v_r)swapM(v_1-v_{i+1})
		-\sum_{i=2}^r swap(Q)(v_2,\ldots,v_r)swapM(v_1-v_i),\eqno(A.5)$$
		which is zero.

		Thus we can drop the second and fifth terms from (A.4), and it reduces to
		$$(A.4)=\sum_{v_2\cdots v_r=abc, b,c\ne \emptyset} swap(Q)\bigl(sh(v_1,a)c\bigr)swapM(b\rfloor)
		+\sum_{v_2\cdots v_r=abc,b\ne \emptyset} swap(Q)\bigl(a\,sh(v_1,c)\bigr)swapM(b\rfloor)$$
		$$-\sum_{v_2\cdots v_r=abc,b\ne \emptyset} swap(Q)\bigl(sh(v_1,a)c\bigr)swapM(\lfloor b)
	-\sum_{v_2\cdots v_r=abc, a,b\ne \emptyset} swap(Q)\bigl(a\,sh(v_1,c)\bigr)swapM(\lfloor b),\eqno(A.6)$$
	which we rewrite as

	$$=\sum_{v_2\cdots v_r=abc, b,c\ne \emptyset} 
	\Biggl(swap(Q)(sh(v_1,ac))swapM(b-c_f)+swap(Q)(av_1c)swapM(b-v_1)\Biggr)$$
$$-\sum_{v_2\cdots v_r=abc,a,b\ne \emptyset} \Biggl(swap(Q)(sh(v_1,ac))swapM(b-a_{last})
		+swap(Q)(av_1c)swapM(b-v_1)\Biggr)$$
	$$+\sum_{i=1}^{r-1} swap(Q)(v_2,\ldots,v_i,v_1)swapM(v_{i+1}-v_1,\ldots,v_r-v_1)$$
	$$-\sum_{i=2}^r swap(Q)(v_1,v_{i+1},\ldots,v_r)swapM(v_2-v_1,\ldots,v_i-v_1)$$
	which after some agreeable cancellation becomes
$$=\sum_{v_2\cdots v_r=abc, a,b,c\ne \emptyset} swap(Q)(sh(v_1,ac))\Biggl(swapM(b-c_f)
		-swapM(b-a_{last})\Biggr)$$
	$$+\sum_{j=2}^{r-1}
	swap(Q)(sh((v_1),(v_{j+1},\ldots,v_r))swapM(v_2-v_{j+1},\ldots,v_j-v_{j+1})$$
			$$-\sum_{i=2}^{r-1}swap(Q)(sh((v_1),(v_2,\ldots,v_i))swapM(v_{i+1}-v_i,\ldots,v_r-v_i),$$
				or more succinctly,
				$$=\sum_{1\le i<j\le r-1} swap(Q)(sh((v_1),(v_2,\ldots,v_i,v_{j+1},\ldots,v_r)))
				swapM(v_{i+1}-v_{j+1},\ldots,v_j-v_{j+1})$$
				$$-\sum_{2\le i<j\le r} swap(Q)(sh((v_1),(v_2,\ldots,v_i,v_{j+1},\ldots,v_r)))
				swapM(v_{i+1}-v_i,\ldots,v_j-v_i).$$

				Now we use the assumption that $dar\,Q$ satisfies the strict relations, and 
				apply (A.1) to simplify the shuffle terms in $swap(Q)$: 
				$$=-\sum_{1\le i<j\le r-1} 
				swap(Q)(v_2-v_1,\ldots,v_i-v_1,v_{j+1}-v_1,\ldots,v_r-v_1,-v_1)
				swapM(v_{i+1}-v_{j+1},\ldots,v_j-v_{j+1})$$
				$$+\sum_{2\le i<j\le r} 
				swap(Q)(v_2-v_1,\ldots,v_i-v_1,v_{j+1}-v_1,\ldots,v_r-v_1,-v_1)
				swapM(v_{i+1}-v_i,\ldots,v_j-v_i).$$
				We make two further modifications that are useful for the proof of the desired equality
				(A.1).  Since $M$ is assumed bialternal, it is $push$-invariant, so
				we change nothing by applying $push$ to the $swapM$ terms in the 
				second line:
				$$=-\sum_{1\le i<j\le r-1} 
swap(Q)(v_2-v_1,\ldots,v_i-v_1,v_{j+1}-v_1,\ldots,v_r-v_1,-v_1)
	swapM(v_{i+1}-v_{j+1},\ldots,v_j-v_{j+1})$$
$$+\sum_{2\le i<j\le r} swap(Q)(v_2-v_1,\ldots,v_i-v_1,v_{j+1}-v_1,\ldots,v_r-v_1,-v_1)
	swapM(v_i-v_j,v_{i+1}-v_j,\ldots,v_{j-1}-v_j).$$
	Finally, we reindex the second sum by $i\mapsto i-1$ and $j\mapsto j-1$:
	$$=-\sum_{1\le i<j\le r-1} 
swap(Q)(v_2-v_1,\ldots,v_i-v_1,v_{j+1}-v_1,\ldots,v_r-v_1,-v_1)
	swapM(v_{i+1}-v_{j+1},\ldots,v_j-v_{j+1})$$
$$+\sum_{1\le i<j\le r-1} swap(Q)(v_2-v_1,\ldots,v_{i+1}-v_1,v_{j+2}-v_1,\ldots,v_r-v_1,-v_1)
	swapM(v_{i+1}-v_{j+1},\ldots,v_{j}-v_{j+1}),$$
	which allows us to combine the two sums as
	$$=\sum_{1\le i<j\le r-1} 
	swapM(v_{i+1}-v_{j+1},\ldots,v_j-v_{j+1})\times$$
	$$\Biggl(swap(Q)(v_2-v_1,\ldots,v_{i+1}-v_1,v_{j+2}-v_1,\ldots,v_r-v_1,-v_1)-
			swap(Q)(v_2-v_1,\ldots,v_i-v_1,v_{j+1}-v_1,\ldots,v_r-v_1,-v_1)\Biggr).$$
	Finally, we separate out the term where $i=1$, to obtain the
	following expression for the first term of (A.3):
		$$\bigl(arit(swapM)\cdot swap(Q)\bigr)\bigl(sh((v_1),(v_2,\ldots,v_r))\bigr)
		   =\sum_{j=2}^{r-1} swapM(v_2-v_{j+1},\ldots,v_j-v_{j+1})\times$$
		   $$\Biggl(swap(Q)(v_2-v_1,v_{j+2}-v_1,\ldots,v_r-v_1,-v_1)-
				   swap(Q)(v_{j+1}-v_1,\ldots,v_r-v_1,-v_1)\Biggr)$$
		   $$+\sum_{2\le i<j\le r-1} 
		   swapM(v_{i+1}-v_{j+1},\ldots,v_j-v_{j+1})\times$$
		   $$\Biggl(swap(Q)(v_2-v_1,\ldots,v_{i+1}-v_1,v_{j+2}-v_1,\ldots,v_r-v_1,-v_1)-
				   swap(Q)(v_2-v_1,\ldots,v_i-v_1,v_{j+1}-v_1,\ldots,v_r-v_1,-v_1)\Biggr).\eqno(A.7)$$
		   \vskip .2cm
		   \noindent {\bf Second term of (A.3).} 
	For the second term of (A.3), we have
$$\Bigl(mu\bigl(swap(Q),swapM\bigr)\Bigr)
	\bigl(sh\bigl((v_1),(v_2,\ldots,v_r)\bigr)\bigr)\qquad\qquad\qquad\qquad\qquad
	\qquad\qquad\qquad$$
	$$=\sum_{v_2\cdots v_r=ab} \Biggl(swap(Q)\bigl(sh(v_1,a)\bigr)swapM(b)+
			swap(Q)(a)swapM\bigl(sh(v_1,b)\bigr)\Biggr).\eqno(A.8)$$

	The second term of the summand disappears when $a=0$ since $swap(Q)\in \overline{ARI}$.
	Because $M$ is bialternal, the second term also disappears for all $b$ except 
	$b=\emptyset$.  Thus with $a=(v_2,\ldots,v_i)$ and $b=(v_{i+1},\ldots,v_r)$,
	we can rewrite (A.8) as
$$=\sum_{i=2}^{r-1} swap(Q)\bigl(sh((v_1),(v_2,\ldots,v_i))\bigr)
	swapM(v_{i+1},\ldots,v_r)$$
	$$+swap(Q)(v_1)swapM(v_2,\ldots,v_r)+swap(Q)(v_2,\ldots,v_r)swapM(v_1).$$

	Since $Q$ satisfies the strict Fay relations, we can use (A.1) in the sum to obtain
$$\Bigl(mu\bigl(swap(Q),swapM\bigr)\Bigr)
	\bigl(sh\bigl((v_1),(v_2,\ldots,v_r)\bigr)\bigr)=
	\qquad \qquad \qquad \qquad \qquad \qquad \qquad \qquad \qquad \qquad
	\qquad\qquad$$
	$$-\sum_{i=2}^{r-1}swap(Q)(v_2-v_1,\ldots,v_i-v_1,-v_1)swapM(v_{i+1},\ldots,v_r)\qquad\qquad\qquad\qquad\qquad$$
	$$+ swap(Q)(v_1)swapM(v_2,\ldots,v_r)+swap(Q)(v_2,\ldots,v_r)swapM(v_1).\eqno(A.9)$$
	\vskip .1cm
	\noindent {\bf Third term of (A.3).} 
	We first write
$$mu(M,Q)\bigr)(u_1,\ldots,u_r)=\sum_{i=1}^{r-1} M(u_1,\ldots,u_i)
	Q(u_{i+1},\ldots,u_r),$$
	so
	$$swap\bigl(mu(M,Q)\bigr)(v_1,\ldots,v_r)=
	\sum_{i=1}^{r-1} M(v_r,v_{r-1}-v_r,\ldots,v_{r-i+1}-v_{r-i+2})Q(v_{r-i}-v_{r-i+1},\ldots, v_1-v_2)$$
$$=\sum_{i=1}^{r-1} swapM(v_{r-i+1},v_{r-i+2},\ldots,v_{r-1},v_r)
	swap(Q)(v_1-v_{r-i+1},v_2-v_{r-i+1},\ldots,v_{r-i}-v_{r-i+1})$$

	Thus,
	$$swap\bigl(mu(M,Q)\bigr)\bigl(sh((v_1),(v_2,\ldots,v_r))\bigr)
	\qquad \qquad \qquad \qquad \qquad \qquad \qquad \qquad \qquad \qquad $$
$$=\sum_{i=1}^{r-1} \Biggl(swapM(v_{r-i+1},\ldots,v_r)
		swap(Q)\bigl(sh((v_1-v_{r-i+1}),(v_2-v_{r-i+1},\ldots,v_{r-i}-v_{r-i+1}))\bigr)$$
		$$+swapM(v_1,v_{r-i+2},\ldots,v_r)swapQ(v_2-v_1,\ldots,v_{r-i+1}-v_1)\Biggr)$$
$$+\sum_{i=2}^{r-1} swapM(v_{r-i+2}\cdot sh((v_1),(v_{r-i+3},\ldots,v_r))
		swap(Q)(v_2-v_{r-i+2},\ldots,v_{r-i+1}-v_{r-i+2}).\eqno(A.10)$$
		Since $Q$ satisfies the strict Fay relations, we can use the equality
		$$swap(Q)\bigl(sh((v_1-v_{r-i+1}),(v_2-v_{r-i+1},\ldots,v_{r-i}-v_{r-i+1}))\bigr)=
		-swap(Q)\bigl(v_2-v_1,\ldots,v_{r-i}-v_1,v_{r-i+1}-v_1)$$
		in the first line.
		Furthermore, since $swapM$ is alternal, we can use the fact that
		$$swapM(v_{r-i+2}\cdot sh((v_1),(v_{r-i+3},\ldots,v_r)\bigr)=
			-swapM(v_1,v_{r-i+2},\ldots,v_r)$$
			in the third line.  This gives
			$$swap\bigl(mu(M,Q)\bigr)\bigl(sh((v_1),(v_2,\ldots,v_r))\bigr)
			\qquad \qquad \qquad \qquad \qquad \qquad \qquad \qquad \qquad \qquad $$
			$$=\sum_{i=1}^{r-1} \Biggl(-swapM(v_{r-i+1},\ldots,v_r)
				swap(Q)(v_2-v_1,\ldots,v_{r-i+1}-v_1)$$
				$$+swapM(v_1,v_{r-i+2},\ldots,v_r)swapQ(v_2-v_1,\ldots,v_{r-i+1}-v_1)\Biggr)$$
			$$-\sum_{i=2}^{r-1} swapM(v_1,v_{r-i+2},\ldots,v_r)
			swap(Q)(v_2-v_{r-i+2},\ldots,v_{r-i+1}-v_{r-i+2})\Biggr)$$
		$$=\sum_{i=2}^{r-1} swapM(v_1,v_{r-i+2},\ldots,v_r)\Biggl(swapQ(v_2-v_1,\ldots,v_{r-i+1}-v_1)
			-swap(Q)(v_2-v_{r-i+2},\ldots,v_{r-i+1}-v_{r-i+2})\Biggr)$$
		$$-\sum_{i=1}^{r-1} swapM(v_{r-i+1},\ldots,v_r)
		swap(Q)(v_2-v_1,\ldots,v_{r-i+1}-v_1)$$
$$+swapM(v_1)swap(Q)(v_2-v_1,\ldots,v_r-v_1).$$
Finally, we reindex the second sum with $i\mapsto r-i$ and separate out the
term for $i=1$:
$$=\sum_{i=2}^{r-1} swapM(v_1,v_{r-i+2},\ldots,v_r)\Biggl(swapQ(v_2-v_1,\ldots,v_{r-i+1}-v_1)
		-swap(Q)(v_2-v_{r-i+2},\ldots,v_{r-i+1}-v_{r-i+2})\Biggr)$$
$$-\sum_{i=2}^{r-1} swapM(v_{i+1},\ldots,v_r)
	swap(Q)(v_2-v_1,\ldots,v_{i+1}-v_1)$$
$$+swapM(v_1)swap(Q)(v_2-v_1,\ldots,v_r-v_1)
	-swapM(v_2,\ldots,v_r)swap(Q)(v_2-v_1).\eqno(A.11)$$

	\vskip .1cm
	\noindent {\bf Leftmost term of (A.1): total.}  The three terms of (A.3) calculated
	in (A.7), (A.9) and (A.11) combine as in (A.3) (with negative signs for the
			first two terms) to yield an explicit expression for the left-hand side of 
	the desired equality (A.1).  
	$$swap\bigl(arat(M)\cdot Q\bigr)\bigl(sh((v_1),(v_2,\ldots,v_r))\bigr)=
	\qquad\qquad \qquad\qquad \qquad\qquad \qquad\qquad$$
	$$\sum_{2\le i<j\le r-1} 
	swapM(v_{i+1}-v_{j+1},\ldots,v_j-v_{j+1})\times$$
$$\Biggl(swap(Q)(v_2-v_1,\ldots,v_i-v_1,v_{j+1}-v_1,\ldots,v_r-v_1,-v_1)
		-swap(Q)(v_2-v_1,\ldots,v_{i+1}-v_1,v_{j+2}-v_1,\ldots,v_r-v_1,-v_1)\Biggr)$$
	$$+\sum_{j=2}^{r-1} swapM(v_2-v_{j+1},\ldots,v_j-v_{j+1})\times$$
	$$\Biggl(swap(Q)(v_{j+1}-v_1,\ldots,v_r-v_1,-v_1)-swap(Q)(v_2-v_1,v_{j+2}-v_1,\ldots,v_r-v_1,-v_1)\Biggr)$$
	$$+\sum_{i=2}^{r-1}swap(Q)(v_2-v_1,\ldots,v_i-v_1,-v_1)swapM(v_{i+1},\ldots,v_r)\qquad\qquad\qquad\qquad\qquad$$
	$$- swap(Q)(v_1)swapM(v_2,\ldots,v_r)-swap(Q)(v_2,\ldots,v_r)swapM(v_1)$$
$$+\sum_{i=2}^{r-1} swapM(v_1,v_{r-i+2},\ldots,v_r)\Biggl(swapQ(v_2-v_1,\ldots,v_{r-i+1}-v_1)
		-swap(Q)(v_2-v_{r-i+2},\ldots,v_{r-i+1}-v_{r-i+2})\Biggr)$$
$$-\sum_{i=2}^{r-1} swapM(v_{i+1},\ldots,v_r)
	swap(Q)(v_2-v_1,\ldots,v_{i+1}-v_1)$$
$$+swapM(v_1)swap(Q)(v_2-v_1,\ldots,v_r-v_1)
	-swapM(v_2,\ldots,v_r)swap(Q)(v_2-v_1).\eqno(A.12)$$
	The third and fifth sums combine to form
	$$=\sum_{2\le i<j\le r-1} 
	swapM(v_{i+1}-v_{j+1},\ldots,v_j-v_{j+1})\times$$
$$\Biggl(swap(Q)(v_2-v_1,\ldots,v_i-v_1,v_{j+1}-v_1,\ldots,v_r-v_1,-v_1)
		-swap(Q)(v_2-v_1,\ldots,v_{i+1}-v_1,v_{j+2}-v_1,\ldots,v_r-v_1,-v_1)\Biggr)$$
	$$+\sum_{j=2}^{r-1} swapM(v_2-v_{j+1},\ldots,v_j-v_{j+1})\times$$
	$$\Biggl(swap(Q)(v_{j+1}-v_1,\ldots,v_r-v_1,-v_1)-swap(Q)(v_2-v_1,v_{j+2}-v_1,\ldots,v_r-v_1,-v_1)\Biggr)$$
	$$+\sum_{i=2}^{r-1}swapM(v_{i+1},\ldots,v_r)\times$$
	$$\Biggl(swap(Q)(v_2-v_1,\ldots,v_i-v_1,-v_1)-swap(Q)(v_2-v_1,\ldots,v_{i+1}-v_1)\Biggr)$$
	$$- swap(Q)(v_1)swapM(v_2,\ldots,v_r)-swap(Q)(v_2,\ldots,v_r)swapM(v_1)$$
$$+\sum_{i=2}^{r-1} swapM(v_1,v_{r-i+2},\ldots,v_r)\Biggl(swapQ(v_2-v_1,\ldots,v_{r-i+1}-v_1)
		-swap(Q)(v_2-v_{r-i+2},\ldots,v_{r-i+1}-v_{r-i+2})\Biggr)$$
$$+swapM(v_1)swap(Q)(v_2-v_1,\ldots,v_r-v_1)
	-swapM(v_2,\ldots,v_r)swap(Q)(v_2-v_1).\eqno(A.13)$$
	The newly formed third sum can now be considered as the $j=r$ term of the first sum, with the
	convention $v_{r+1}=0$.
	$$=\sum_{2\le i<j\le r} 
	swapM(v_{i+1}-v_{j+1},\ldots,v_j-v_{j+1})\times$$
$$\Biggl(swap(Q)(v_2-v_1,\ldots,v_i-v_1,v_{j+1}-v_1,\ldots,v_r-v_1,-v_1)
		-swap(Q)(v_2-v_1,\ldots,v_{i+1}-v_1,v_{j+2}-v_1,\ldots,v_r-v_1,-v_1)\Biggr)$$
	$$+\sum_{j=2}^{r-1} swapM(v_2-v_{j+1},\ldots,v_j-v_{j+1})\times$$
	$$\Biggl(swap(Q)(v_{j+1}-v_1,\ldots,v_r-v_1,-v_1)-swap(Q)(v_2-v_1,v_{j+2}-v_1,\ldots,v_r-v_1,-v_1)\Biggr)$$
	$$- swap(Q)(v_1)swapM(v_2,\ldots,v_r)-swap(Q)(v_2,\ldots,v_r)swapM(v_1)$$
$$+\sum_{i=2}^{r-1} swapM(v_1,v_{r-i+2},\ldots,v_r)\Biggl(swapQ(v_2-v_1,\ldots,v_{r-i+1}-v_1)
		-swap(Q)(v_2-v_{r-i+2},\ldots,v_{r-i+1}-v_{r-i+2})\Biggr)$$
$$+swapM(v_1)swap(Q)(v_2-v_1,\ldots,v_r-v_1)
	-swapM(v_2,\ldots,v_r)swap(Q)(v_2-v_1).\eqno(A.14)$$

	The terms
	$$swapM(v_1)swap(Q)(v_2-v_1,\ldots,v_r-v_1) -swap(Q)(v_2,\ldots,v_r)swapM(v_1)$$
	can be considered as the term $i=1$ term of the third sum,
	whereas using $-swap(Q)(v_1)=swap(Q)(-v_1)$ thanks to the oddness of $swap(Q)$ in depth $1$, 
	the terms
	$$swap(Q)(-v_1)swapM(v_2,\ldots,v_r) -swapM(v_2,\ldots,v_r)swap(Q)(v_2-v_1)$$
	can be considered as the $j=r$ term of the second sum:
	$$=\sum_{2\le i<j\le r} 
	swapM(v_{i+1}-v_{j+1},\ldots,v_j-v_{j+1})\times$$
$$\Biggl(swap(Q)(v_2-v_1,\ldots,v_i-v_1,v_{j+1}-v_1,\ldots,v_r-v_1,-v_1)
		-swap(Q)(v_2-v_1,\ldots,v_{i+1}-v_1,v_{j+2}-v_1,\ldots,v_r-v_1,-v_1)\Biggr)$$
	$$+\sum_{j=2}^r swapM(v_2-v_{j+1},\ldots,v_j-v_{j+1})\times$$
	$$\Biggl(swap(Q)(v_{j+1}-v_1,\ldots,v_r-v_1,-v_1)-swap(Q)(v_2-v_1,v_{j+2}-v_1,\ldots,v_r-v_1,-v_1)\Biggr)$$
$$+\sum_{i=1}^{r-1} swapM(v_1,v_{r-i+2},\ldots,v_r)\Biggl(swapQ(v_2-v_1,\ldots,v_{r-i+1}-v_1)
		-swap(Q)(v_2-v_{r-i+2},\ldots,v_{r-i+1}-v_{r-i+2})\Biggr)\eqno(A.15)$$

	Finally, since $M$ is bialternal, $swapM$ is $push$-invariant, so we can apply the $push$ to the $swapM$
	terms in the third sum without changing the value, obtaining
	$$=\sum_{2\le i<j\le r} 
	swapM(v_{i+1}-v_{j+1},\ldots,v_j-v_{j+1})\times$$
$$\Biggl(swap(Q)(v_2-v_1,\ldots,v_i-v_1,v_{j+1}-v_1,\ldots,v_r-v_1,-v_1)
		-swap(Q)(v_2-v_1,\ldots,v_{i+1}-v_1,v_{j+2}-v_1,\ldots,v_r-v_1,-v_1)\Biggr)$$
	$$+\sum_{j=2}^r swapM(v_2-v_{j+1},\ldots,v_j-v_{j+1})\times$$
	$$\Biggl(swap(Q)(v_{j+1}-v_1,\ldots,v_r-v_1,-v_1)-swap(Q)(v_2-v_1,v_{j+2}-v_1,\ldots,v_r-v_1,-v_1)\Biggr)$$
$$+\sum_{i=1}^{r-1} swapM(v_{r-i+2}-v_1,\ldots,v_r-v_1,-v_1)\Biggl(swapQ(v_2-v_1,\ldots,v_{r-i+1}-v_1)
		-swap(Q)(v_2-v_{r-i+2},\ldots,v_{r-i+1}-v_{r-i+2})\Biggr)\eqno(A.16)$$

	This concludes our computation of the left-hand side of the desired equality (A.1).
	We will now compute the right-hand side of (A.1), in four steps.
	\vskip .2cm
	\noindent {\bf Step 1. $\Bigl(arat(M)\cdot Q\Bigr)(u_1,\ldots,u_r)$.} For $w=(u_1,\ldots,u_r)$, we have 
	$$\Bigl(arat(M)\cdot Q\Bigr)(u_1,\ldots,u_r)=\sum_{w=abc} 
	Q(a\rceil c)M(b)-Q(a\lceil c)M(b).$$
	Writing $a=(u_1,\ldots,u_i)$, $b=(u_{i+1},\ldots,u_j)$, $c=(u_{j+1},\ldots,u_r)$, this becomes
	$$\Bigl(arat(M)\cdot Q\Bigr)(u_1,\ldots,u_r)=
	\sum_{1\le i<j\le r-1} M(u_{i+1},\ldots,u_j)\times$$
$$\Biggl(Q(u_1,\ldots,u_{i-1},u_i+u_{i+1}+\cdots+u_j,u_{j+1}, \ldots,u_r)
		-Q(u_1,\ldots,u_i,u_{i+1}+\cdots+u_{j+1},u_{j+2},\ldots, u_r)\Biggr)$$
$$+\sum_{1\le i\le r-1} M(u_1,\ldots,u_i)
	\Biggl(Q(u_{i+1},\ldots,u_r)-Q(u_1+\cdots+u_{i+1},u_{i+2},\ldots, u_r)\Biggr)$$
$$+\sum_{1\le i\le r-1} M(u_{i+1},\ldots,u_r)
	\Biggl(Q(u_1,\ldots,u_{i-1},u_i+u_{i+1}+\cdots+u_r) -Q(u_1,\ldots,u_i)\Biggr),$$
	where the first sum contains the decompositions $w=abc$ in which $a,b,c\ne \emptyset$, the
	second sum contains the terms with $a=\emptyset$ and the third those with $c=\emptyset$.
	\vskip .2cm
	\noindent {\bf Step 2. $swap\Bigl(arat(M)\cdot Q\Bigr)(v_1,\ldots,v_r)$.} We apply the swap to obtain
	$$swap\Bigl(arat(M)\cdot Q\Bigr)(v_1,\ldots,v_r)=
	\sum_{1\le i<j\le r-1} M(v_{r-i}-v_{r-i+1},\ldots,v_{r-j+1}-v_{r-j+2})\times$$
	$$\Biggl(
			Q(v_r,v_{r-1}-v_r,\ldots,v_{r-i+2}-v_{r-i+3},v_{r-j+1}-v_{r-i+2},v_{r-j}-v_{r-j+1},\ldots, v_1-v_2)
			$$
			$$\qquad\qquad 
			-Q(v_r,v_{r-1}-v_r,\ldots,v_{r-i+1}-v_{r-i+2},v_{r-j}-v_{r-i+1},v_{r-j-1}-v_{r-j}, \ldots,v_1-v_2)
			\Biggr)$$
$$+\sum_{1\le i\le r-1} M(v_r,v_{r-1}-v_r,\ldots,v_{r-i+1}-v_{r-i+2})
	\Biggl(Q(v_{r-i}-v_{r-i+1}, \ldots,v_1-v_2)-
			Q(v_{r-i},v_{r-i-1}-v_{r-i},\ldots, v_1-v_2)\Biggr)$$
$$+\sum_{1\le i\le r-1} M(v_{r-i}-v_{r-i+1},\ldots,v_1-v_2)
	\Biggl(Q(v_r,v_{r-1}-v_r,\ldots,v_{r-i+2}-v_{r-i+3},v_1-v_{r-i+2})$$
			$$\qquad\qquad\qquad\qquad\qquad\qquad\qquad\qquad -Q(v_r,v_{r-1}-v_r,\ldots,v_{r-i+1}-v_{r-i+2})\Biggr).$$
	\noindent {\bf Step 3. $swap\bigl(arat(M)\cdot Q\bigr)(v_1,\ldots,v_r)$ in terms
		of $swapM$ and $swap(Q)$.}
		$$=\sum_{1\le i<j\le r-1} swapM(v_{r-j+1}-v_{r-i+1},\ldots,v_{r-i}-v_{r-i+1})\times$$
		$$\Biggl(swap(Q)(v_1,\ldots,v_{r-j+1},v_{r-i+2},\ldots,v_r)-
				swap(Q)(v_1,\ldots,v_{r-j},v_{r-i+1},\ldots,v_r)\Biggr)$$
$$+\sum_{1\le i\le r-1} swapM(v_{r-i+1},v_{r-i+2},\ldots,v_r)
	\Biggl(swap(Q)(v_1-v_{r-i+1},\ldots,v_{r-i}-v_{r-i+1})-
			swap(Q)(v_1,\ldots,v_{r-i})\Biggr)$$
	$$+\sum_{1\le i\le r-1} swapM(v_1-v_{r-i+1},v_2-v_{r-i+1},\ldots,v_{r-i}-v_{r-i+1})\times$$
	$$\qquad\qquad\qquad\qquad\qquad \Biggl(swap(Q)(v_1,v_{r-i+2},\ldots,v_r))-swap(Q)(v_{r-i+1},\ldots,v_r)\Biggr).$$
	\vskip .3cm
	\noindent {\bf Step 4. The inverse $push$.}  We have
	$$swap\Bigl(arat(M)\cdot Q\Bigr)(v_2-v_1,\ldots,v_r-v_1,-v_1)=$$
	$$=\sum_{1\le i<j\le r-1} swapM(v_{r-j+2}-v_{r-i+2},\ldots,v_{r-i+1}-v_{r-i+2})\times$$
$$\Biggl(swap(Q)(v_2-v_1,\ldots,v_{r-j+2}-v_1,v_{r-i+3}-v_1,\ldots,v_r-v_1,-v_1)
		\qquad\qquad\qquad\qquad$$
		$$-swap(Q)(v_2-v_1,\ldots,v_{r-j+1}-v_1,v_{r-i+2}-v_1,\ldots,v_r-v_1,-v_1)\Biggr)$$
	$$+\sum_{1\le i\le r-1} swapM(v_{r-i+2}-v_1,\ldots,v_r-v_1,-v_1)\times$$
	$$\Biggl(swap(Q)(v_2-v_{r-i+2},\ldots,v_{r-i+1}-v_{r-i+2})-
			swap(Q)(v_2-v_1,\ldots,v_{r-i+1}-v_1)\Biggr)$$
	$$+\sum_{1\le i\le r-1} swapM(v_2-v_{r-i+2},\ldots,v_{r-i+1}-v_{r-i+2})\times$$
	$$\Biggl(swap(Q)(v_2-v_1,v_{r-i+3}-v_1,\ldots,v_r-v_1,-v_1))-
	swap(Q)(v_{r-i+2}-v_1,\ldots,v_r-v_1,-v_1)\Biggr).$$
	We reindex the double sum by $j\mapsto r-i+1$, $i\mapsto r-j+1$, and we reindex the
	final sum by $j=r-i+1$:
	$$=\sum_{2\le i<j\le r} swapM(v_{i+1}-v_{j+1},\ldots,v_j-v_{j+1})\times$$
$$\Biggl(swap(Q)(v_2-v_1,\ldots,v_{i+1}-v_1,v_{j+2}-v_1,\ldots,v_r-v_1,-v_1)
		\qquad\qquad\qquad\qquad$$
		$$-swap(Q)(v_2-v_1,\ldots,v_i-v_1,v_{j+1}-v_1,\ldots,v_r-v_1,-v_1)\Biggr)$$
	$$+\sum_{i=1}^{r-1} swapM(v_{r-i+2}-v_1,\ldots,v_r-v_1,-v_1)\times$$
	$$\Biggl(swap(Q)(v_2-v_{r-i+2},\ldots,v_{r-i+1}-v_{r-i+2})-
			swap(Q)(v_2-v_1,\ldots,v_{r-i+1}-v_1)\Biggr)$$
	$$+\sum_{j=2}^r swapM(v_2-v_{j+1},\ldots,v_{j}-v_{j+1})\times$$
	$$\Biggl(swap(Q)(v_2-v_1,v_{j+2}-v_1,\ldots,v_r-v_1,-v_1))-
	swap(Q)(v_{j+1}-v_1,\ldots,v_r-v_1,-v_1)\Biggr).\eqno(A.17)$$
	This is thus the expression for $push^{-1}.swap\bigl(arat(M)\cdot Q\bigr)$.  Since the right-hand
	side of the desired equality (A.1) is $-push^{-1}.swap\bigl(arat(M)\cdot Q\bigr)$, we conclude the
	proof by noting that the negative of (A.17) is equal to the left-hand side as given 
	in (A.16).\hfill{$\diamondsuit$}
	\vfill\eject
	\vfill\eject
	It is simpler to start first by computing 
$$\eqalign{\bar{a}(\tau):=\log\,\bar{A}(\tau)
	&=Dgarit\bigl(\bar{E}(\tau)\bigr)\cdot T_{01}\cr
		&=Dgarit\bigl(\bar{E}(\tau)\bigr)\cdot (-a-{{1}\over{2}}U_1+T'_{01})\cr
		&=-\bar{E}(\tau)+1-a+Dgarit\bigl(\bar{E}(\tau)\bigr)\cdot T'_{01}\cr
		&=-\bar{E}(\tau)+1-a+\sum_{n\ge 0} {{1}\over{n!}} Darit\bigl(\bar{\be}(\tau)\bigr)^n\cdot T'_{01}.}$$
		Since $T'_{01}$ starts in depth 2, we have 
		$$\cases{\bar{a}(\tau)(\emptyset)=-a\cr
			\bar{a}(\tau)(u_1)=-\bar{E}(\tau)(u_1)
				=-\sum_{{{n\ge 3}\atop{n\ {\rm odd}}}} \zeta(n)u_1^{n+1}\cr
				\bar{a}(\tau)(u_1,u_2)
				=-\bar{E}(\tau)(u_1,u_2)+T'_{01}(u_1,u_2)
				=-\bar{E}(\tau)(u_1,u_2)+{{1}\over{12}}u_1-{{1}\over{12}}u_2.}$$
				We then see that 
				$$\F\bigl(\bar{a}(\tau)\bigr)(u_1,u_2)=
				-\F\bigl(\bar{E}(\tau)\bigr)(u_1,u_2)+{{1}\over{12u_1u_2}}\F(u_1-u_2)=0.$$
				Thus since
				$$\eqalign{\bar{A}(\tau)&=\exp\bigl(\bar{a}(\tau)\bigr)\cr
					&=1+\bar{a}(\tau)+ {{1}\over{2}}\bar{a}(\tau)^2+\cdots\cr
						&=1+\bigl(-a+\bar{a}(\tau)(u_1)+\bar{a}(\tau)(u_1,u_2)\bigr)+ 
						{{1}\over{2}}\bigl(-a+\bar{a}(\tau)(u_1)+\bar{a}(\tau)(u_1,u_2)\bigr)^2+\cdots\cr
					&=\exp(-a)\ \ \ \ \ {\rm (in\ depth\ 0)}\cr
						&\qquad +\bar{a}(\tau)(u_1)-{{1}\over{2}}a\,\bar{a}(\tau)(u_1)
						-{{1}\over{2}}\bar{a}(\tau)(u_1)a\ \ \ \ \ {\rm (in\ depth\ 1)}\cr
						&\qquad +\bar{a}(\tau)(u_1,u_2)
						-{{1}\over{2}}a\,\bar{a}(\tau)(u_1,u_2)
						-{{1}\over{2}}\bar{a}(\tau)(u_1,u_2)a+\cdots
						+{{1}\over{2}}\bar{a}(\tau)(u_1)\bar{a}(\tau)(u_2),}$$
						using the evenness of $\bar{a}(\tau)(u_1)=-\bar{E}(\tau)(u_1)$, we have
$$\eqalign{\F\bigl(&\bar{A}(\tau)\bigr)(u_1,u_2)
	=\F\bigl(\bar{a}(\tau)\bigr)(u_1,u_2)
		-{{1}\over{2}}a\F\bigl(\bar{a}(\tau)\bigr)(u_1,u_2)
		-{{1}\over{2}}\F\bigl(\bar{a}(\tau)\bigr)(u_1,u_2)a
		+\F\Bigl(\bar{a}(\tau)(u_1)\bar{a}(\tau)(u_2)\Bigr)\cr
		&=\F\Bigl(\bar{a}(\tau)(u_1)\bar{a}(\tau)(u_2)\Bigr)\cr
		&={{1}\over{u_1u_2}}\bar{a}(\tau)(u_1)\bar{a}(\tau)(u_2)
		-{{1}\over{u_1(u_1+u_2)}}\bar{a}(\tau)(u_1)\bar{a}(\tau)(u_1+u_2)
		-{{1}\over{u_2(u_1+u_2)}}\bar{a}(\tau)(u_2)\bar{a}(\tau)(u_1+u_2)
}$$

\vfill\eject
Write $\bar{A}'(\tau)=\bar{A}(\tau)-1$, and
$$\overline{a}(\tau)=\log\bar{A}(\tau)=\bar{A}'(\tau)-{{1}\over{2}}\bar{A}'(\tau)^2+
{{1}\over{3}}\bar{A}'(\tau)^3+\cdots$$
So
$$\F(\bar{a}(\tau))(u_1,u_2)=\F(\bar{A}'(\tau))(u_1,u_2)+\F(\bar{A}'(\tau)^2)(u_1,u_2)$$
$$=\bar{A}'(\tau)(u_1,u_2)+\bar{A}'(\tau)(-u_1,u_1+u_2)+\bar{A}'(\tau)(u_2,-u_1-u_2)
+\bar{A}'(\tau)^2(u_1,u_2)+\bar{A}'(\tau)^2(-u_1,u_1+u_2)+\bar{A}'(\tau)^2(u_2,-u_1-u_2)$$
$$=\bar{A}'(\tau)(u_1,u_2)+\bar{A}'(\tau)(-u_1,u_1+u_2)+\bar{A}'(\tau)(u_2,-u_1-u_2)
+\bar{A}'(\tau)(u_1)\bar{A}'(\tau)(u_2)+\bar{A}'(\tau)(-u_1)\bar{A}'(\tau)(u_1+u_2)
+\bar{A}'(\tau)(u_2)\bar{A}'(\tau)(-u_1-u_2)$$
$$=\bar{a}(\tau)(u_1)\bar{a}(\tau)(u_2)+\bar{a}(\tau)(-u_1)\bar{a}(\tau)(u_1+u_2)
	+\bar{a}(\tau)(u_2)\bar{a}(\tau)(-u_1-u_2).$$
	\vfill\eject
	\vskip .3cm
	\noindent {\bf Proposition 2.5.} {\it Let $P\in GARI$ satisfy the strict 
		group-like elliptic double shuffle relations.  Let $T\in ARI$ be 
			even in depth 1 and such that $dar^{-1}T$ satisfies the strict Fay relations.
			Then $dar^{-1}\, Dgarit(P)\cdot T$ satisfies the strict Fay relations.}
			\vskip .2cm
			\noindent {\bf Proposition 2.6.} {\it Let $P\in GARI$ satisfy the 
				strict group-like elliptic double shuffle relations.  Then the mould
					$dar^{-1}\, Dgarit(P)\cdot ma(t_{01})$ satisfies the strict Fay relations.}
					\vskip .2cm
					\noindent {\bf Proposition 2.7.} {\it Let $P$ be a mould satisfying the group-like
						elliptic double shuffle relations.  Then the mould $dar^{-1}Dgarit(P)\cdot
							ma(t_{01})$ satisfies the Fay relations.}
							\vskip .3cm
							Since $\bar{\ba}(\tau)=Dgarit\bigl(\bar{E}(\tau)\bigr)\cdot ma(t_{01})$
							by (1.8) and $\bar{E}(\tau)$ satisfies the group-like elliptic double shuffle
							relations by Theorem 1.7, Proposition 2.7 shows that 
							$dar^{-1}\bar{\ba}(\tau)$ satisfies the Fay relations, which completes 
							the proof of Theorem 2.2.\hfill{$\diamondsuit$}
							\vskip .3cm
							We now proceed to the proofs of the three propositions.  In the final paragraph
							of this section, we explicitly compute the Fay relations satisfied by 
							$\bar{\ba}(\tau)$ in low depths.
							\vskip .5cm
							\noindent {\bf \S 2.1. Proof of Proposition 2.5.}
							We use three preliminary lemmas, the first of which is a useful reformulation of the
							Fay relation in terms of mould swaps.
							\vskip .3cm
							\noindent {\bf Lemma 2.9.} {\it (i) Let $M$ be a strictly bialternal mould
								and $Q$ a mould that is odd in depth 1 and satisfies the
									strict Fay relations, and set $P:=arat(M)\cdot Q$.  Then $P$ satisfies
									the strict Fay relations.
									\vskip .1cm
									\noindent (ii) Let $M$ be a non-strictly bialternal mould, and let $C=(c_r)_{r\ge 2}$ (with $c_0=c_1=0$)
									be the constant-valued mould such that $swapM+C$ is alternal.  Let $Q$ satisfy the strict Fay relations.
									Then $P:=arat(M)\cdot Q$ satisfies the corrected Fay relations.}
									\vskip .2cm
									Although essentially straightforward and computational, the proof of this key lemma is particularly arduous, 
									and we have banished it to Appendix B below.
									\vskip .2cm
									\noindent {\bf Lemma 2.10.} {\it (i) Let $M\in ARI$ be even in depth 1 and strictly bialternal.  Let $N=\Delta M$ and
										$P:=exp_{Dari}(N)$.  Then $durM$, $dar^{-1}(N)$ and $dar^{-1}(P)$ all satisfy the strict Lie-like Fay relations.
											\vskip .1cm
											\noindent (ii) Let $M\in ARI$ be even in depth 1 and non-strictly bialternal and let $C=(c_r)_{r\ge 0}$ with $c_0=c_1=0$ 
											be the constant-valued mould such that $swapM+C$ is alternal.  Let $N=\Delta M$ and $P:=exp_{Dari}(N)$.  
											Then $durM$ and $dar^{-1}(N)$ satisfy the Lie-like Fay relations for the family $(c_r)_{r\ge 2}$. 
											The mould $dar^{-1}(P)$ satisfies the Fay relations FILL THIS IN.} 
											\vskip .2cm
											\noindent {\bf Proof.}  (i) By Lemma 2.8, $durM$ satisfies the strict Fay relations if and only if 
											$$swap\,durM\bigl(sh((v_1),(v_2,\ldots,v_r))\bigr)+push^{-1}\cdot swap\,durM(v_1,\ldots,v_r)=0.\eqno(2.5)$$
											For the first term $(v_1,\ldots,v_r)$ in the shuffle, we have
											$$swap\,durM(v_1,\ldots,v_r)=swap\bigl((u_1+\ldots+u_r)M(u_1,\ldots,u_r)\bigr)=
											v_1swapM(v_1,\ldots,v_r),\eqno(2.5)$$
											so for the remaining terms $(v_2,\ldots,v_i,v_1,v_{i+1},\ldots,v_r)$ with $i=2,\ldots,r$, we have
											$$swap\,durM(v_2,\ldots,v_i,v_1,v_{i+1},\ldots,v_r)=
											v_2swapM(v_2,\ldots,v_i,v_1,v_{i+1},\ldots,v_r).$$
											Equation (2.5) also allows us to compute the $push$ term:
											$$push^{-1}\cdot swap\,durM(v_1,\ldots,v_r)=
											(v_2-v_1)push^{-1}\cdot swapM(v_1,\ldots,v_r)=(v_2-v_1)\,swapM(v_1,\ldots,v_r),$$
											since $M$ being strictly bialternal, $swapM$ is $push$-invariant.
											Thus the terms in (2.5) add up to  
											$$v_2\,swapM\bigl(sh((v_1),(v_2,\ldots,v_r))\bigr),$$ 
											which
											is zero since $swapM$ is alternal.  Thus $durM$ satisfies the relations (2.5) for $r\ge 2$, so
											by Lemma 2.8, it satisfies the strict Fay relations.

											For $dar^{-1}(N)$, we simply note that $dar^{-1}(N)=durM$ by the definition of $dar$, $dur$ and $\Delta$. Consider
											now the case of $P$.
											We have
											\vskip .3cm
											$$\eqalign{dar^{-1}P&=dar^{-1}exp_{Dari}(N)\cr
												&=1+dar^{-1}\sum_{n\ge 0} {{1}\over{n!}}Darit(N)^n\cdot N\cr
													&=1+dar^{-1}\sum_{n\ge 0} {{1}\over{n!}}\ dar\cdot arat(\Delta^{-1}N)^n\cdot dar^{-1}N\cr
													&=1+\sum_{n\ge 0} {{1}\over{n!}}\ arat(M)^n\cdot durM.}\eqno(2.6)$$
													The mould $durM$ is odd in depth 1, and we just saw that $durM$ satisfies the strict Fay relations.  Thus,
													since $M$ is strictly bialternal, Lemma 2.9 shows that $arat(M)\cdot durM$ satisfies the strict Fay relations, 
													and then successively so does $arat(M)^n\cdot durM$ for all $n\ge 1$.  Thus by (2.6), $dar^{-1}(P)$ satisfies
													the strict Fay relations.
													\vskip .2cm
													\noindent (ii) To show that $durM=dar^{-1}(N)$ satisfies the Fay relations for the family $(c_r)_{r\ge 2}$, we
													will use the equivalent formulation (2.4).  Since bialternal moulds are push-invariant even when the bialternality
													is not strict, the same argument as above shows that the left-hand side of (2.4) adds
													up to $v_2\,swapM\bigl(sh((v_1),(v_2,\ldots,v_r))\bigr)$, but now we have
													$$swapM\bigl(sh((v_1),(v_2,\ldots,v_r))\bigr)=-rc_r,$$
													so the left-hand side of (2.4) is equal to $-rc_rv_2$, proving that $durM$ satisfies the Fay relations in the form
													(2.4).

													The situation is more complicated for $dar^{-1}(P)$.  As above, we have
													$$dar^{-1}P=1+\sum_{n\ge 0}{{1}\over{n!}}arat(M)^n\cdot durM.$$
													\hfill{$\diamondsuit$}
													\vskip .4cm
													Using Lemmas 2.8, 2.9 and 2.10, we can now proceed to the proof of 
													Proposition 2.5.
													\vskip .2cm
													\noindent {\bf Proof of Proposition 2.5.}  Since $P$ satisfies the
													strict group-like elliptic double shuffle relations, $log_{Dari}(P)$
													satisfies the strict Lie-like elliptic double shuffle relations, which means
													by definition that the mould $M:=\Delta^{-1}log_{Dari}(P)$ is strictly
													bialternal.  We have
													$$\eqalign{dar^{-1}\cdot Dgarit(P)\cdot T&=dar^{-1}\cdot \sum_{n\ge 0} {{1}\over{n!}}Darit\bigl(log_{Dari}(P)\bigr)^n\cdot T\cr
														&=dar^{-1}\cdot \sum_{n\ge 0} {{1}\over{n!}}\ dar\cdot arat\bigl(\Delta^{-1}log_{Dari}(P)\bigr)^n\cdot dar^{-1}T\cr
															&=\sum_{n\ge 0}{{1}\over{n!}}\ arat(M)^n\cdot Q,}\eqno(2.5)$$
															where $Q:=dar^{-1}T$ satisfies the strict Fay relations by the hypothesis on 
															$T$.  Since $T$ is assumed even in depth 1, $Q$ is odd in depth 1,
															Therefore we can apply Lemma 2.9 to conclude that $arat(M)\cdot Q$ satisfies 
															the strict Fay relations. Furthermore, $arat(M)\cdot Q$ is zero in depth 1.
															Therefore we can apply Lemma 2.9 again successively to see that in fact
															$arat(M)^n\cdot Q$ satisfies the strict Fay relations for
															all $n\ge 1$.  This proves the result.\hfill{$\diamondsuit$} 
															\vskip .4cm
															\noindent {\bf \S 2.2. Proof of Proposition 2.6.}
															Again we need some preliminary lemmas. 
															\vskip .3cm
															\noindent {\bf Proof of Proposition 2.6.}  Let $U_1=ma([b,a])$ as above, so that $U_1$
															is the mould concentrated in depth $1$ such that $U_1(u_1)=u_1$.  
															Since
															$$t_{01}=-a+{{1}\over{2}}ad(b)(a)+t'_{01},$$
															we have the mould equality
															$$ma(t_{01})=-a+{{1}\over{2}}\,U_1+T'_{01}.$$
															Since $Dgarit(P)$ is an automorphism of $ARI^a_{lu}$, it is additive, 
															so we can write
															$$dar^{-1}Dgarit(P)\cdot ma(t_{01})=
															-dar^{-1}Dgarit(P)\cdot a
															+{{1}\over{2}}\ dar^{-1}Dgarit(P)\cdot U_1+
															dar^{-1}Dgarit(P)\cdot T'_{01}.$$
															By definition of $Dgarit$, we have
															$Dgarit(P)\,U_1=U_1$ and $Dgarit(P)\cdot a=P-1+a$.  Set $B=dar^{-1}U_1$; in fact
															$B$ is the mould concentrated in depth 1 such that $B(u_1)=1$, i.e.~$B=ma(b)$.
															Thus we have 
															$$dar^{-1}Dgarit(P)\cdot ma(t_{01})=
															-dar^{-1}(P-1+a)+{{1}\over{2}}B+dar^{-1}Dgarit(P)\cdot T'_{01}.\eqno(2.15)$$
															The Fay relations only concern terms of depth $r\ge 2$, so to check that 
															$dar^{-1}Dgarit(P)\cdot ma(t_{01})$ satisfies the strict Fay relations, we can
															ignore the terms in depths $0$ and $1$, and restrict our attention to the mould
															$$-dar^{-1}\,P+dar^{-1}Dgarit(P)\cdot T'_{01}.$$
															We will show that each of these two terms separately satisfies the strict
															Fay equations.

															The fact that $-dar^{-1}\,P$ satisfies the strict Fay relations follows from
															Lemma 2.13.  The fact that the term $dar^{-1}Dgarit(P)\cdot T'_{01}$ satisfies the strict
															Fay relations follows from Lemma 2.11, which shows that $dar^{-1}T'_{01}$
															satisfies the strict Fay relations, and Proposition 2.5 with $T=t'_{01}$.
															This concludes the proof of Proposition 2.6.\hfill{$\diamondsuit$}
															\vskip .4cm
															\noindent {\bf \S 2.3.  Proof of Proposition 2.7.}  Recall that if $P\in GARI$
															satisfies the non-strict group-like elliptic double shuffle relations, then
															$log_{Dari}(P)$ is non-strictly bialternal, i.e.~there exists a 
															constant-valued mould $C$ given by $C(\emptyset)=C(v_1)=0$,
															$C(v_1,\ldots,v_r)=c_r$ for $r\ge 2$, such that $swap\bigl(log_{Dari}(P)\bigr)+C$ is 
															alternal.  The proof of Proposition
															2.7 comes down to showing that for such a mould $P$, the correction term $C$
															carries over directly to the mould $dar^{-1}Dgarit(P)\cdot ma(t_{01})$ in the
															sense that this mould satisfies the non-strict Fay relations for the family of 
															constants $(c_r)_{r\ge 2}$.

															We begin by computing
															$$arat(M)\cdot ma(t_{01})$$
															for a bialternal mould $M$ with correction term $C$.
															$$swap(N')\bigl(sh\bigl((v_1), (v_2,\ldots,v_r)\bigr)\bigr)
															+push^{-1}\cdot swap(N')(v_1,\ldots,v_r)=push^{-1}\cdot swap(C_{N'})(v_1,\ldots,v_r).$$
															\vskip .3cm
															\noindent {\bf Corollary 2.13.} The mould $\bar{\ba}_\tau$ satisfies the
															corrected Fay relations (2.13).

															\vfill\eject
															We have
															$$Q(a,b)= -\sum_{i\ge 0} b^i\, a\, \partial^i_b\bigl(p_a(b,a)\bigr).$$
															$$Q(b,a)= -\sum_{i\ge 0} a^i\, b\, \partial^i_a\bigl(p_a(a,b)\bigr).$$
															Let us write $\underline{m}=(m_0,\ldots,m_r)$ and $\overline{m}=
															(m_0,\ldots,m_r)$ with $m_r>0$.
															$$p(a,b)=\sum_{r\ge 1}\sum_{\underline{m}} c_{\underline{m}}a^{m_0}b\cdots b a^{m_r}.$$
															$$p_a(a,b)=\sum_{r\ge 1}\sum_{\overline{m}} c_{\overline{m}}a^{m_0}b\cdots b a^{m_r-1}.$$
															$$p_a(b,a)=\sum_{r\ge 1}\sum_{\overline{m}} c_{\overline{m}}b^{m_0}a\cdots a b^{m_r-1}.$$
															$$\partial^i_a\bigl(p_a(b,a)\bigr)=
															\sum_{r\ge 1}\sum_{\overline{m}}
															\sum_{{{k_0+\cdots+k_r=i}\atop{k_j\le m_j,k_r\le m_r-1}}}
															\bigl({{m_0}\atop{k_0}}\bigr) \cdots
\bigl({{m_{r-1}}\atop{k_{r-1}}}\bigr) \bigl({{m_r-1}\atop{k_r}}\bigr)
	c_{\overline{m}}b^{m_0-k_0}a\cdots ab^{m_{r-1}-k_{r-1}}ab^{m_r-1-k_r}.$$
	$$Q(b,a)=-\sum_{i\ge 0} 
	\sum_{r\ge 1}\sum_{\overline{m}}
	\sum_{{{k_0+\cdots+k_r=i}\atop{k_j\le m_j,k_r\le m_r-1}}}
	\!\!\!\!\!{{(-1)^i}\over{i!}} 
	\bigl({{m_0}\atop{k_0}}\bigr) \cdots
\bigl({{m_{r-1}}\atop{k_{r-1}}}\bigr) \bigl({{m_r-1}\atop{k_r}}\bigr)
	c_{\overline{m}}
	a^i\,b^{m_0-k_0+1}ab^{m_1-k_1}a\cdots ab^{m_{r-1}-k_{r-1}}ab^{m_r-1-k_r}.$$

	\vfill\eject
	Let $EE=Delta(ad_{ari}(invpal)(f5))$, and let $e$ be the key part of this
	polynomial, namely the weight 11, depth 5 part. Write
	$$e_yy=\sum_{{\bf a}} t_{{\bf a}}x^{i_1}yx^{i_2}yx^{i_3}yx^{i_4}x^{i_5}y.$$
	Then
	$$e=\sum_{i\ge 0} {{(-1)^i}\over{i!}}\partial^i_x(e_y)yx^i$$
	and
	$$e_xx=\sum_{i\ge 1} {{(-1)^i}\over{i!}}\partial^i_x(e_y)yx^i.$$
	The partner is given by
	$$f=-\sum_{i\ge 0} {{(-1)^i}\over{i!}} x^iy\partial^i_x(e_x).$$
	Let us now write  $\overline{e}$ for $e(y,x)$ etc.  We have
	$$\overline{e}(x,y)=e(y,x)=e_x(y,x)y+e_y(y,x)x=\overline{e_x}y+\overline{e_y}x,$$
	so $\overline{e}_x=\overline{e_y}$ and $\overline{e}_y=\overline{e_x}$.
	Thus in particular
	$$\overline{\partial^i_x(e_x)}=\partial^i_y\bigl(e_x(y,x)\bigr)=
	\partial^i_y(\overline{e_x})=\partial^i_y(\overline{e}_y).$$
	Thus we have
	$$\overline{f}=f(y,x)=-\sum_{i\ge 0} {{(-1)^i}\over{i!}} y^ix\partial^i_y(\bar{e}_y).$$
	Since 
	$$\overline{e}_y=\overline{e_x}=\sum_{j\ge 1} {{(-1)^j}\over{j!}}\partial^j_y(\overline{e_y})xy^{j-1}$$
	we have
	$$\partial^i_y(\overline{e}_y)=\partial^i_y(\overline{e_x})=
	\sum_{j\ge 1} {{(-1)^j}\over{j!}}\partial^{i+j}_y(\overline{e_y})xy^{j-1}
	+\sum_{j\ge 1} {{(-1)^j}\over{j!}}\partial^j_y(\overline{e_y})x\partial^i_y(y^{j-1})$$
	and since
	$$\partial^i_y(y^{j-1})=\cases{(j-1)y^{j-2}&if $i\le j-1$\cr
		0&if $i>j$,}$$
		we have
		$$\overline{f}=\sum_{i\ge 0} {{(-1)^i}\over{i!}}y^ix
		\Biggl(\sum_{j\ge 1} {{(-1)^j}\over{j!}}\partial^{i+j}_y(\overline{e_y})xy^{j-1}
				+\sum_{j\ge 1} {{(-1)^j}\over{j!}}\partial^j_y(\overline{e_y})x\partial^i_y(y^{j-1})\Biggr)$$
		Let's write the sum up to $d+1$ where $d$ is the key depth (of $e$).
		$$\overline{f}=\sum_{i\ge 0} {{(-1)^i}\over{i!}}y^ix
		\Biggl(\sum_{j=1}^d {{(-1)^j}\over{j!}}\partial^{i+j}_y(\overline{e_y})xy^{j-1}
				+\sum_{j=1}^d {{(-1)^j}\over{j!}}\partial^j_y(\overline{e_y})x\partial^i_y(y^{j-1})\Biggr)$$
		$$=\sum_{j=1}^d \sum_{i=0}^{d-j}{{(-1)^{i+j}}\over{i!j!}}y^ix\partial^{i+j}_y(\overline{e_y})xy^{j-1}
		+\sum_{j=1}^d \sum_{i=0}^{j-1}{{(-1)^{i+j}}\over{i!j!}}y^ix\partial^j_y(\overline{e_y})x\partial^i_y(y^{j-1})$$
		$$=\sum_{j=1}^d \sum_{i=0}^{d-j}{{(-1)^{i+j}}\over{i!j!}}y^ix\partial^{i+j}_y(\overline{e_y})xy^{j-1}
		+\sum_{j=1}^d \sum_{i=0}^{j-1}{{(-1)^{i+j}(j-1)\cdots(j-i)}\over{i!j!}}y^ix\partial^j_y(\overline{e_y})xy^{j-i-1}$$
		Thus we have
		$$\overline{f}_yy=
		\sum_{j=2}^d \sum_{i=0}^{d-j}{{(-1)^{i+j}}\over{i!j!}}y^ix\partial^{i+j}_y(\overline{e_y})xy^{j-1}
		+\sum_{j=1}^d \sum_{i=0}^{j-2}{{(-1)^{i+j}(j-1)\cdots(j-i)}\over{i!j!}}y^ix\partial^j_y(\overline{e_y})xy^{j-i-1}$$
		\vfill\eject
		$$e=-{{1}\over{4}}\Bigl(yyxxxyx+xyyyxxx+yyxyxxx+xyxxxyy+yxxxyxy+yxyxxxy+xxxyxyy+xxxyyyx\Bigr)$$
		$$+{{1}\over{2}}\Bigl(yyxxyxx+yxyyxxx+xxxyyxy+xxyxxyy+yxxyxxy+xyyxxxy+xxyyyxx+yxxxyyx\Bigr)$$
		$$-{{13}\over{12}}\Bigl(xxyxyxy+xyxyxxy+yxyxyxx+xyyxxyx+xxyyxyx+yxxyxyx+xyxyyxx+xyxxyyx\Bigr)$$
		$$+{{5}\over{6}}\Bigl(xyxxyxy+xyyxyxx+yxyxxyx+xxyxyyx\Bigr)$$
		$$+{{1}\over{3}}\Bigl(xxyyxxy+yxxyyxx\Bigr)$$
		$$+{{8}\over{3}}xyxyxyx.$$
		Then
		$$e_yy=-{{1}\over{4}}\Bigl(xyxxxyy+yxxxyxy+yxyxxxy+xxxyxyy\Bigr)$$
		$$+{{1}\over{2}}\Bigl(xxxyyxy+xxyxxyy+yxxyxxy+xyyxxxy\Bigr)$$
		$$-{{13}\over{12}}\Bigl(xxyxyxy+xyxyxxy\Bigr)$$
		$$+{{5}\over{6}}\Bigl(xyxxyxy\Bigr)$$
		$$+{{1}\over{3}}\Bigl(xxyyxxy\Bigr)$$
$$e_yy=-{{1}\over{4}}\Bigl(v_1v_2^3+v_2^3v_3+v_2v_3^3+v_1^3v_2\Bigr)
	+{{1}\over{2}}\Bigl(v_1^3v_3+v_1^2v_2^2+v_2^2v_3^2+v_1v_3^3\Bigr)$$
$$-{{13}\over{12}}\Bigl(v_1^2v_2v_3+v_1v_2v_3^2\Bigr)
	+{{5}\over{6}}v_1v_2^2v_3
	+{{1}\over{3}}v_1^2v_3^2.$$
	Write 
	$$e=\sum_{a+b+c+d=4}t_{abcd} x^ayx^byx^cyx^d=\sum_{i=0}^2 {{(-1)^i}\over{i!}}
	\partial^i_x(e_yy)x^i,$$
	$$e_yy=\sum_{a+b+c=4}t_{abc0} x^ayx^byx^cy
	=\sum_{a+b+c=4}t_{abc0} v_1^av_2^bv_3^c,$$
	$$e_xx=\sum_{d>0}t_{abcd} x^ayx^byx^cyx^d=
	\sum_{i=1}^3 {{(-1)^i}\over{i!}} \partial^i_x(e_yy)x^i
	=-\partial_x(e_yy)x+{{1}\over{2}}\partial^2_x(e_yy)x^2
	-{{1}\over{6}}\partial^3_x(e_yy)x^3$$
	so
	$$e_x=\sum_{d>0}t_{abcd} x^ayx^byx^cyx^{d-1}=
	\sum_{i=1}^2 {{(-1)^i}\over{i!}} \partial^i_x(e_yy)x^{i-1}
	=-\partial_x(e_yy)+{{1}\over{2}}\partial^2_x(e_yy)x-{{1}\over{6}}\partial^3_x(eyy)x^2.$$
	For the partner, we have
$$f=-\sum_{i=0}^2 {{(-1)^i}\over{i!}}x^iy\partial^i_x(e_x)
	=-ye_x+xy\partial_x(e_x)-{{1}\over{2}}x^2y\partial^2_x(e_x)$$
$$=y\partial_x(e_yy)
	-{{1}\over{2}}y\partial^2_x(e_yy)x
-{{1}\over{2}}xy\partial^2_x(e_yy)
	+{{1}\over{6}}xy\partial^3_x(e_yy)x
	+{{1}\over{6}}y\partial^3_x(e_yy)x^2
	+{{1}\over{6}}x^2y\partial^3_x(e_yy).$$
	Now we can start switching $x$ and $y$. Letting $\bar{e_y}=e_y(y,x)$, we have
$$\overline{f}=f(y,x)=x\partial_y(\overline{e_y}x)
	-{{1}\over{2}}x\partial^2_y(\overline{e_y}x)y
-{{1}\over{2}}yx\partial^2_y(\overline{e_y}x)
	+{{1}\over{6}}yx\partial^3_y(\overline{e_y}x)y
	+{{1}\over{6}}x\partial^3_y(\overline{e_y}x)y^2
	+{{1}\over{6}}y^2x\partial^3_y(\overline{e_y}x).$$
	Thus
	$$\overline{f}_yy=
	-{{1}\over{2}}x\partial^2_y(\overline{e_y}x)y
	+{{1}\over{6}}yx\partial^3_y(\overline{e_y}x)y
	+{{1}\over{6}}x\partial^3_y(\overline{e_y}x)y^2$$
	which we can even write as
	$$\overline{f}_yy=
	-{{1}\over{2}}x\partial^2_y(\overline{e_yy})y
	+{{1}\over{6}}yx\partial^3_y(\overline{e_yy})y
	+{{1}\over{6}}x\partial^3_y(\overline{e_yy})y^2$$
	if we really want to relate it to $e_yy$, or else as
	$$\overline{f}_yy=
	-{{1}\over{2}}x\partial^2_y(\overline{e_y})xy
	+{{1}\over{6}}yx\partial^3_y(\overline{e_y})xy
	+{{1}\over{6}}x\partial^3_y(\overline{e_y})xy^2$$
	if we want to separate the terms according to the power of $y$ on the right.

	Using
	$$e_yy=\sum_{a+b+c=4}t_{abc0} x^ayx^byx^cy$$ 
	we compute
	$$e_y=\sum_{a+b+c=4}t_{abc0} x^ayx^byx^c$$ 
	$$\overline{e_y}=\sum_{a+b+c=4}t_{abc0} y^axy^bxy^c$$ 
	$$\partial_y(\overline{e_y})=\sum_{a+b+c=4}t_{abc0} 
	\Bigl(ay^{a-1}xy^bxy^c+by^axy^{b-1}xy^c+cy^axy^bxy^{c-1}\Bigr)$$
	$$\partial^2_y(\overline{e_y})=\sum_{a+b+c=4}t_{abc0} 
	\Bigl(a(a-1)y^{a-2}xy^bxy^c+aby^{a-1}xy^{b-1}xy^c+acy^{a-1}xy^bxy^{c-1}
			+aby^{a-1}xy^{b-1}xy^c$$
			$$+b(b-1)y^axy^{b-2}xy^c+bcy^axy^{b-1}xy^{c-1}+acy^{a-1}xy^bxy^{c-1}+bcy^axy^{b-1}xy^{c-1}+c(c-1)y^axy^bxy^{c-2}\Bigr)$$
	$$=\sum_{a+b+c=4}t_{abc0} 
	\Bigl(a(a-1)y^{a-2}xy^bxy^c
			+b(b-1)y^axy^{b-2}xy^c
			+c(c-1)y^axy^bxy^{c-2}$$
			$$+2aby^{a-1}xy^{b-1}xy^c
			+2acy^{a-1}xy^bxy^{c-1}
			+2bcy^axy^{b-1}xy^{c-1}\Bigr)$$
	$$\partial^3_x(\overline{e_y})=
	\sum_{a+b+c=4}t_{abc0} 
	\Bigl(a(a-1)(a-2)y^{a-3}xy^bxy^c
			+b(b-1)(b-2)y^axy^{b-3}xy^c
			+c(c-1)(c-2)y^axy^bxy^{c-3}$$
			$$+3a(a-1)by^{a-2}xy^{b-1}xy^c
			+3a(a-1)cy^{a-2}xy^bxy^{c-1}
			+3b(b-1)cy^axy^{b-2}xy^{c-1}
			+3bc(c-1)y^axy^{b-1}xy^{c-2}$$
			$$+3ab(b-1)y^{a-1}xy^{b-2}xy^c
			+3ac(c-1)y^{a-1}xy^bxy^{c-2}
			+6abcy^{a-1}xy^{b-1}xy^{c-1}\Bigr)$$
	Now we can compute the three terms of $\overline{f}_yy$:
	$$-{{1}\over{2}}x\partial^2_y(\overline{e_y})xy
	=\sum_{a+b+c=4}t_{abc0} 
	\Bigl(a(a-1)xy^{a-2}xy^bxy^cxy
			+b(b-1)xy^axy^{b-2}xy^cxy
			+c(c-1)xy^axy^bxy^{c-2}xy$$
			$$+2abxy^{a-1}xy^{b-1}xy^cxy
			+2acxy^{a-1}xy^bxy^{c-1}xy
			+2bcxy^axy^{b-1}xy^{c-1}xy\Bigr),$$
	so the coefficient of a single word of the form $xy^axy^bxy^cxy$ with $a+b+c=2$
	(coming from the $\partial^2$ term) is given by
$$t_{a+2,bc0} (a+2)(a+1)
+t_{a,b+2,c0}(b+2)(b+1)
	+t_{ab,c+2,0}(c+2)(c+1)$$
$$+2t_{a+1,b+1,c0}(a+1)(b+1)
+2t_{a+1,b,c+1,0}(a+1)(c+1)
	+2t_{a,b+1,c+1,0}(b+1)(c+1).$$
	So for the different $(a,b,c)$ with $a+b+c=2$ we have
	$$\cases{
		t_{1120}=6t_{3100}+6t_{1300}+2t_{1120}+8t_{2200}+4t_{2110}+4t_{1210}&$(1,1,0)$\cr
			t_{1210}=6t_{3010}+2t_{1210}+6t_{1030}+6t_{2110}+8t_{2020}+4t_{1120}&$(1,0,1)$\cr
			t_{2110}=2t_{2110}+6t_{0310}+6t_{0130}+4t_{1210}+4t_{1120}+8t_{0220}&$(0,1,1)$\cr
			t_{1030}=12t_{4000}+2t_{2200}+2t_{2020}+6t_{3100}+6t_{3010}+2t_{2110}&$(2,0,0)$\cr
			t_{2020}=2t_{2200}+12t_{0400}+2t_{0220}+6t_{1300}+2t_{1210}+6t_{0310}&$(0,2,0)$\cr 
			t_{3010}=2t_{2020}+2t_{0220}+12t_{0040}+2t_{1120}+6t_{1030}+6t_{0130}&$(0,0,2)$.}$$
			\vfill\eject
			$$s(P)=\sum_{i\ge 0} {{(-1)^i}\over{i!}}\partial^i_x(P)yx^i.$$
			$$t(Q)=\sum_{i\ge 0} {{(-1)^i}\over{i!}}\partial^i_y(Q)xy^i.$$
			$$s'(P)=\sum_{i\ge 0} {{(-1)^i}\over{i!}}x^iy\partial^i_x(P).$$

			Let $f$ be a Lie polynomial with no words starting or ending in $x$.  Write
			$f_yy=xPy$.  Then $f=[x,s(P)]$.
			\vskip .1cm
			Let $f$ be a Lie polynomial with no words starting or ending in $x$.  Write
			$f_xx=yQx$.  Then $f=[y,t(Q)]$.
			\vskip .1cm

			Let $f$ be a Lie polynomial with no words starting and ending in either $x$ or $y$.
			Then $f=[x,s(P)]=[y,t(Q)]$ so $[t(Q),y]+[x,s(P)]=0$ and $E=t(Q)$, $F=s(P)$ are partners.
			What is the relationship between $s(P)(y,x)$ and $t(Q)$?
			We have
			$$s(P)=\sum_{i\ge 0} {{(-1)^i}\over{i!}}\partial^i_x(P)yx^i$$
			so
			$$s(P)(y,x)=\sum_{i\ge 0} {{(-1)^i}\over{i!}}\partial^i_y(P(y,x))xy^i$$
			while 
			$$t(Q)=\sum_{i\ge 0} {{(-1)^i}\over{i!}}\partial^i_y(Q)xy^i.$$
			This means that
			$$s(P)(y,x) = t(\overline{P})$$
			where $\overline{P}=P(y,x)$.
			So the relationship between the two comes directly from the relationship between
			$\overline{P}$ and $Q$. We have $f_x=yQ$ and $f_y=xP$, so $f_y(y,x)=yP(y,x)$.  
			Writing
			$\overline{f}=f(y,x)$, we have 
			$$\overline{f}=\overline{f}_xx+\overline{f}_yy=f(y,x)=f_y(y,x)x+f_x(y,x)y$$
			so $f_y(y,x)=\overline{f}_x=y\overline{P}$, i.e.~the relationship is 
			$$f_x=yQ,\ \ \overline{f}_x=y\overline{P}.$$

			Suppose that $t(Q)$ is $\Delta$-bialternal.  We need to show that
			$s(P)(y,x)$ is as well, but in view of what we just showed, it suffices
			to show that $t(\overline{P})$ is $\Delta$-bialternal:
			$$t(\overline{P})=\sum_{i\ge 0}{{(-1)^i}\over{i!}}\partial^i_y(\overline{P})xy^i.$$

			Another way to think about this is that there is a natural involution $dig$
			on $ARI^{pol,push}$ sending $p(a,b)\mapsto -q(b,a)$, or if one prefers, 
			sending the derivation $a\mapsto p(a,b)$, $b\mapsto q(a,b)$ to the derivation 
			$a\mapsto -q(b,a)$, $b\mapsto -p(b,a)$, and
			one wants to know if the involution $dig$ actually preserves the subspace of 
			$\Delta(ARI^{al/al})$ of $ARI^{push}$. 

			Equivalently, we need to know 
			whether the involution $\Delta^{-1}\circ dig\circ\Delta$ of
			$ARI^{\Delta,push}$ preserves the subspace $ARI^{al/al}$.

			\vfill\eject
			\noindent {\bf Lemma.} {\it Suppose that $swap(N)$ is push/push.
				Then $Fay(N')=0$ is equivalent to 
					$$swap(N')\bigl(sh\bigl((v_1), (v_2,\ldots,v_r)\bigr)\bigr)
					=-push^{-1}\cdot swap(N')(v_1,\ldots,v_r)\eqno(1)$$
					which is in turn 
					equivalent to the first alternality relation of $swap(N/\Delta)$.}
					\vskip .2cm
					\noindent Proof. Let $N'=dar^{-1}(N)$ so
					$$N'(u_1,\ldots,u_r)={{N(u_1,\ldots,u_r}\over{u_1\cdots u_r}}$$
							and
							$$swap(N')(v_1,\ldots,v_r)={{swap(N)(v_1,\ldots,v_r}\over{(v_1-v_2)\cdots
								(v_{r-1}-v_r)v_r}}.$$
								Let us first show that $Fay(N')=0$ is equivalent to
								$$swap(N')\bigl(sh\bigl((v_1), (v_2,\ldots,v_r)\bigr)\bigr)
								=-push^{-1}\cdot swap(N')(v_1,\ldots,v_r)\eqno(1)$$
								where 
								$$swap(N')(v_1,\ldots,v_r)={{swap(N)(v_1,\ldots,v_r)}\over
								{(v_1-v_2)\cdots(v_{r-1}-v_r)v_r}}.$$
								Then
								$$push^{-1}swap(N')(v_1,\ldots,v_r)=
								-{{swap(N)(v_2-v_1,\ldots,v_r-v_1,-v_1)}\over {v_1(v_2-v_3)\cdots(v_{r-1}-v_r)v_r}}.$$
								The fact that $swap(N/\Delta)$ is push-invariant means that $swap(N)$ is
								push-invariant, so
								$$-push^{-1}swap(N')(v_1,\ldots,v_r)=
								{{swap(N)(v_1,\ldots,v_r)}\over {v_1(v_2-v_3)\cdots(v_{r-1}-v_r)v_r}}.$$
								Replacing this in the RHS of (1) gives the equivalent form 
								$$swap(N')\bigl(sh\bigl((v_1), (v_2,\ldots,v_r)\bigr)\bigr)
								={{swap(N)(v_1,\ldots,v_r)}\over {v_1(v_2-v_3)\cdots(v_{r-1}-v_r)v_r}}\eqno(2)$$
								i.e.
								$$swap(N')\bigl(sh\bigl((v_1), (v_2,\ldots,v_r)\bigr)\bigr)
								={{(v_1-v_2)}\over{v_1}}swap(N')(v_1,\ldots,v_r)\eqno(3)$$
							i.e.
							$${{v_2}\over{v_1}}swap(N')(v_1,v_2,\ldots,v_r)+swap(N')(v_2,v_1,\ldots,v_r)
							+\cdots+swap(N')(v_2,\ldots,v_r,v_1)=0.\eqno(4)$$
							Writing
							$$swap(N')(v_1,\ldots,v_r)=dur\cdot swap(N/\Delta)(v_1,\ldots,v_r)=
							v_1\,swap(N/\Delta)(v_1,\ldots,v_r),$$
							we rewrite (4) in terms of $swap(N/\Delta)$ as
							$$v_2\,swap(N')(v_1,v_2,\ldots,v_r)+v_2\,swap(N/\Delta)(v_2,v_1,\ldots,v_r)
							+\cdots+v_2\,swap(N/\Delta)(v_2,\ldots,v_r,v_1)=0.\eqno(5)$$
							which is equivalent to the first alternality relation on $swap(N/\Delta)$.
							\hfill{$\diamondsuit$}
							\vfill\eject
							$$Fay(N')(u_1,u_2,u_3,u_4)=N'(u_1,u_2,u_3,u_4)
							+N'(-u_1,u_1+u_2,u_3,u_4)
							+N'(u_2,-u_1-u_2,u_1+u_2+u_3,u_4)$$
							$$+N'(u_2,u_3,-u_1-u_2-u_3,u_1+u_2+u_3+u_4)
							+N'(u_2,u_3,u_4,-u_1-u_2-u_3-u_4)=0.$$
							Apply swap.
							$$swap\Bigl(Fay(N')\Bigr)(v_1,v_2,v_3,v_4)=N'(v_4,v_3-v_4,v_2-v_3,v_1-v_2)
							+N'(-v_4,v_3,v_2-v_3,v_1-v_2)
							+N'(v_3-v_4,-v_3,v_2,v_1-v_2)$$
							$$+N'(v_3-v_4,v_2-v_3,-v_2,v_1)
							+N'(v_3-v_4,v_2-v_3,v_1-v_2,-v_1)=0.$$
							Write it in terms of $swap(N')$.
							$$swap\Bigl(Fay(N')\Bigr)(v_1,v_2,v_3,v_4)=swap(N')(v_1,v_2,v_3,v_4)
							+swap(N')(v_1-v_4,v_2-v_4,v_3-v_4,-v_4)
							+swap(N')(v_1-v_4,v_2-v_4,-v_4,v_3-v_4)$$
							$$+swap(N')(v_1-v_4,-v_4,v_2-v_4,v_3-v_4)
							+swap(N')(-v_4,v_1-v_4,v_2-v_4,v_3-v_4)=0.$$
							Apply $push^{-1}$.
							$$push^{-1}\cdot swap\Bigl(Fay(N')\Bigr)(v_1,v_2,v_3,v_4)=
							swap(N')(v_2-v_1,v_3-v_1,v_4-v_1,-v_1)
							+swap(N')(v_2,v_3,v_4,v_1)
							+swap(N')(v_2,v_3,v_1,v_4)$$
							$$+swap(N')(v_2,v_1,v_3,v_4)
							+swap(N')(v_1,v_2,v_3,v_4)=0.$$
							Now use that $swap(N')$ is push-invariant for only the first term of the RHS.
							$$push^{-1}\cdot swap\Bigl(Fay(N')\Bigr)(v_1,v_2,v_3,v_4)=
							swap(N')(v_1,v_2,v_3,v_4)
							+swap(N')(v_2,v_3,v_4,v_1)
							+swap(N')(v_2,v_3,v_1,v_4)$$
							$$+swap(N')(v_2,v_1,v_3,v_4)
							+swap(N')(v_1,v_2,v_3,v_4)=0,$$
							i.e.
							$$push^{-1}\cdot swap\Bigl(Fay(N')\Bigr)(v_1,v_2,v_3,v_4)=
							swap(N')\Bigl((v_1),(v_2,v_3,v_4)\Bigr).$$
							\vfill\eject
							If $A$ is circ-neutral and push-invariant then the first alternal relation holds.  This is true for moulds in the $u_i$, but is it true for moulds in the 
							$v_i$?  What is the relation between the first alternality relation in the
							$v_i$ and $circ$-neutrality?

							$$A(u_1,u_2,u_3)+A(u_2,u_3,u_1)+A(u_3,u_1,u_2)=0$$

							$$pushA(u_1,u_2,u_3)=A(-u_1-u_2-u_3,u_1,u_2)$$
							So by push-invariance of $A$, we have
							$$A(u_1,u_2,u_3)=A(-u_1-u_2-u_3,u_1,u_2)=A(u_3,-u_1-u_2-u_3,u_1)=
							A(u_2,u_3,-u_1-u_2-u_3),$$
							i.e.
							$$A(u_1,u_2,u_3)=A(u_0,u_1,u_2)=A(u_3,u_0,u_1)=A(u_2,u_3,u_0),$$

$$circpushA(u_1,u_2,u_3)=pushA(u_2,u_3,u_1)=A(-u_1-u_2-u_3,u_2,u_3)
	=A(u_0,u_2,u_3)$$
	$$circ2push^2A(u_1,u_2,u_3)=push^2A(u_3,u_1,u_2)=A(u_2,u_0,u_3)$$
	Therefore
	$$A(u_1,u_2,u_3)+circA(u_1,u_2,u_3)+circ2A(u_1,u_2,u_3)=$$
	$$push^3A(u_1,u_2,u_3)+circpushA(u_1,u_2,u_3)+circ2push2A(u_1,u_2,u_3)=$$
	$$=A(u_2,u_3,u_0)+A(u_0,u_2,u_3)+A(u_2,u_0,u_3)=0$$
	which is 1st alternality.

	$$A(v_1,v_2,v_3)+A(v_2,v_3,v_1)+A(v_3,v_1,v_2)=0$$
	$$A(v_1,v_2,v_3)+circA(v_1,v_2,v_3)+circ2A(v_1,v_2,v_3)=0$$
	$$push(A)(v_1,v_2,v_3)=A(-v_3,v_1-v_3,v_2-v_3)$$
	So push-invariance means
$$A(v_1,v_2,v_3)=A(-v_3,v_1-v_3,v_2-v_3)
	=A(v_3-v_2,-v_2,v_1-v_2)=A(v_2-v_1,v_3-v_1,-v_1)$$
	$$circpushA(v_1,v_2,v_3)=pushA(v_2,v_3,v_1)=
	A(-v_1,v_2-v_1,v_3-v_1)$$
	$$circ2push2A(v_1,v_2,v_3)=push2A(v_3,v_1,v_2)=
	A(v_2-v_1,-v_1,v_3-v_1)$$
	so
	$$push^3A(v_1,v_2,v_3)+circpushA(v_1,v_2,v_3)+circ2push2A(v_1,v_2,v_3)=$$
	$$A(v_2-v_1,v_3-v_1,-v_1)+
A(-v_1,v_2-v_1,v_3-v_1)
	A(v_2-v_1,-v_1,v_3-v_1)=0,$$
	which is just 1st alternality.
	\vfill\eject
	By definition, a mould $N$ satisfyies the Fay relations if $\F(N')(u_1,\ldots,
			u_r)=0$.

	$$A(v_1,v_2,v_3,v_4)=A(-v_4,v_1-v_4,v_2-v_4,v_3-v_4)$$
	$$=A(v_4-v_3,-v_3,v_1-v_3,v_2-v_3)$$
	$$=A(v_3-v_2,v_4-v_2,-v_2,v_1-v_2)$$
	$$=A(v_2-v_1,v_3-v_1,v_4-v_1,-v_1)$$
	\vfill\eject
	\noindent {\bf Appendix A: Mould basics}
	\vskip .2cm
	In this appendix, we introduce the basic mould definitions needed for this
	article.  For the purposes of this article, {\it moulds} will refer only to
	rational function moulds over $\Q$.
	\vskip .2cm
	\noindent {\bf \S A.1. Definition of moulds and the swap operator.}
A {\it mould} is a tuple 
$\bigl(M_r(u_1,\ldots,u_r)\bigr)_{r\ge 0}$ in which $M_0\in \Q$ and 
for each $r\ge 1$, $M_r$ is a rational function of
$r$ commutative variables.
\vskip .3cm
We will introduce several properties and operators on moulds.  The key operator
is the $swap$, which acts on a mould $M=(M_r)$ as follows:
$$swapM_r(v_1,\ldots,v_r)=M_r(v_r,v_{r-1}-v_r,\ldots,v_1-v_2).$$
We usually write the commutative variables of the swap of a mould as
$v_i$ rather than $u_i$, simply because the different variable names make 
it easy for us to recognize right away that we are dealing with the swap
of a mould.  We also generally write $\overline{M}$ for $swapM$, and
when the number of variables is specified, we drop the subscript $r$ (for
instance, it is understood that $M(u_1,u_2)=M_2(u_1,u_2)$).

Let $ARI$ denote the vector space of moulds with constant term $0$, and
let $\overline{ARI}=swap(ARI)$.  In fact $\overline{ARI}$ is isomorphic to $ARI$
since the swap-operator is invertible; it is nothing other than the vector 
space of (rational-function) moulds of the commutative variables $v_i$ with 
constant term $0$.  However, we will see below that it is useful to 
distinguish between the two spaces.
\vskip .3cm
\noindent {\bf \S A.2. Standard multiplication and Lie bracket on $ARI$.}
There are three different Lie brackets that one can put on the
space $ARI$.  We begin by introducing the standard
mould multiplication that \'Ecalle denotes $mu(A,B)$:
$$mu(A,B)(u_1,\ldots,u_r)=\sum_{i=0}^r A(u_1,\ldots,u_i)B(u_{i+1},\ldots,
u_r).$$
The associated Lie bracket $lu$ is defined by $lu(A,B)=mu(A,B)-mu(B,A)$.
We write $ARI_{lu}$ for $ARI$ viewed as a Lie algebra for the $lu$-bracket.
The identical formulas yield a multiplication and Lie algebra (also called
$mu$ and $lu$) on $\overline{ARI}$.

Let $c_i=ad(a)^{i-1}(b)$, $i\ge 1$.  Power series in the $c_i$ form a subspace 
of $\Qab$, namely the kernel of the derivation $\partial_x$ defined by 
$x\mapsto 1$, $y\mapsto 0$.  The standard map from power series 
in the $c_i$ to polynomial-valued moulds is given by linearly
extending the map on monomials
$$ma:c_{a_1}\cdots c_{a_r}\mapsto u_1^{a_1-1}\cdots u_r^{a_r-1}.$$

If $f$ and $g$ are power series in $\QC$ and $A=ma(f)$, $B=ma(g)$, then
the multiplication $mu$ is nothing other than the mould translation 
of the usual non-commutative multiplication,
and $lu$ the usual Lie bracket:
$$mu(A,B)=ma(fg),\ \ \ \ lu(A,B)=ma\bigl([f,g]\bigr).$$
\vskip .1cm
\noindent {\bf \S A.3. A few more mould operators.}  We define the operators
$$\cases{durP(u_1,\ldots,u_r)=(u_1+\cdots+u_r)\ P(u_1,\ldots,u_r)\cr
dar(P)(u_1,\ldots,u_r)=(u_1\cdots u_r)\,P(u_1,\ldots,u_r)\cr
\Delta P(u_1,\ldots,u_r)=(u_1\cdots u_r)u_1\cdots u_r\,P(u_1,\ldots,u_r).}$$
We have the following useful basic result from [S2].  Let $ARI^a$ denote the 
vector space generated by $ARI$ together with the mould $ma(a)$, which is 
concentrated in depth $0$ with value $a$.  By a slight abuse of notation, we
continue to refer to this mould as $a$, the meaning of the symbol (as variable
or as mould) being clear from the context. We make $ARI^a$ into a
Lie algebra $ARI^a_{lu}$, inside of which $ARI_{lu}$ forms a Lie ideal,
by extending the Lie bracket $lu$ to all of $ARI^a$ via the relation
$$lu(P,a)=durP\ \ {\rm for}\ \ P\in ARI.$$
This identity is nothing other than an extension to all moulds $P\in ARI$ 
the identity ([ref??]) for the case $P=ma(p)$ where $p\in \Qab$:
$$ma([p,a])=dur\bigl(ma(p)\bigr).$$
\vskip .3cm
\noindent {\bf \S A.4. The $ari$-bracket on $ARI$.}
In order to define \'Ecalle's $ari$-bracket, we first introduce three
derivations of $ARI_{lu}$ associated to a given mould $A\in ARI$. It
is non-trivial to prove that these operators are actually derivations
(cf.~[S2, Prop. 2.2.1]).
\vskip .3cm
\noindent {\bf Definition.} ([Ec]) Let $B\in ARI$.  Then the derivation $amit(B)$ of
$ARI_{lu}$ is given by
$$\bigl(amit(B)\cdot A\bigr)(u_1,\ldots,u_r)=
\sum_{1\le i<j<r} A(u_1,\ldots,u_i, u_{i+1}+\cdots+u_{j+1},u_{j+2},\ldots,u_r)
B(u_{i+1},\ldots,u_j),$$
and the derivation $anit(B)$ is given by
$$\bigl(anit(B)\cdot A\bigr)(u_1,\ldots,u_r)=
\sum_{1<i<j\le r} A(u_1,\ldots,u_i+\cdots+u_j,u_{j+1},\ldots,u_r)
B(u_{i+1},\ldots,u_j).$$

\noindent We also have corresponding derivations $\overline{amit}(B)$ and
$\overline{anit}(B)$ of $\overline{ARI}_{lu}$ 
for $B\in\overline{ARI}$, given by the formulas 
$$\bigl(\overline{amit}(B)\cdot A\bigr)(v_1,\ldots,v_r)=
\sum_{1\le i<j<r} A(u_1,\ldots,u_i,u_{j+1},\ldots,u_r)B(u_{i+1}-u_{j+1},
\ldots,u_j-u_{j+1}),$$
$$\bigl(\overline{anit}(B)\cdot A\bigr)(v_1,\ldots,v_r)=
\sum_{1<i<j\le r} A(u_1,\ldots,u_i,u_{j+1},\ldots,u_r)B(u_{i+1}-u_i,
\ldots,u_j-u_i).$$
Finally, \'Ecalle defines the derivation $arit(B)$ on $ARI_{lu}$ by
$$arit(B)=amit(B)-anit(B),$$
and the $ari$-bracket on $ARI$ by
$$ari(A,B)=arit(B)\cdot A-arit(A)\cdot B+lu(A,B),$$
as well as the derivation $\overline{arit}$ on $\overline{ARI}_{lu}$
and the bracket $\overline{ari}$ on $\overline{ARI}$ by the same
formulas with overlines.
\vskip .2cm
The definitions of $amit$, $anit$, $arit$ and $ari$ are generalizations
to all moulds of familiar derivations of $\QC$.  Indeed,
if $b,b'\in \QC$ and $A=ma(b)$, $B=ma(b')$, then
$$amit(B)\cdot A=ma\bigl(D^l_g(f)\bigr)$$
where $D^l_g$ is defined by $x\mapsto 0$, $y\mapsto b'y$,
$$anit(B)\cdot A=ma\bigl(D^r_{b'}(b)\bigr)$$
where $D^r_{b'}$ is defined by $x\mapsto 0$, $y\mapsto yb'$, and thus
$$arit(B)\cdot A=ma\bigl(-d_{b'}(b)\bigr)$$
where $d_g$ is the Ihara derivation $x\mapsto 0$, $y\mapsto [y,b']$
(see [Iharader]), and
$$ ari(A,B)=ma\bigl([b,b']+d_b(b')-d_{b'}(b)\bigr)=ma\bigl(\{b,b'\}\bigr).
\eqno(A.1)$$
corresponds to the Ihara or Poisson Lie bracket 
on the free Lie algebra on two ${\rm Lie}[a,b]$ (see [S2, Corollary 3.3.4]). 
\vskip .3cm
\noindent {\bf \S A.5. The $Dari$-bracket.} 
We now pass to the $Dari$-bracket, which is the Lie bracket on 
$ARI$ obtained by transfer by the $\Delta$-operator: it is given by
$$Dari(A,B) = \Delta\Bigl(ari\bigl(\Delta^{-1}(A),\Delta^{-1}(B)\bigr)
\Bigr).
$$
This means that $\Delta$ gives an isomorphism of Lie algebras
$$\Delta:ARI_{ari}\buildrel\sim\over\rightarrow ARI_{Dari}.$$

It is shown in [S3, Prop. 3.2.1] that we have a second definition
for the $Dari$-bracket, which is more complicated but sometimes very
useful in certain proofs.  Let $dar$ denote the mould operator
defined by $dar(A)(u_1,\ldots,u_r)=u_1\cdots u_r\,A(u_1,\ldots,u_r)$.
We begin by introducing,
for each $A\in ARI$, an associated derivation $Darit(A)$ of $ARI_{lu}$
by the following formula:
$$Darit(A)=dar\cdot\Bigl(-arit\bigl(\Delta^{-1}(A)\bigr)+ad\bigl(\Delta^{-1}(A)\bigr)\Bigr)\cdot dar^{-1},$$
where $ad(A)\cdot B=lu(A,B)$.  Then $Dari$ corresponds to
the bracket of derivations, in the sense that
$$Dari(A,B)=Darit(A)\cdot B-Darit(B)\cdot A.$$
\noindent {\bf \S A.6. The $ari$-exponential.}
\vskip .2cm
Recall that a mould with constant term $1$ is said to be bisymmetral, or 
$as*as$, if it is symmetral and its swap is symmetral up to multiplication by 
a constant mould; if the swap is already symmetral the mould is said to 
be {\it strictly bisymmetral} or $as/as$.  Similarly, a mould with 
constant term $0$ is said to be {\it bialternal} or $al*al$ if it is alternal 
and its swap is alternal up to addition of a constant, and strictly alternal 
or $al/al$ if the swap is already alternal.  
\vskip .3cm
\noindent {\bf \S A.7. The $Dgarit$ automorphism.}
For every $Q\in GARI$, there is an automorphism of
$ARI^a_{lu}$ denoted $Dgarit(Q)$, defined simply as the exponential of the
derivation $Darit(P)$, where $Q=exp_{ari}(P)$.  In particular, the action
of this automorphism on the mould $a$ is given by
$$Dgarit(Q)\cdot a=Q-1+a.$$
\vskip .3cm
\'Ecalle showed the following.
\vskip .3cm
\noindent {\bf Theorem 2.2.} {\it Let $P$ be a symmetral mould with constant 
term $1$, even in depth $1$ and whose swap is symmetril up to addition of
a constant mould.
Then $Q=Ad_{gari}(invpal)\cdot P$ is
bisymmetral.  Furthermore, if $P,Q\in GARI$ are both bisymmetral and even
in depth 1, then $garit(P)\cdot Q$ is bisymmetral.}
\vskip .3cm
INTRODUCE $\Delta^*$, $Ad_{gari}$, etc.
\vskip .5cm
\noindent {\bf Proof.}  The Drinfel'd associator $\Phi_{KZ}$ satisfies the group-like double shuffle relations, which translate into moulds as being 
symmetral with swap that is symmetril up to multiplication by a constant mould.
When we pass to $\overline{\Phi}_{KZ}$,
the resulting mould becomes even in depth 1 (which is not the case for 
$\Phi_{KZ}$ due to the term $\zeta(2)[x,y]$).  Thus \'Ecalle's Main Theorem
(\S A.7) shows that
that $\Psi_{KZ}=Ad_{gari}(invpal)\cdot \overline{\Phi}_{KZ}$ is bisymmetral.
It is also well-known that $(\Delta^*)^{-1}g(\tau)$ is bisymmetral
(cf.~[B], [BS]), so again by Proposition 4.2, 
$garit\bigl((\Delta^*)^{-1}g(\tau)\bigr)\cdot \Psi_{KZ}$ is bisymmetral.

The automorphism $garit$ is related to $Dgarit$ by the equality
$$Dgarit(P)\circ\Delta^*=\Delta^*\circ garit\bigl((\Delta^*)^{-1}P\bigr)\eqno(2.2)$$ 
for all $P\in GARI$ (cf.~[S2]).  Thus we have 
$$\bar{E}(\tau)= Dgarit(g(\tau))\cdot\Delta^*\Psi_{KZ}=\Delta^*\Bigl(garit\bigl((\Delta^*)^{-1}g(\tau)\bigr)
\cdot \Psi_{KZ}\Bigr).$$ 
Thus $(\Delta^*)^{-1}\bar{E}_\tau$ is bisymmetral, which concludes the 
proof.\hfill{$\diamondsuit$}
\vskip .3cm
It is shown in
[S2] that for all push-invariant moulds $Q\in ARI$, there is a unique 
derivation $Darit(Q)$ of $ARI^a_{lu}$ such that $Darit(Q)\cdot a=Q$
and $Darit(Q)\cdot ma([a,b])=0$.  We define the $Dari$-exponential as
$$exp_{Dari}(Q)=1+\sum_{n\ge 0} {{1}\over{n!}}Darit(Q)^n\cdot Q.$$
The Darit derivations are related to the Dgarit automorphisms by the formula
$$Dgarit\bigl(exp_{Dari}(Q)\bigr)=exp\bigl(Darit(Q)\bigr).$$
We also recall that $Darit$ is related to the $arit$ derivations studied by
\'Ecalle via the formula
$$Darit(Q)=dar\circ arat\bigl(\Delta^{-1}Q\bigr)\circ dar^{-1},$$
where
$$arat(Q)=-arit(Q)+ad(Q),$$
with $ad(Q)\cdot P=lu(Q,P)$.
\vskip .3cm

We are now armed to attack Lemma ???
whose statement we recall.

\vskip .2cm
\noindent {\bf Lemma 21.} {\it The space $\overline{ARI}_{circneut}$
of circ-neutral moulds $A\in \overline{ARI}$ forms a Lie algebra under
the $\overline{ari}$-bracket.}  

\vskip .2cm
\noindent {\bf Proof.} Let $A,B\in \overline{ARI}_{circneut}$.
We need to show that 
$$\sum_{i=1}^r \overline{ari}(A,B)(v_i,\ldots,v_r,v_1,\ldots,v_{i-1})=0,$$
where the formula for the $\overline{ari}$-bracket on $\overline{ARI}$ is given
as in 4.1.3 by the expression 
$$\eqalign{\overline{ari}(A,B)&=
lu(A,B)+\overline{arit}(B)\cdot A-\overline{arit}(A)\cdot B\cr
&=lu(A,B)+\overline{amit}(B)\cdot A-\overline{anit}(B)\cdot A
-\overline{amit}(A)\cdot B+\overline{anit}(A)\cdot B.
}$$
We will show that this expression is circ-neutral because in fact, each of the five terms in 
the sum is individually circ-neutral.  Let us start by showing this for
the first term, $lu(A,B)$.

Let $\varphi$ denote the cyclic permutation of $\{1,\ldots,r\}$ defined by
$$\varphi(i)=i+1\ \ {\rm for}\ \ 1\le i\le r-1,\ \ 
\varphi(r)=1.$$
By additivity, since the circ-neutrality property is depth-by-depth, 
we may assume that $A$ is concentrated in depth $s$ and $B$ in
depth $t$, with $s\le t$, $s+t=r$.  
In this simplifed situation, we have
$$lu(A,B)(v_1,\ldots,v_r)=A(v_1,\ldots,v_s)B(v_{s+1},\ldots,v_r)
-B(v_1,\ldots,v_t)A(v_{t+1},\ldots,v_r).$$
We have 
$$\eqalign{&\sum_{i=0}^{r-1}lu(A,B)(v_{\varphi^{i}(1)},
\ldots,v_{\varphi^{i}(r)})\cr
&= \sum_{i=0}^{r-1} A(v_{\varphi^i(1)},\ldots,v_{\varphi^i(s)})
B(v_{\varphi^i(s+1)},\ldots,v_{\varphi^i(i)})-
B(v_{\varphi^i(1)},\ldots,v_{\varphi^i(t)})
A(v_{\varphi^i(t+1)},\ldots,v_{\varphi^i(i-1)}),\cr
&= \sum_{i=0}^{r-1} A(v_{\varphi^i(1)},\ldots,v_{\varphi^i(s)})
B(v_{\varphi^i(s+1)},\ldots,v_{\varphi^i(i)})
-A(v_{\varphi^{i+t}(1)},\ldots,v_{\varphi^{i+t}(s)})
B(v_{\varphi^{i+t}(s+1)},\ldots,v_{\varphi^{i+t}(r)})\cr
&=0}$$
as the terms cancel out pairwise.

We now prove that the second term
$$\bigl(\overline{amit}(B)\cdot A\bigr)(v_1,\ldots,v_r)=\sum_{i=1}^s
A(v_1,\ldots,v_{i-1}, v_{i+t},\ldots,v_r)B(v_i-v_{i+t},\ldots,v_{i+t-1}-v_{i+t})$$
is circ-neutral.
%
%
Fix $j\in \{1,\ldots,s\}$ and consider the term
$$A(v_1,\ldots,v_{j-1},v_{j+t},\ldots,v_r)
B(v_j-v_{j+t},\ldots,v_{j+t-1}-v_{j+t}).$$
Thus for each of the other terms 
$$A(v_1,\ldots,v_{i-1},v_{i+t},\ldots,v_r)
B(v_i-v_{i+t},\ldots,v_{i+t-1}-v_{i+t})$$
in the sum, with $i\in \{1,\ldots,s\}$, there is exactly one cyclic 
permutation, namely $\varphi^{j-i}$, that maps this term to
$$A(v_{\varphi^{j-i}(1)},\ldots,v_{\varphi^{j-i}(i-1)},v_{\varphi^{j-i}(i+t},\ldots,v_{\varphi^{j-i}(r)})
B(v_j-v_{j+t},\ldots,v_{j+t-1}-v_{j+t}).$$
For fixed $j\in \{1,\ldots,s\}$, the values of $k=j-i$ mod $s$ as 
$i$ runs through $\{1,\ldots,s\}$ are exactly $\{0,\ldots,s-1\}$.
Therefore, the coefficient of the
term $B(v_j-v_{j+t},\ldots,v_{j+t-1}-v_{j+t})$ in the sum of
the cyclic permutations of $\overline{amit}(B)\cdot A$ is equal to
$$\sum_{k=0}^{s-1} A(v_{\varphi^{k}(1)},\ldots,v_{\varphi^{k}(i-1)},v_{\varphi^{k}(i+t)},\ldots,v_{\varphi^{k}(r)}),$$
which is zero due to the circ-neutrality of $A$.
Thus the coefficient of the term
$B(v_j-v_{j+t},\ldots,v_{j+t-1}-v_{j+t})$ in the sum of the cyclic
permutations of $\overline{amit}(B)\cdot A$ is zero, and this
holds for $1\le j\le s$, so the entire sum is $0$, 
i.e.~$\overline{amit}(B)\cdot A$ is circ-neutral.
The proof of the circ-neutrality of the term $\overline{anit}(B)\cdot 
A$ is analogous. By exchanging $A$ and $B$, this also shows that 
$\overline{amit}(A)\cdot B$ and $\overline{anit}(A)\cdot B$ are circ-neutral, 
which concludes the proof of the lemma.  \hfill{$\square$}
\vskip 1cm
\bye